\documentclass[graybox]{svmult}

\usepackage{mathptmx}       
\usepackage{helvet}         
\usepackage{courier}        
\usepackage{type1cm}        
\usepackage{makeidx}         
\usepackage{graphicx}        
\usepackage{multicol}        
\usepackage[bottom]{footmisc}

\usepackage[numbers, square, comma, sort&compress]{natbib}
\usepackage{amssymb}
\usepackage{mathtools}

\usepackage{bm}

\usepackage{amsmath}
\usepackage{amsfonts}
\usepackage{subfig}
\usepackage{xcolor}

\usepackage[all,cmtip]{xy}
\usepackage{soul}

\newcommand{\diffusiontensor}{\mathbb{K}}

\newcommand{\moperator}[1]{\boldsymbol{\mathsf{#1}}}
\newcommand{\transpose}{\mathsf{T}}

\newcommand \hone {H^{1}(\Omega)}

\renewcommand \u {\bm{u}}
\newcommand \bv {\bm{v}}

\newcommand \f {\mbox{\boldmath$f$}}
\definecolor{not-used}{gray}{0.0}

\newcommand{\VARUN}[1]{\textcolor{black}{#1}}
\makeindex             

\begin{document}

\title*{Mimetic Spectral Element Method for Anisotropic Diffusion}
\author{Marc Gerritsma, Artur Palha, Varun Jain and Yi Zhang}
\institute{Marc Gerritsma \at Faculty of Aerospace Engineering, Delft University of Technology, Kluyverweg 1, Delft, The Netherlands, \email{m.i.gerritsma@tudelft.nl}
\and Artur Palha \at Department of Mechanical Engineering, Eindhoven University of Technology, de Zaale, Eindhoven, The Netherlands, \email{a.palha@tue.nl}
\and Varun Jain \at Faculty of Aerospace Engineering, Delft University of Technology, Kluyverweg 1, Delft, The Netherlands, \email{v.jain@tudelft.nl}
\and Yi Zhang \at Faculty of Aerospace Engineering, Delft University of Technology, Kluyverweg 1, Delft, The Netherlands, \email{y.zhang-14@tudelft.nl}}
%
%
\maketitle

\abstract*{This paper addresses the topological structure of steady, anisotropic, inhomogeneous diffusion problems. Differential operators are represented by sparse incidence matrices, while weighted mass matrices play the role of metric-dependent Hodge matrices. The resulting mixed formulation is point-wise divergence-free if the right hand side function $\f=0$. The method is inf-sup stable; no stabilization is required and the method displays optimal convergence on orthogonal and deformed grids.}

\abstract{This paper addresses the topological structure of steady, anisotropic, inhomogeneous diffusion problems. \VARUN{Two discrete formulations: a) mixed and b) direct formulations are discussed}. Differential operators are represented by sparse incidence matrices, while weighted mass matrices play the role of metric-dependent Hodge matrices. The resulting mixed formulation\VARUN{s} \VARUN{are} point-wise divergence-free if the right hand side function $\f=0$. The method is inf-sup stable and displays optimal convergence on orthogonal and \VARUN{non-affine} grids.}

%


\section{Introduction} \label{sec:introduction}
	Anisotropic and inhomogeneous diffusion appears in many applications such as heat transfer \cite{Incropera2011}, flow through porous media \cite{Neuman1977}, turbulent fluid flow \cite{Taylor1938}, image processing \cite{Perona1990} \VARUN{or} plasma physics \cite{Sovinec2004}. In 2D, steady, anisotropic diffusion is governed by the following elliptic partial differential equation
\begin{equation}
- \nabla \cdot \left ( \mathbb{K} \nabla  p \right ) = f\;,
\label{eq:anisotropic_diffusion_intro}
\end{equation}
\VARUN{Here, $p$ is the flow potential, $f$ the source term,} with $ p=\bar{p}$ along $\Gamma_p$ and $\left ( \mathbb{K}\nabla p, \bm{n} \right ) = \bar{u}_n$ along $\Gamma_u$. Here, for all $\bm{x}$, $\mathbb{K}(\bm{x})$ is a symmetric, positive definite tensor.
	
	In the presence of \emph{strong anisotropy}, i.e. large ratio between the smallest and largest eigenvalues of the diffusion tensor, the construction of robust and efficient discretizations becomes particularly challenging. Under these conditions, the convergence rates of the discretization error can be considerably reduced; this effect is commonly referred in the literature as \emph{locking effect}, see for example \cite{Manzini2007,Babuska1992a,Aavatsmark1998,Aavatsmark1998a}. For sufficiently refined discretizations, the deterioration of the convergence rates eventually disappears. Unfortunately, this may occur only when the grid cell size is prohibitively small.
	
	Another important aspect is mesh flexibility. In many applications of diffusion equations, particularly in porous media flow, typical grids are highly irregular. In many of these situations the results obtained are strongly dependent on the grid type, see \cite{Aziz1993} for a discussion of the use and properties of different grids in reservoir modelling.
	
	\subsection{Overview of standard discretizations}
		In order to overcome these limitations and improve the efficiency and robustness of the discretization of the anisotropic diffusion equations, several approaches have been proposed.
		
		The discretization of the anisotropic diffusion equations in complex media in many situations is still a trade-off between, e.g. \cite{Nilsen2012}:
		\begin{itemize}
			\item Accuracy in the representation of the medium (complex grids).
			\item Accuracy in the discretization of the equations.
		\end{itemize}
		
		The need for such a choice is rooted in the use of numerical schemes based on \emph{two-point flux approximations} (TPFA), see for example, \cite{Wu2009,Aavatsmark2007,Nilsen2012}. These methods produce good approximations on orthogonal grids when the diffusion tensor $\diffusiontensor$ is diagonal, but are known to introduce significant discretization errors in the presence of a non-diagonal diffusion tensor. This introduces severe limitations into the possible grid choices. Under these conditions, the geometric flexibility introduced by \emph{perpendicular bisector} (PEBI) grids, \cite{Heinemann1991,Aziz1993,Palagi1994}, is considerably limited, for example.
		
		It has been known that the discretization error is related to the misalignment between the grid and the principal directions of the diffusion tensor $\diffusiontensor$. In fact, Aavatsmark showed in \cite{Aavatsmark2007} that for TPFA this misalignment leads to the discretization of the wrong diffusion tensor.
		
		These ideas initially led to the construction of grids aligned with the principal axis of the diffusion tensor, so called $\diffusiontensor$\emph{-orthogonal grids}, see for example \cite{Heinemann1991,Gunasekera1997}. This approach significantly improves the performance of the numerical method but substantially limits the geometric flexibility.
		
		More recently, \emph{multipoint flux-approximation} (MPFA) schemes have been introduced specifically to address these limitations, see e.g. the initial works by Aavatsmark \cite{Aavatsmark1998,Aavatsmark1998a} or a more recent presentation \cite{Aavatsmark2002}, and by Edwards and Rogers \cite{Edwards1998}. This method is based on a cell-centred finite volume formulation and introduces a dual grid in order to generate shared sub-cells and sub-interfaces. This in turn produces a discretization of the flux between two cells that involves a linear combination of several adjacent cells. This method is robust and locally conservative but does not guarantee a resulting symmetric discrete diffusion operator. More recently, this work has been connected to the mixed finite element method, \cite{Edwards2002}.
		
		Alternative approaches based on the finite element formulation have also been proposed by several authors. We briefly mention the work on the control-volume finite element discretization by Forsyth \cite{Forsyth1990} and Durlofsky \cite{Durlofsky1994}, on nodal Galerkin finite elements by Young \cite{Young1999}, and on mixed finite elements by \VARUN{Durlofsky} \cite{Durlofsky1993}.
	
	\subsection{Overview of mimetic discretizations}
		Over the years, the development of numerical schemes that preserve some of the structures of the differential models they approximate has been identified as an important ingredient of numerical analysis. One of the contributions of the formalism of mimetic methods is to identify differential geometry as the proper language in which to encode these structures/symmetries. Another novel aspect of mimetic discretizations is the identification and separation of physical field laws into two sets: (i) topological relations (metric-free), and (ii) constitutive relations (metric dependent). Topological relations are intimately related to conservation laws and can (and should) be exactly represented on the computational grid. Constitutive relations include all material properties and therefore are approximate relations. For this reason, all numerical discretization error should be included in these equations. A general introduction and overview of spatial and temporal mimetic/geometric methods can be found in \cite{Christiansen2011,perot43discrete,Budd2003,Hairer2006}.

		The relation between differential geometry and algebraic topology in physical theories was first established by Tonti \cite{tonti1975formal}. Around the same time Dodziuk \cite{Dodziuk76} set up a finite difference framework for harmonic functions based on Hodge theory. Both Tonti and Dodziuk introduce differential forms and cochain spaces as the building blocks for their theory. The relation between differential forms and cochains is established by the Whitney map ($k$-cochains $\rightarrow$ $k$-forms) and the de Rham map ($k$-forms $\rightarrow$ $k$-cochains). The interpolation of cochains to differential forms on a triangular grid was already established by Whitney, \cite{Whitney57}. These generalized interpolatory forms are now known as {\em Whitney forms}.

		Hyman and Scovel \cite {HymanScovel90} set up the discrete framework in terms of cochains, which are the natural building blocks of finite volume methods. Later\VARUN{,} Bochev and Hyman \cite{bochev2006principles} extended this work and derived discrete operators such as the discrete wedge product, the discrete codifferential, \VARUN{and} the discrete inner products.

		Robidoux, Hyman, Steinberg and Shashkov, \cite{HymanShashkovSteinberg97,HymanShashkovSteinberg2002,HYmanSteinberg2004,RobidouxAdjointGradients1996,RobidouxThesis,bookShashkov,Steinberg1996,SteibergZingano2009, Kikinzon2017} used symmetry considerations to construct discretizations on rough grids, within the finite difference/volume setting . In a more recent paper by Robidoux and Steinberg \cite{Robidoux2011} a finite difference discrete vector calculus is presented. In that work, the differential operators grad, curl and div are exactly represented at the discrete level and the numerical approximations are all contained in the constitutive relations, which are already polluted by modeling and experimental error. For mimetic finite differences, see also the work of Brezzi et al. \cite{BrezziBuffaLipnikov2009,brezzi2010} and Beir\~{a}o da Veiga et al. \cite{Veiga_mimetic_finite_differences}.

		The application of mimetic ideas to unstructured triangular staggered grids has been extensively studied by Perot, \cite{Perot2000,ZhangSchmidtPerot2002,perot2006mimetic,PerotSubramanian2007a,PerotSubramanian2007}, specially in \cite{perot43discrete} where the rationale of preserving symmetries in numerical algorithms is well described. The most \emph{geometric approach} is presented in the work by Desbrun et al. \cite{desbrun2005discrete,ElcottTongetal2007,MullenCraneetal2009,PavlovMullenetal2010} and the thesis by Hirani \cite{Hirani_phd_2003}.
		
		The \emph{Japanese papers}  by Bossavit, \cite{bossavit_japan_computational_1,bossavit_japan_computational_2,bossavit_japan_computational_3,bossavit_japan_computational_4,bossavit_japan_computational_5}, serve as an excellent introduction and motivation for the use of differential forms in the description of physics and the use in numerical modeling. The field of application is electromagnetism, but these papers are sufficiently general to extend to other physical theories.
		
		In a series of papers by Arnold, Falk and Winther, \cite{arnold:Quads,arnold2006finite,arnold2010finite}, a finite element exterior calculus framework is developed. Higher order methods are described by Rapetti \cite{Rapetti2007,Rapetti2009} and Hiptmair \cite{hiptmair2001}. Possible extensions to spectral methods were described by Robidoux, \cite{robidoux-polynomial}. A different approach for constructing arbitrary order mimetic finite elements has been proposed by the authors \cite{Palha2014, gerritsmaicosahom2012,bouman::icosahom2009,palha::icosahom2009}, with applications to advection problems \cite{palhaAdvectionIcosahom2014}, Stokes' flow \cite{kreeft::stokes}, MHD equilibrium \cite{Palha2016}, Navier-Stokes \cite{Palha2017}, and within a Least-Squares finite element formulation \cite{Gerritsma,palha:lssc2009,bochev_mimetic_ls_2014,gerritsma_mimetic_ls_stokes_2016}.
		
		Extensions of these ideas to polyhedral meshes have been proposed by Ern, Bonelle and co-authors in \cite{Bonelle2015,Bonelle2014,Bonelle2015a,Bonelle2016}, by di Pietro and co-authors in \cite{DiPietro2014,DiPietro2015}, by Brezzi and co-authors in \cite{Brezzi2014}, and by Beir\~{a}o da Veiga and co-authors in \cite{BeiraodaVeiga2014,BeiraodaVeiga2016,DaVeiga2015,DaVeiga2015a}. These approaches provide more geometrical flexibility while maintaining fundamental structure preserving properties.
		
		Mimetic isogeometric discretizations have been introduced by Buffa et al. \cite{BuffaDeFalcoSangalli2011}, Evans and Hughes \cite{Evans2013}, and Hiemstra et al. \cite{Hiemstra2014}.
	
		Another approach \VARUN{to develop} a discretization of the physical field laws \VARUN{is} based on a discrete variational principle for the discrete Lagrangian action. This approach has been used in the past to construct variational integrators for Lagrangian systems, e.g. \cite{Kouranbaeva2000,Marsden2003}. Kraus and Maj \cite{Kraus2015} have used the method of formal Lagrangians to derive generalized Lagrangians for non-Lagrangian systems of equations. This allows to apply variational techniques to construct structure preserving discretizations on a much wider range of systems.
\VARUN{Recently, Bauer and Gay-Balmaz presented variational integrators for elastic and pseudo-incompressible flows \cite{Bauer2017}.}

		Due to the inherent challenges in discretizing the diffusion equations with anisotropic diffusion tensor $\diffusiontensor$, several authors have explored different mimetic discretizations of these equations. Focussing on generalized diffusion equations we highlight \cite{Hermeline2000,Bastian2011a,HymanShashkovSteinberg97,HymanShashkovSteinberg2002,HYmanSteinberg2004,Kikinzon2017,RobidouxAdjointGradients1996,RobidouxThesis,bookShashkov,Steinberg1996,SteibergZingano2009,PerotSubramanian2007a} for a finite-difference/finite-volume setting, \cite{Brezzi2006,Brezzi2007,Brezzi2005,Bonelle2015} for polyhedral discretizations, and \cite{BoffiGastaldi2009,bochev:RehabQuad,Rebelo2014,Palha2014,Younes2010} for a finite element/mixed finite element setting. For applications to Darcy flow equations and reservoir modelling see for example \cite{Alpak2010a,Alpak2007,Lie2012,Nilsen2012,Aarnes2008,Hirani2015,Dziubek2016}.
	
	\subsection{Outline of chapter}
In Section~\ref{sec:theory} the topological structure of anisotropic diffusion problems is discussed. In Section~\ref{sec:MSEM} spectral basis functions are introduced which are compatible with the topological structure introduced in Section~\ref{sec:theory}. In Section~\ref{sec:transformation} transformation to curvilinear elements is discussed. Results of the proposed method are presented in Section~\ref{sec:numerical_results}.


\section{Anisotropic diffusion / Darcy problem}\label{sec:theory}
Let $\Omega \subset \mathbb{R}^d$ be a contractible domain with Lipschitz continuous boundary $\partial \Omega = \Gamma_p \cup \Gamma_u$, $ \Gamma_p \cap \Gamma_u = \varnothing$. The steady anisotropic diffusion problem is given by
\begin{equation}
- \nabla \cdot \left ( \mathbb{K} \nabla  p \right ) = f\;,
\label{eq:anisotropic_diffusion}
\end{equation}
with $ p=\bar{p}$ along $\Gamma_p$ and $\left ( \mathbb{K}\nabla p, \bm{n} \right ) = \bar{u}_n$ along $\Gamma_u$. Here, for all $\bm{x}$, $\mathbb{K}(\bm{x})$ is a symmetric, positive definite tensor, i.e. there exist constants $\alpha,C>0$ such that
\[ \alpha \bm{\xi}^T \bm{\xi} \leq \bm{\xi}^T \mathbb{K}(\bm{x}) \bm{\xi} \leq C \bm{\xi}^T \bm{\xi} \;.\]
If $\Gamma_p \neq \varnothing$, then (\ref{eq:anisotropic_diffusion}) has a unique solution. If $\Gamma_p = \varnothing$ then (\ref{eq:anisotropic_diffusion}) only possesses solutions if
\[ \int_{\partial \Omega} \bar{u}_n \,\mathrm{d}S = \int_{\Omega} f \, \mathrm{d} \Omega \;,
\]
in which case the solution, $p$, is determined up to a constant.

An equivalent first order system is obtained by introducing $\bm{u} = -\mathbb{K} \nabla p$ in which case (\ref{eq:anisotropic_diffusion}) can be written as
\begin{equation}
\left \{ \begin{array}{ll}
\bm{u} + \mathbb{K} \nabla p = 0 \quad & \mbox{in } \Omega \\[1ex]
\nabla \cdot \bm{u} = f \quad  &\mbox{in } \Omega
\end{array} \right .
\quad \mbox{with} \quad
\left \{
\begin{array}{ll}
(\bm{u},\bm{n}) = \bar{u}_n \quad & \mbox{along } \Gamma_u \\[1ex]
p = \bar{p} \quad & \mbox{along } \Gamma_p
\end{array} \right . \;.
\label{eq:Darcy_problem}
\end{equation}

An alternative first-order formulation is given by
\begin{equation}
\left \{ \begin{array}{ll}
\bm{v} -  \nabla p = 0 \quad & \mbox{in } \Omega \\[1ex]
\bm{u} + \mathbb{K} \bm{v} = 0 \quad & \mbox{in } \Omega \\[1ex]
\nabla \cdot \bm{u} = f \quad  &\mbox{in } \Omega
\end{array} \right .
\quad \mbox{with} \quad
\left \{
\begin{array}{ll}
(\bm{u},\bm{n}) = \bar{u}_n \quad & \mbox{along } \Gamma_u \\[1ex]
p = \bar{p} \quad & \mbox{along } \Gamma_p
\end{array} \right . \;.
\label{eq:extended_Darcy_problem}
\end{equation}

Formulation (\ref{eq:Darcy_problem}) is generally referred to as the {\em Darcy problem}, while the relation $\bm{u} = -\mathbb{K} \nabla p$ is called {\em Darcy's law}, \cite{Neuman1977}. The Darcy problem plays an important role in reservoir engineering. In this case $\bm{u}$ is the flow velocity in a porous medium and $p$ denotes the pressure.

\VARUN{While the formulations (\ref{eq:anisotropic_diffusion}), (\ref{eq:Darcy_problem}) and (\ref{eq:extended_Darcy_problem}) are equivalent}, (\ref{eq:anisotropic_diffusion})  only has $1$ unknown, $p$, (\ref{eq:Darcy_problem}) has $(d+1)$ unknowns, $p$ and the $d$ components of $\bm{u}$\VARUN{,} and (\ref{eq:extended_Darcy_problem}) has $(2d+1)$ unknowns. \VARUN{Formulation} (\ref{eq:extended_Darcy_problem}) is of special interest, because it decomposes the anisotropic diffusion problem into two topological conservation laws and one constitutive law\footnote{
An even more extended system is, see for instance \cite{bochev_mimetic_ls_2014}
\begin{equation}
\left \{ \begin{array}{ll}
\bm{v} -  \nabla p = 0 \quad & \mbox{in } \Omega \\[1ex]
\bm{u} + \mathbb{K} \bm{v} = 0 \quad & \mbox{in } \Omega \\[1ex]
\VARUN{\nabla \cdot \bm{u} - \psi} = 0 \quad  &\mbox{in } \Omega \\[1ex]
\VARUN{\psi = f} \quad  &\mbox{in } \Omega
\end{array} \right .
\quad \mbox{with} \quad
\left \{
\begin{array}{ll}
(\bm{u},\bm{n}) = \bar{u}_n \quad & \mbox{along } \Gamma_u \\[1ex]
p = \bar{p} \quad & \mbox{along } \Gamma_p
\end{array} \right . \;.
\label{eq:very_extended_Darcy_problem}
\end{equation}
This seems an unnecessarily complicated system. If we eliminate $\psi$ from (\ref{eq:very_extended_Darcy_problem}) we obtain (\ref{eq:extended_Darcy_problem}). The usefulness of this system lies in the fact that by introducing $\psi$, the conservation $\nabla \cdot \bm{u}=f$ becomes independent of the data of the PDE, in this case the right hand side function. A similar situation occurs when $\mathbb{K}=\mathbb{I}$, the identity tensor, then the equation $\bm{u}\VARUN{+}\mathbb{K}\bm{v}=0$ in (\ref{eq:extended_Darcy_problem}) seems redundant, but we have good reason \VARUN{to} keep this seemingly redundant equation as we will show in this paper.
}.
By making a suitable choice {\em where} and {\em how} to represent the unknowns on a grid, the topological relations, $\bm{v} -  \nabla p = 0$ and $\nabla \cdot \bm{u} = f$ reduce to extremely simple algebraic relations which \VARUN{depend only on the topology of the mesh and are} independent of the mesh size, independent of the shape of the mesh\VARUN{,} and independent of the order of the numerical scheme. We will refer to such discretizations as {\em exact discrete representations}.

\subsection{Gradient relation}\label{sec:gradient}
Consider two points $A,B\in \Omega$ and a curve $\mathcal{C}$ which connects these two points, then
\[
\bm{v} -  \nabla p = 0   \quad \Longrightarrow \quad \bar{\bm{v}}_{\mathcal{C}} := \int_{\mathcal{C}} \bm{v} \cdot \mathrm{d}l = \int_A^B \bm{v} \cdot \mathrm{d}l = \int_A^B \nabla p \cdot \mathrm{d}l = p(B) - p(A) \;,
\]
where $\mathrm{d}l$ is a small increment along the curve $\mathcal{C}$.

Suppose that we take another curve $\tilde{\mathcal{C}}$ which connects the two points $A$ and $B$ then we also have
\begin{equation} \bar{\bm{v}}_{\tilde{\mathcal{C}}} \VARUN{:}= \int_{\tilde{\mathcal{C}}} \bm{v} \cdot \mathrm{d}l  = p(B) - p(A) \;,
\label{eq:velocity_gradient_relation}
\end{equation}
\begin{figure}[h!]
	\centering
	\includegraphics[width=0.8\linewidth]{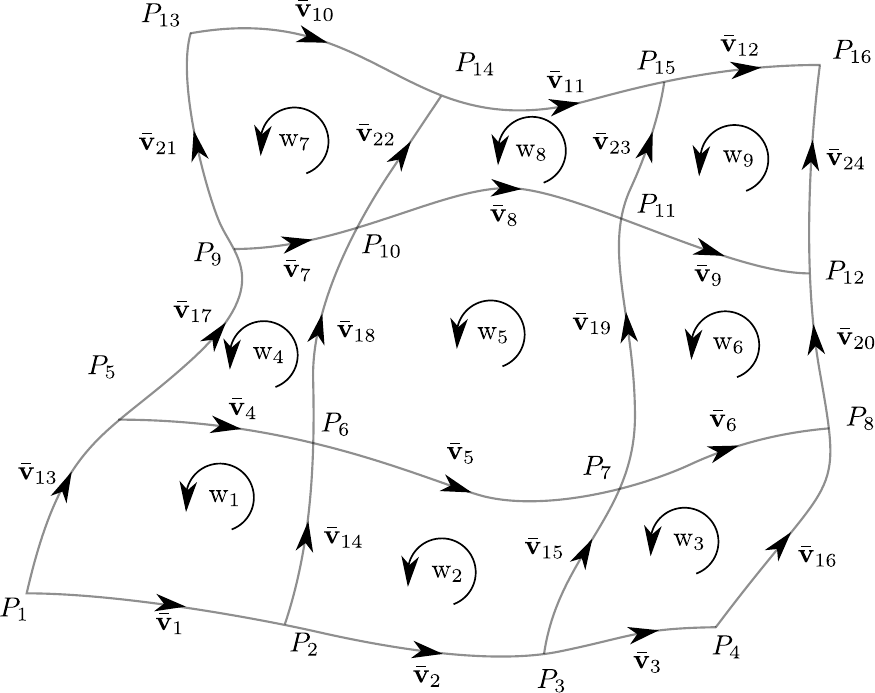}
	\caption{Relation between pressure in points, integrated velocity along line segments and vorticity in surfaces.}
	\label{fig:inner_oriented_p_v}
\end{figure}

The integral along $\mathcal{C}$ is equal to the integral along $\tilde{\mathcal{C}}$. We will refer to $\bar{\bm{v}}$ as an {\em integral value}, since it denotes an integral and not a point-wise evaluation of $\bm{v}$. The advantages of integral values are:
\begin{enumerate}
\item The velocity-gradient relation is exact. It is not obtained by \VARUN{truncated} Taylor-series expansions or does not depend on the choice of \VARUN{basis functions/interpolations}.
\item Does not depend on mesh parameters. The mesh size $h$ does not appear in (\ref{eq:velocity_gradient_relation}). Whether the curve which connects two points is straight or curved is irrelevant in this relation, therefore this relation is directly applicable on curved domains.
\item Integral quantities are additive.
\end{enumerate}

%
%

Consider the points and lines segments as shown in Figure~\ref{fig:inner_oriented_p_v}. In this figure the arrow along the curves indicates the direction in which $\bm{v}$ is integrated\footnote{The points in the grid shown in Figure~\ref{fig:inner_oriented_p_v} are also `oriented', in the sense that when we `move into a point following the integration direction' we assign a positive value and when we `leave a point' we assign a negative value. That is why we have {\em plus} $P(B)$ and {\em minus} $P(A)$ in (\ref{eq:velocity_gradient_relation}). This is just a convention. Without loss of generality we could change this sign convention.}. Application of (\ref{eq:velocity_gradient_relation}) shows, for instance, that
\[ \bar{\bm{v}}_{14} = P_6 - P_2 \;.\]
The additivity property implies that
\begin{align*} P_7 - P_2 & = \bar{\bm{v}}_2 + \bar{\bm{v}}_{15} = P_3 - P_2 + P_7 - P_3 \\
 & = \bar{\bm{v}}_{14} + \bar{\bm{v}}_5 = P_6 - P_2 + P_7 - P_6\;,
\end{align*}
and even more paths can be constructed that connect $P_2$ and $P_7$. The independence of the path depends critically on the assumption that the space is contractible, i.e. there are no holes in the domain (Poincar\'{e}'s Lemma).

A special case is the curve from a point to itself, say $P_2 \rightarrow P_2$ in Figure~\ref{fig:inner_oriented_p_v}. This integral is zero and if the integral is independent of the path this implies that, for instance,
\begin{equation}
0=\bar{\bm{v}}_2 + \bar{\bm{v}}_{15} - \bar{\bm{v}}_5 - \bar{\bm{v}}_{14} = \oint \bm{v} \cdot \mathrm{d}l = \iint \nabla \times \bm{v} \cdot \mathrm{d}\bm{S} = \VARUN{\mathsf{w}}_2 \;,
\label{eq:vorticity_relation}
\end{equation}
where we once again use the additivity property. We see that the circulation vanishes if $\bm{v}$ is a potential flow, which in turn implies that the \VARUN{circulation} of the velocity field over \VARUN{the boundary of} any surface vanishes. Or, using Stokes' theorem, the integrated vorticity \VARUN{$\mathsf{w}$} vanishes. \VARUN{Here the vorticity $\mathsf{w}$ is represented as the integral over a surface}.

We can collect all the integrated velocity fields and pressures in Figure~\ref{fig:inner_oriented_p_v} in the following form
{\scriptsize
\[
\left (
\begin{array}{c}
\bar{\bm{v}}_{1} \\
\bar{\bm{v}}_{2} \\
\bar{\bm{v}}_{3} \\
\bar{\bm{v}}_{4} \\
\bar{\bm{v}}_{5} \\
\bar{\bm{v}}_{6} \\
\bar{\bm{v}}_{7} \\
\bar{\bm{v}}_{8} \\
\bar{\bm{v}}_{9} \\
\bar{\bm{v}}_{10} \\
\bar{\bm{v}}_{11} \\
\bar{\bm{v}}_{12} \\
\bar{\bm{v}}_{13} \\
\bar{\bm{v}}_{14} \\
\bar{\bm{v}}_{15} \\
\bar{\bm{v}}_{16} \\
\bar{\bm{v}}_{17} \\
\bar{\bm{v}}_{18} \\
\bar{\bm{v}}_{19} \\
\bar{\bm{v}}_{20} \\
\bar{\bm{v}}_{21} \\
\bar{\bm{v}}_{22} \\
\bar{\bm{v}}_{23} \\
\bar{\bm{v}}_{24}
\end{array} \right ) =
\left ( \begin{array}{cccccccccccccccc}
-1 & 1 & 0 & 0 & 0 & 0 & 0 & 0 & 0 & 0 & 0 & 0 & 0 & 0 & 0 & 0 \\
0 & -1 & 1 & 0 & 0 & 0 & 0 & 0 & 0 & 0 & 0 & 0 & 0 & 0 & 0 & 0 \\
0 & 0 & -1 & 1 & 0 & 0 & 0 & 0 & 0 & 0 & 0 & 0 & 0 & 0 & 0 & 0 \\
0 & 0 & 0 & 0 & -1 & 1 & 0 & 0 & 0 & 0 & 0 & 0 & 0 & 0 & 0 & 0 \\
0 & 0 & 0 & 0 & 0 & -1 & 1 & 0 & 0 & 0 & 0 & 0 & 0 & 0 & 0 & 0 \\
0 & 0 & 0 & 0 & 0 & 0 & -1 & 1 & 0 & 0 & 0 & 0 & 0 & 0 & 0 & 0 \\
0 & 0 & 0 & 0 & 0 & 0 & 0 & 0 & -1 & 1 & 0 & 0 & 0 & 0 & 0 & 0 \\
0 & 0 & 0 & 0 & 0 & 0 & 0 & 0 & 0 & -1 & 1 & 0 & 0 & 0 & 0 & 0 \\
0 & 0 & 0 & 0 & 0 & 0 & 0 & 0 & 0 & 0 & -1 & 1 & 0 & 0 & 0 & 0 \\
0 & 0 & 0 & 0 & 0 & 0 & 0 & 0 & 0 & 0 & 0 & 0 & -1 & 1 & 0 & 0 \\
0 & 0 & 0 & 0 & 0 & 0 & 0 & 0 & 0 & 0 & 0 & 0 & 0 & -1 & 1 & 0 \\
0 & 0 & 0 & 0 & 0 & 0 & 0 & 0 & 0 & 0 & 0 & 0 & 0 & 0 & -1 & 1 \\
-1 & 0 & 0 & 0 & 1 & 0 & 0 & 0 & 0 & 0 & 0 & 0 & 0 & 0 & 0 & 0 \\
0 & -1 & 0 & 0 & 0 & 1 & 0 & 0 & 0 & 0 & 0 & 0 & 0 & 0 & 0 & 0 \\
0 & 0 & -1 & 0 & 0 & 0 & 1 & 0 & 0 & 0 & 0 & 0 & 0 & 0 & 0 & 0 \\
0 & 0 & 0 & -1 & 0 & 0 & 0 & 1 & 0 & 0 & 0 & 0 & 0 & 0 & 0 & 0 \\
0 & 0 & 0 & 0 & -1 & 0 & 0 & 0 & 1 & 0 & 0 & 0 & 0 & 0 & 0 & 0 \\
0 & 0 & 0 & 0 & 0 & -1 & 0 & 0 & 0 & 1 & 0 & 0 & 0 & 0 & 0 & 0 \\
0 & 0 & 0 & 0 & 0 & 0 & -1 & 0 & 0 & 0 & 1 & 0 & 0 & 0 & 0 & 0 \\
0 & 0 & 0 & 0 & 0 & 0 & 0 & -1 & 0 & 0 & 0 & 1 & 0 & 0 & 0 & 0 \\
0 & 0 & 0 & 0 & 0 & 0 & 0 & 0 & -1 & 0 & 0 & 0 & 1 & 0 & 0 & 0 \\
0 & 0 & 0 & 0 & 0 & 0 & 0 & 0 & 0 & -1 & 0 & 0 & 0 & 1 & 0 & 0 \\
0 & 0 & 0 & 0 & 0 & 0 & 0 & 0 & 0 & 0 & -1 & 0 & 0 & 0 & 1 & 0 \\
0 & 0 & 0 & 0 & 0 & 0 & 0 & 0 & 0 & 0 & 0 & -1 & 0 & 0 & 0 & 1
\end{array}
\right)
\left ( \begin{array}{c}
P_1 \\
P_2 \\
P_3 \\
P_4 \\
P_5 \\
P_6 \\
P_7 \\
P_8 \\
P_9 \\
P_{10} \\
P_{11} \\
P_{12} \\
P_{13} \\
P_{14} \\
P_{15} \\
P_{16}
\end{array} \right )\VARUN{\;.}
\]
}
If we store all $\bar{\bm{v}}_i$ in a vector $\mathsf{v}$ and all $P_j$ in a vector $\mathsf{P}$ and denote the matrix by $\mathbb{E}^{1,0}$, we have
\[ \mathsf{v} = \mathbb{E}^{1,0} \mathsf{P} \;.\]
\VARUN{If we now also collect all the integrated vorticities, $\mathsf{w} _i$, we can relate them to  the integrated velocities in the following way}
{\scriptsize
\[
\left ( \begin{array}{c}
\VARUN{\mathsf{w}}_1 \\
\VARUN{\mathsf{w}}_2 \\
\VARUN{\mathsf{w}}_3 \\
\VARUN{\mathsf{w}}_4 \\
\VARUN{\mathsf{w}}_5 \\
\VARUN{\mathsf{w}}_6 \\
\VARUN{\mathsf{w}}_7 \\
\VARUN{\mathsf{w}}_8 \\
\VARUN{\mathsf{w}}_9
\end{array} \right ) =
\left ( \begin{array}{cccccccccccccccccccccccc}
1 & 0 & 0 & -1 & 0 & 0 & 0 & 0 & 0 & 0 & 0 & 0 & -1 & 1 & 0 & 0 & 0 & 0 & 0 & 0 & 0 & 0 & 0 & 0 \\
0 & 1 & 0 & 0 & -1 & 0 & 0 & 0 & 0 & 0 & 0 & 0 & 0 & -1 & 1 & 0 & 0 & 0 & 0 & 0 & 0 & 0 & 0 & 0 \\
0 & 0 & 1 & 0 & 0 & -1 & 0 & 0 & 0 & 0 & 0 & 0 & 0 & 0 & -1 & 1 & 0 & 0 & 0 & 0 & 0 & 0 & 0 & 0 \\
0 & 0 & 0 & 1 & 0 & 0 & -1 & 0 & 0 & 0 & 0 & 0 & 0 & 0 & 0 & 0 & -1 & 1 & 0 & 0 & 0 & 0 & 0 & 0\\
0 & 0 & 0 & 0 & 1 & 0 & 0 & -1 & 0 & 0 & 0 & 0 & 0 & 0 & 0 & 0 & 0 & -1 & 1 & 0 & 0 & 0 & 0 & 0\\
0 & 0 & 0 & 0 & 0 & 1 & 0 & 0 & -1 & 0 & 0 & 0 & 0 & 0 & 0 & 0 & 0 & 0 & -1 & 1 & 0 & 0 & 0 & 0\\
0 & 0 & 0 & 0 & 0 & 0 & 1 & 0 & 0 & -1 & 0 & 0 & 0 & 0 & 0 & 0 & 0 & 0 & 0 & 0 & -1 & 1 & 0 & 0\\
0 & 0 & 0 & 0 & 0 & 0 & 0 & 1 & 0 & 0 & -1 & 0 & 0 & 0 & 0 & 0 & 0 & 0 & 0 & 0 & 0 & -1 & 1 & 0\\
0 & 0 & 0 & 0 & 0 & 0 & 0 & 0 & 1 & 0 & 0 & -1 & 0 & 0 & 0 & 0 & 0 & 0 & 0 & 0 & 0 & 0 & -1 & 1
\end{array} \right )
\left (
\begin{array}{c}
\bar{\bm{v}}_{1} \\
\bar{\bm{v}}_{2} \\
\bar{\bm{v}}_{3} \\
\bar{\bm{v}}_{4} \\
\bar{\bm{v}}_{5} \\
\bar{\bm{v}}_{6} \\
\bar{\bm{v}}_{7} \\
\bar{\bm{v}}_{8} \\
\bar{\bm{v}}_{9} \\
\bar{\bm{v}}_{10} \\
\bar{\bm{v}}_{11} \\
\bar{\bm{v}}_{12} \\
\bar{\bm{v}}_{13} \\
\bar{\bm{v}}_{14} \\
\bar{\bm{v}}_{15} \\
\bar{\bm{v}}_{16} \\
\bar{\bm{v}}_{17} \\
\bar{\bm{v}}_{18} \\
\bar{\bm{v}}_{19} \\
\bar{\bm{v}}_{20} \\
\bar{\bm{v}}_{21} \\
\bar{\bm{v}}_{22} \\
\bar{\bm{v}}_{23} \\
\bar{\bm{v}}_{24}
\end{array} \right )\VARUN{\;.}
\]
}
If we store all vorticity integrals, $\VARUN{\mathsf{w}}_i$ in the vector $\vec{\VARUN{\mathsf{w}}}$, then we can write this as
\begin{equation}
\mathsf{w} = \mathbb{E}^{2,1} \mathsf{v} \;.
\label{eq:topological_curl_v}
\end{equation}
The matrices $\mathbb{E}^{1,0}$ and $\mathbb{E}^{2,1}$ are called {\em incidence matrices}. We have $\mathbb{E}^{2,1} \cdot \mathbb{E}^{1,0} \equiv 0$. This identity holds for this particular case, but is generally true; it holds when we would have used triangles or polyhedra instead of quadrilaterals and it holds in any space dimension $d$. If $\mathbb{E}^{1,0}$ represents the gradient operation and $\mathbb{E}^{2,1}$ represents the curl operation, then $\mathbb{E}^{2,1} \cdot \mathbb{E}^{1,0} \equiv 0$ is the discrete analogue of the vector identity $\nabla \times \nabla \equiv 0$, \cite{Bonelle2014,Bonelle2015a,bossavit_japan_computational_1,bossavit_japan_computational_2,desbrun2005discrete,Nicolaides,Robidoux2011}.

If boundary conditions for $p$ are prescribed along $\partial \Omega$, then these degrees of freedom can be removed from the grid in Figure~\ref{fig:inner_oriented_p_v}.
\begin{figure}[h!]
	\centering
	\includegraphics[width=0.6\linewidth]{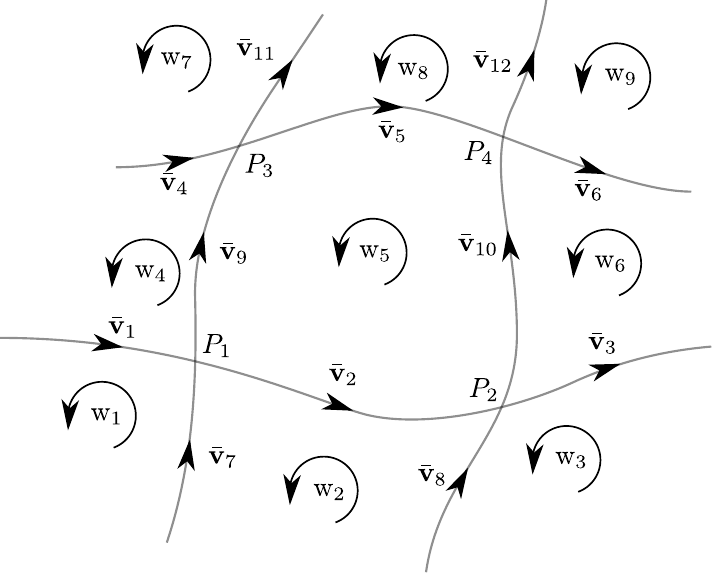}
	\caption{Relation between pressure in points and integrated velocity along line segments in case $\Gamma_p = \partial \Omega$.}
	\label{fig:inner_oriented_p_v_bc}
\end{figure}

If $p$ is known along the boundary then \VARUN{the integral of} $\bm{v}$ is also known along the boundary, so the degrees of freedom for $\mathsf{v}$ can also be removed. Relabeling the remaining unknowns gives the geometric degrees of freedom as shown in Figure~\ref{fig:inner_oriented_p_v_bc}.

{\scriptsize
\begin{equation}
\left (
\begin{array}{c}
\bar{\bm{v}}_{1} \\
\bar{\bm{v}}_{2} \\
\bar{\bm{v}}_{3} \\
\bar{\bm{v}}_{4} \\
\bar{\bm{v}}_{5} \\
\bar{\bm{v}}_{6} \\
\bar{\bm{v}}_{7} \\
\bar{\bm{v}}_{8} \\
\bar{\bm{v}}_{9} \\
\bar{\bm{v}}_{10} \\
\bar{\bm{v}}_{11} \\
\bar{\bm{v}}_{12}
\end{array} \right ) =
\left ( \begin{array}{cccc}
1 & 0 & 0 & 0 \\
-1 & 1 & 0 & 0 \\
0 & -1 & 0 & 0 \\
0 & 0 & 1 & 0 \\
0 & 0 & -1 & 1 \\
0 & 0 & 0 & -1 \\
1 & 0 & 0 & 0 \\
0 & 1 & 0 & 0 \\
-1 & 0 & 1 & 0 \\
0 & -1 & 0 & 1 \\
0 & 0 & -1 & 0 \\
0 & 0 & 0 & -1
\end{array}
\right)
\left ( \begin{array}{c}
P_1 \\
P_2 \\
P_3 \\
P_4
\end{array} \right )\VARUN{\;.}
\label{eq:velocity_pressure_with_p_prescribed}
\end{equation}
}
\subsection{Divergence relation}\label{sec:divergence}
Consider a bounded, contractible volume $\mathcal{V} \subset \Omega$ then we have
\[
\nabla \cdot \bm{u} = f \quad \Longrightarrow \quad \int_{\partial \mathcal{V}} \bm{u}\cdot \bm{n}\,\mathrm{d}S = \int_{\mathcal{V}} f \,\mathrm{d}\mathcal{V} \;.
\]
\begin{figure}[h!]
	\centering
	\includegraphics[width=0.8\linewidth]{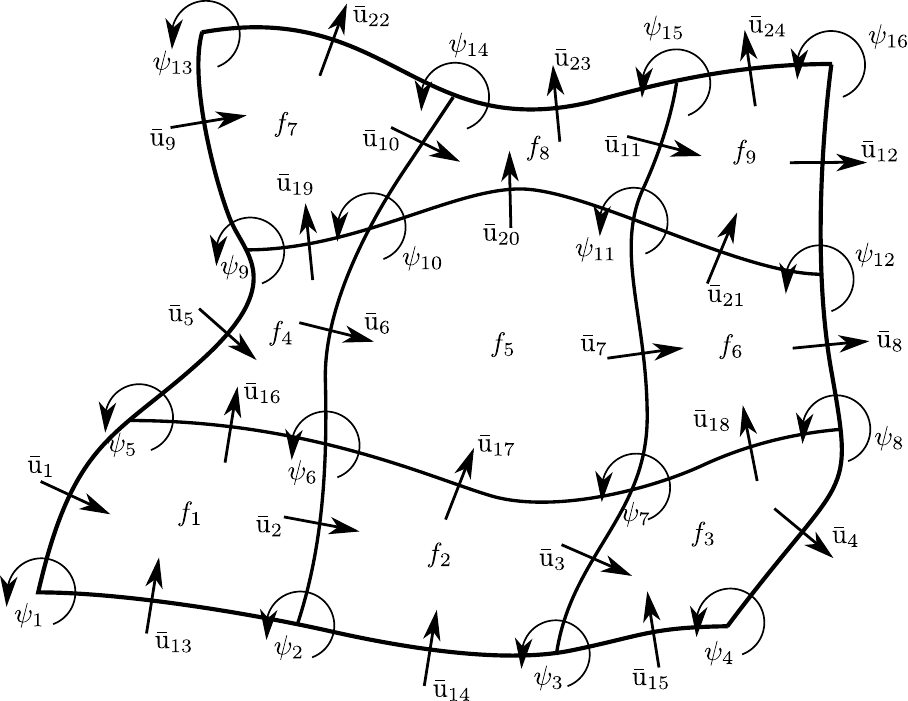}
	\caption{Stream function, fluxes and the divergence degrees of freedom.}
	\label{fig:outer_oriented_u}
\end{figure}

If the boundary $\partial \mathcal{V}$ can be partitioned into $n$ sub-boundaries, \VARUN{$\partial \mathcal{V} = \bigcup\limits_{i} \Gamma _i$ and $\bigcap\limits_{i} \Gamma _i =0$}, we have
\[ \pm \sum_{i=1}^n \bar{\u}_i = \pm \sum_{i=1}^n \int_{\VARUN{\Gamma}_i} \bm{u}\cdot \bm{n}\,\mathrm{d}S = \int_{\mathcal{V}} f \,\mathrm{d}\mathcal{V} =: f_\mathcal{V} \;,
\]
where we have the convention that the fluxes, $ \bar{\u}_i$, are {\em positive} when the flow leaves the volume and {\em negative} when the flow enters the volume.
For a 2D case the integral flux degrees of freedom, $\bar{\u}_i$ are depicted in Figure~\ref{fig:outer_oriented_u}.
The arrow in this figure indicates the positive default direction of the fluxes.
The integrated values of source function \VARUN{$f$} are shown in the 2D volumes in Figure~\ref{fig:outer_oriented_u} as \VARUN{$f_i$}. The topological relation between the fluxes and the integrated \VARUN{source values} $f_i$, for the situation shown in Figure~\ref{fig:outer_oriented_u}, is given by

{
	\scriptsize
\[
\left ( \begin{array}{cccccccccccccccccccccccc}
-1 & 1 & 0 & 0 & 0 & 0 & 0 & 0 & 0 & 0 & 0 & 0 & -1 & 0 & 0 & 1 & 0 & 0 & 0 & 0 & 0 & 0 & 0 & 0 \\
0 & -1 & 1 & 0 & 0 & 0 & 0 & 0 & 0 & 0 & 0 & 0 & 0 & -1 & 0 & 0 & 1 & 0 & 0 & 0 & 0 & 0 & 0 & 0 \\
0 & 0 & -1 & 1 & 0 & 0 & 0 & 0 & 0 & 0 & 0 & 0 & 0 & 0 & -1 & 0 & 0 & 1 & 0 & 0 & 0 & 0 & 0 & 0 \\
0 & 0 & 0 & 0 & -1 & 1 & 0 & 0 & 0 & 0 & 0 & 0 & 0 & 0 & 0 & -1 & 0 & 0 & 1 & 0 & 0 & 0 & 0 & 0 \\
0 & 0 & 0 & 0 & 0 & -1 & 1 & 0 & 0 & 0 & 0 & 0 & 0 & 0 & 0 & 0 & -1 & 0 & 0 & 1 & 0 & 0 & 0 & 0 \\
0 & 0 & 0 & 0 & 0 & 0 & -1 & 1 & 0 & 0 & 0 & 0 & 0 & 0 & 0 & 0 & 0 & -1 & 0 & 0 & 1 & 0 & 0 & 0 \\
0 & 0 & 0 & 0 & 0 & 0 & 0 & 0 & -1 & 1 & 0 & 0 & 0 & 0 & 0 & 0 & 0 & 0 & -1 & 0 & 0 & 1 & 0 & 0 \\
0 & 0 & 0 & 0 & 0 & 0 & 0 & 0 & 0 & -1 & 1 & 0 & 0 & 0 & 0 & 0 & 0 & 0 & 0 & -1 & 0 & 0 & 1 & 0 \\
0 & 0 & 0 & 0 & 0 & 0 & 0 & 0 & 0 & 0 & -1 & 1 & 0 & 0 & 0 & 0 & 0 & 0 & 0 & 0 & -1 & 0 & 0 & 1 \\
\end{array} \right ) \left ( \begin{array}{c}
\bar{u}_1 \\
\bar{u}_2 \\
\bar{u}_3 \\
\bar{u}_4 \\
\bar{u}_5 \\
\bar{u}_6 \\
\bar{u}_7 \\
\bar{u}_8 \\
\bar{u}_9 \\
\bar{u}_{10} \\
\bar{u}_{11} \\
\bar{u}_{12} \\
\bar{u}_{13} \\
\bar{u}_{14} \\
\bar{u}_{15} \\
\bar{u}_{16} \\
\bar{u}_{17} \\
\bar{u}_{18} \\
\bar{u}_{19} \\
\bar{u}_{20} \\
\bar{u}_{21} \\
\bar{u}_{22} \\
\bar{u}_{23} \\
\bar{u}_{24}
\end{array} \right ) =
\left ( \begin{array}{c}
f_1 \\
f_2 \\
f_3 \\
f_4 \\
f_5 \\
f_6 \\
f_7 \\
f_8 \\
f_9
\end{array} \right) \;.
\]
}
Collecting all fluxes and \VARUN{source} terms in vectors $\mathsf{u}$ and $\mathsf{f}$, respectively, we can write this equation as
\begin{equation}
\tilde{\mathbb{E}}^{2,1} \mathsf{u} = \mathsf{f} \;.
\label{eq:topological_div}
\end{equation}
The matrix $\tilde{\mathbb{E}}^{2,1}$ is the {\em incidence matrix} which represents the divergence operator, not to be confused with $\mathbb{E}^{2,1}$ in (\ref{eq:topological_curl_v}) which represents the curl operator.

If, in the 2D case, the flow field is divergence-free, i.e. $f=0$, we know that a stream function $\boldsymbol{\psi}$ exists which is connected to $\u$ by
\[ u_x = \frac{\partial \boldsymbol{\psi}}{\partial y}\;,\quad u_y = - \frac{\partial \boldsymbol{\psi}}{\partial x} \;.\]
If we represent the stream function in the nodes of the grid shown in Figure~\ref{fig:outer_oriented_u}, then we have the exact topological equation
{\scriptsize
\[
\left ( \begin{array}{c}
\bar{u}_1 \\
\bar{u}_2 \\
\bar{u}_3 \\
\bar{u}_4 \\
\bar{u}_5 \\
\bar{u}_6 \\
\bar{u}_7 \\
\bar{u}_8 \\
\bar{u}_9 \\
\bar{u}_{10} \\
\bar{u}_{11} \\
\bar{u}_{12} \\
\bar{u}_{13} \\
\bar{u}_{14} \\
\bar{u}_{15} \\
\bar{u}_{16} \\
\bar{u}_{17} \\
\bar{u}_{18} \\
\bar{u}_{19} \\
\bar{u}_{20} \\
\bar{u}_{21} \\
\bar{u}_{22} \\
\bar{u}_{23} \\
\bar{u}_{24}
\end{array} \right ) =
\left ( \begin{array}{cccccccccccccccc}
-1 & 0 & 0 & 0 & 1 & 0 & 0 & 0 & 0 & 0 & 0 & 0 & 0 & 0 & 0 & 0 \\
0 & -1 & 0 & 0 & 0 & 1 & 0 & 0 & 0 & 0 & 0 & 0 & 0 & 0 & 0 & 0 \\
0 & 0 & -1 & 0 & 0 & 0 & 1 & 0 & 0 & 0 & 0 & 0 & 0 & 0 & 0 & 0 \\
0 & 0 & 0 & -1 & 0 & 0 & 0 & 1 & 0 & 0 & 0 & 0 & 0 & 0 & 0 & 0 \\
0 & 0 & 0 & 0 & -1 & 0 & 0 & 0 & 1 & 0 & 0 & 0 & 0 & 0 & 0 & 0 \\
0 & 0 & 0 & 0 & 0 & -1 & 0 & 0 & 0 & 1 & 0 & 0 & 0 & 0 & 0 & 0 \\
0 & 0 & 0 & 0 & 0 & 0 & -1 & 0 & 0 & 0 & 1 & 0 & 0 & 0 & 0 & 0 \\
0 & 0 & 0 & 0 & 0 & 0 & 0 & -1 & 0 & 0 & 0 & 1 & 0 & 0 & 0 & 0 \\
0 & 0 & 0 & 0 & 0 & 0 & 0 & 0 & -1 & 0 & 0 & 0 & 1 & 0 & 0 & 0 \\
0 & 0 & 0 & 0 & 0 & 0 & 0 & 0 & 0 & -1 & 0 & 0 & 0 & 1 & 0 & 0 \\
0 & 0 & 0 & 0 & 0 & 0 & 0 & 0 & 0 & 0 & -1 & 0 & 0 & 0 & 1 & 0 \\
0 & 0 & 0 & 0 & 0 & 0 & 0 & 0 & 0 & 0 & 0 & -1 & 0 & 0 & 0 & 1 \\
1 & -1 & 0 & 0 & 0 & 0 & 0 & 0 & 0 & 0 & 0 & 0 & 0 & 0 & 0 & 0 \\
0 & 1 & -1 & 0 & 0 & 0 & 0 & 0 & 0 & 0 & 0 & 0 & 0 & 0 & 0 & 0 \\
0 & 0 & 1 & -1 & 0 & 0 & 0 & 0 & 0 & 0 & 0 & 0 & 0 & 0 & 0 & 0 \\
0 & 0 & 0 & 0 & 1 & -1 & 0 & 0 & 0 & 0 & 0 & 0 & 0 & 0 & 0 & 0 \\
0 & 0 & 0 & 0 & 0 & 1 & -1 & 0 & 0 & 0 & 0 & 0 & 0 & 0 & 0 & 0 \\
0 & 0 & 0 & 0 & 0 & 0 & 1 & -1 & 0 & 0 & 0 & 0 & 0 & 0 & 0 & 0 \\
0 & 0 & 0 & 0 & 0 & 0 & 0 & 0 & 1 & -1 & 0 & 0 & 0 & 0 & 0 & 0 \\
0 & 0 & 0 & 0 & 0 & 0 & 0 & 0 & 0 & 1 & -1 & 0 & 0 & 0 & 0 & 0 \\
0 & 0 & 0 & 0 & 0 & 0 & 0 & 0 & 0 & 0 & 1 & -1 & 0 & 0 & 0 & 0 \\
0 & 0 & 0 & 0 & 0 & 0 & 0 & 0 & 0 & 0 & 0 & 0 & 1 & -1 & 0 & 0 \\
0 & 0 & 0 & 0 & 0 & 0 & 0 & 0 & 0 & 0 & 0 & 0 & 0 & 1 & -1 & 0 \\
0 & 0 & 0 & 0 & 0 & 0 & 0 & 0 & 0 & 0 & 0 & 0 & 0 & 0 & 1 & -1
\end{array} \right ) \left ( \begin{array}{c}
\psi_1 \\
\psi_2 \\
\psi_3 \\
\psi_4 \\
\psi_5 \\
\psi_6 \\
\psi_7 \\
\psi_8 \\
\psi_9 \\
\psi_{10} \\
\psi_{11} \\
\psi_{12} \\
\psi_{13} \\
\psi_{14} \\
\psi_{15} \\
\psi_{16}
\end{array} \right ) \;.
\]
}
We can write this in terms of incidence matrices as\footnote{Note that if we performed the same steps in 3D, then the divergence relation (\ref{eq:topological_div}) would be
\[ \tilde{\mathbb{E}}^{3,2} \mathsf{u} = \mathsf{f} \;,\]
and the \VARUN{2D} stream function becomes the \VARUN{3D} stream vector field and we would have
\[ \mathsf{u} = \tilde{\mathbb{E}}^{2,1} \vec{\psi} \;.\]
So clearly the incidence matrices $\tilde{\mathbb{E}}$ depend on the dimension of the space $d$ in which the problem is posed. Note that this is not the case for the incidence matrices $\mathbb{E}$. Alternatively, we could refer to the dimension-dependent incidence matrices as
\[ \mathbb{E}^{d,d-1} = \left \{ \begin{array}{ll}
\tilde{\mathbb{E}}^{2,1} \quad & \mbox{if } d=2 \\[1ex]
\tilde{\mathbb{E}}^{3,2} \quad & \mbox{if } d=3
\end{array} \right . \quad \mbox{and} \quad
\mathbb{E}^{d-1,d-2} = \left \{ \begin{array}{ll}
\tilde{\mathbb{E}}^{1,0} \quad & \mbox{if } d=2 \\[1ex]
\tilde{\mathbb{E}}^{2,1} \quad & \mbox{if } d=3
\end{array} \right . \;,
\]
in which case it is immediately clear  that these matrices depend on the $d$. From now on we will use the incidence matrices with the $d$, because then the results are valid for any space dimension $d$.
\label{footnote:incidence_matrices}
}
\begin{equation}
\mathsf{u} = \tilde{\mathbb{E}}^{1,0} \vec{\psi} \;.
\label{eq:flux_streamfunction_relation}
\end{equation}
If the flux $\u$ is prescribed along the $\Gamma_u$ the associated edges (2D) or surfaces (3D) can be eliminated from the system $\mathbb{E}^{d,d-1} \VARUN{\mathsf{u}} = \VARUN{\mathsf{f}}$ and transferred to the right hand side.

For the discretization of (\ref{eq:extended_Darcy_problem}) the first and last equation in that system can be represented on the mesh by
\[
\left \{ \begin{array}{l}
\mathsf{v} - \mathbb{E}^{1,0} \mathsf{p} = 0 \\[1ex]
\VARUN{\tilde{\mathbb{E}}}^{d,d-1} \mathsf{u} = \mathsf{f}
\end{array} \right .\;.
\]
Prescription of boundary conditions $p$ along $\Gamma_u$ and $\u$ along $\Gamma_u$ can be done strongly. The degrees of freedom can be eliminated and transferred to the right hand side. The equation between $p$ and $\bv$ is exact on any grid and the discrete divergence relation between $\u$ and $f$ is exact on any grid. Note the $(\bv,p)$-grid is not necessarily the $(\u,f)$-grid, so in principle we can use different grids for both equations.

Unfortunately, neither of the two problems, $\bv = \nabla p$ and $\nabla \cdot \u =f$ has a unique solution on their respective grids. It is the final equation in (\ref{eq:extended_Darcy_problem}), $\u = -\mathbb{K} \bv$, that couples the solution on the two grids and renders a unique solution. It is also in this equation that the numerical approximation is made; the more accurate we approximate this algebraic equation, the more accurate the solution to the first order system (\ref{eq:extended_Darcy_problem}) \VARUN{will be}.

For many numerical methods\footnote{A notable exception is the class of least-squares formulation\VARUN{s} which aim\VARUN{s} to {\em minimize} the expression $\u + \mathbb{K}\bv$ \cite{BochevGunzbuger}.} well-posedness requires that the number of discrete degrees of freedom $\bar{\bv}_i$ is equal to the discrete number of degrees of freedom $\bar{\u}_j$, or more geometrically, that the number of $k$-dimensional geometric objects on one grid is equal to the number of $(d-k)$-dimensional geometric objects on the other grid. Here $k=0$ refers to points in the grid, $k=1$ to edges in the grid, $k=2$ the faces in the grid\VARUN{,} and $k=3$ the volumes \VARUN{in the grid}.

The requirement $\#k=\#(d-k)$ cannot be accomplished on a single grid, so this requires two different grids which are constructed in such a way that $\#k=\#(d-k)$ holds, \cite{Bonelle2014,Bonelle2015a,desbrun2005discrete,Kreeft2011,Nicolaides,Robidoux2011}.
\begin{figure}[h!]
	\centering
	\includegraphics[width=0.6\linewidth]{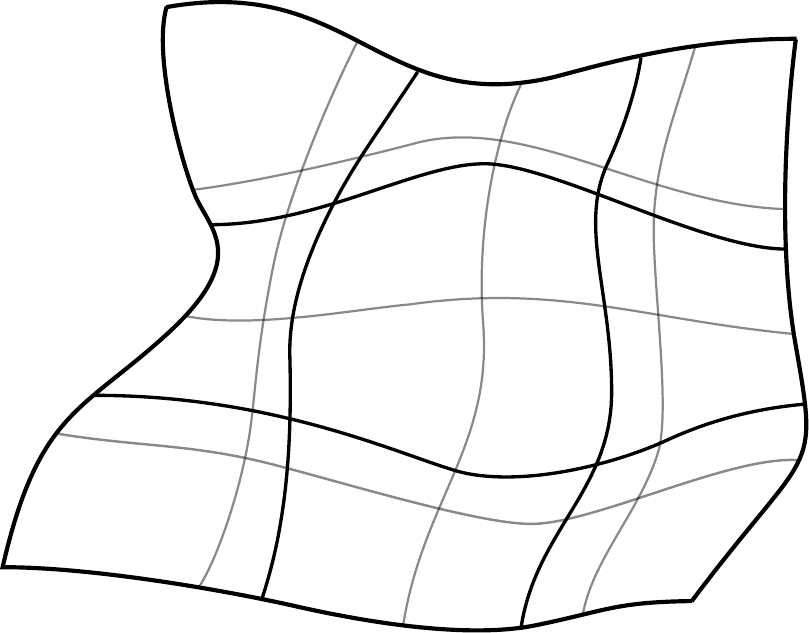}
	\caption{The primal grid (thin gray) where $(\bv,p)$ are represented and the dual grid (thick black) where $(\u,f)$ are represented. Note that $\Gamma_p=\partial \Omega$ and consequently $\Gamma_u=\varnothing$.}
	\label{fig:primal_dual}
\end{figure}

A dual grid complex is shown in Figure~\ref{fig:primal_dual}. The integral quantities $(\bv,p)$ can be represented on the gray grid. If $p$ is prescribed along the entire boundary, then those degrees of freedom are eliminated (including the \VARUN{gray} edges along the boundary for which the integral value $\bv$ is then known also), see for instance Figure~\ref{fig:inner_oriented_p_v_bc}. In that case flux $\u$ along the boundary cannot be prescribed. \VARUN{In Figure~\ref{fig:primal_dual}, the number of points in the gray grid, $9$, equals the number of surfaces in the black grid, the number of edges in the grey grid is equal to the number of edges on the black grid, $24$\VARUN{,} and the number of surfaces on the gray grid equals the number of points in the black grid, $16$, therefore, we have $\#k=\#(d-k)$ for $d=2$}.

Alternatively, we could have represented $(\u,f)$ on the gray grid with $\u$ and the stream function $\psi$ prescribed and $(\bv,p)$ on the black grid.In this case $\Gamma_u = \partial \Omega$ and $\Gamma_p = \varnothing$.

\subsection{Dual grids}
If dual grids, such as described above, are employed then we have two properties:
\begin{enumerate}
\item There exists a square, invertible matrix $\mathbb{H}_{\mathbb{K}}^{d-1,1}$ such that $\mathsf{u} = \mathbb{H}_{\mathbb{K}}^{d-1,1} \mathsf{v}$.
\item The incidence matrices on the primal and dual grid satisfy\footnote{This relation is true if the orientations on primal and dual grid agree. This is not always the case and then the relation reads $\mathbb{E}^{d-k,d-k-1} = -\mathbb{E}^{{k,k-1}^T}$. A well known example is the duality between grad and div.}
    \[  \mathbb{E}^{d-k,d-k-1} = \VARUN{\left( \mathbb{E}^{{k,k-1}} \right)} ^T \;.\]
\end{enumerate}
If we use dual grids and these properties hold, we can write (\ref{eq:extended_Darcy_problem}) as
\begin{equation}
\left \{ \begin{array}{l}
\mathsf{v} - \mathbb{E}^{1,0} \mathsf{p} = 0\\[1ex]
\mathsf{u} - \mathbb{H}_{\mathbb{K}}^{d-1,1} \mathsf{v} = 0  \\[1ex]
\mathbb{E}^{{1,0}^T} \mathsf{u} = \mathsf{f}
\end{array} \right . \;,
\label{eq:discrete_extended_Darcy_problem}
\end{equation}
where the vectors $\mathsf{p}$, $\mathsf{v}$, $\mathsf{u}$ and $\mathsf{f}$ contain the integral quantities in the mesh as discussed in the previous sections.

In the diagram below, we place the various integral values in appropriate `spaces'
\[
		\xymatrix{
			 \mathsf{p} \in \textcolor{not-used}{H_{0}} \ar@[not-used][r]^-{\textcolor{not-used}{\mathbb{E}^{1,0}}} \ar@_{->}@[not-used][d]<-0.25ex>_-{\textcolor{not-used}{\mathbb{H}^{d,0}}}  & \mathsf{v}\in H_{1} \ar[r]^-{\mathbb{E}^{2,1}} \ar@_{->}[d]<-0.25ex>_-{\mathbb{H}_{\mathbb{K}}^{d-1,1}} & \xi \in H_{2} \\ 
			\mathsf{f} \in \tilde{H}_d \ar@_{->}@[not-used][u]<-0.25ex>_-{\textcolor{not-used}{\mathbb{H}^{0,d}}} & \ar[l]^-{\mathbb{E}^{{1,0}^T}} \mathsf{u}\in \tilde{H}_{d-1}  \ar@_{->}[u]<-0.25ex>_-{\mathbb{H}_{\mathbb{K}^{-1}}^{1,d-1}}  & \ar@[not-used][l]^-{\color{not-used}{\mathbb{E}^{{2,1}^T}}} \color{not-used}{\psi \in \tilde{H}_{d-2}} 
		}
\label{eq:doubleDeRham}
	\]
Here $H_k$ denotes the space of values assigned to $k$-dimensional objects in the $H$-grid for $k=0,1,2$. If $\tilde{H}$ denotes the dual grid, then $\tilde{H}_l$ is the space of values assigned to $l$-dimensional objects in the $\tilde{H}$-grid.

For dual grids the number of points in the $H$-grid is equal to the number of $d$-dimensional volumes in the dual grid $\tilde{H}$. Let $\mathbb{H}^{d,0}$ and $\mathbb{H}^{0,d}$ be square, invertible matrices which map between $H_0$ and $\tilde{H}_d$ as shown in the diagram above.

If we eliminate $\mathsf{v}$ and $\mathsf{u}$ from (\ref{eq:discrete_extended_Darcy_problem}) we have
\begin{equation}
\mathbb{E}^{{1,0}^T} \mathbb{H}_{\mathbb{K}}^{d-1,1} \mathbb{E}^{1,0} \mathsf{p} = \mathsf{f} \;.
\label{eq:generic_direct_formulation}
\end{equation}
This discretization corresponds to (\ref{eq:anisotropic_diffusion}). We will refer to this formulation as the {\em direct formulation}.

If $\VARUN{\mathsf{p}}\in \tilde{H}_d$ we  can set up the diffusion problem as
\begin{equation}
\left \{ \begin{array}{lcl}
- \mathbb{H}_{\mathbb{K}^{-1}}^{1,d-1} \mathsf{u} + \mathbb{E}^{{d,d-1}^T} \mathbb{H}^{0,d} \VARUN{\mathsf{p}} & = & 0 \\[1ex]
\mathbb{H}^{0,d} \mathbb{E}^{{d,d-1}} \mathsf{u} & = & 
\mathsf{f}
\end{array} \right . \;.
\label{eq:generic_mixed_formulation}
\end{equation}
This formulation, where we solve for $\VARUN{p}$ and $\u$ simultaneously, resembles (\ref{eq:Darcy_problem}), and will be called the {\em mixed formulation}, \cite{brezzi1991mixed}.

\section{Mimetic spectral element method}\label{sec:MSEM}

The incidence matrices introduced in the previous section are generic and only depend on the grid topology. The matrices $\mathbb{H}$ which switch between the primal and the dual grid representation explicitly depend on the numerical method that is used. In this section we will introduce spectral element functions which interpolate the integral values in a grid. With these functions we can construct the $\mathbb{H}$-matrices, which turn out to be (weighted) finite element mass matrices. The derivation in this section will be on an orthogonal grid. The extension to curvilinear grids will \VARUN{be} discussed in the next section.

\subsection{One dimensional spectral basis functions}
Consider the interval $[-1,1]\subset \mathbb{R}$ and the Legendre polynomials, $L_N(\xi)$\VARUN{,} of degree $N$, $\xi \in [-1,1]$. The $(N+1)$ roots, $\xi_i$, of the polynomial $(1-\xi^2)L_N'(\xi)$ satisfy $-1 \leq \xi_i \leq 1$. Here $L_N'(\xi)$ is the derivative of the Legendre polynomial. The roots $\xi_i$, $i=0,\ldots,N$, are called the {\em Gauss-Lobatto-Legendre (GLL) points}, \cite{Canuto_et_al}. Let $h_i(\xi)$ be the Lagrange polynomial through the GLL points such that
\begin{equation}
h_i(\xi_j) = \left \{ \begin{array}{ll}
1 \quad \quad & \mbox{if } i=j \\
 & \\
0 \quad \quad & \mbox{if } i \neq j
\end{array} \right . \quad \quad i,j = 0,\ldots N \;.
\label{eq:Kronecker_delta}
\end{equation}
The explicit form of the Lagrange polynomials in terms of the Legendre polynomials is given by
\begin{equation}
h_i(\xi) = \frac{(1-\xi^2) L_N'(\xi)}{N(N+1)L_N(\xi_i)(\xi_i - \xi)} \;.
\end{equation}
Let $f(\xi)$ be a function defined for $\xi \in [-1,1]$ by
\begin{equation}
f(\xi) = \sum_{i=0}^N a_i h_i(\xi) \;.
\label{eq:nodal_expansion}
\end{equation}
Using property (\ref{eq:Kronecker_delta}) we see that $f(\xi_j) = a_j$, so the expansion coefficients in (\ref{eq:nodal_expansion}) coincide with the value of $f$ in the GLL nodes. We will refer to this expansion as a {\em nodal expansion}, because the expansion coefficients, $a_i$ in (\ref{eq:nodal_expansion}) are the value of $f(\xi)$ in the {\em nodes} $\xi_i$. The basis functions $h_i(\xi)$ are polynomials of degree $N$.

From the nodal basis functions, define the polynomials $e_i(\xi)$ by
\begin{equation}
e_i(\xi) = - \sum_{k=0}^{i-1} \frac{dh_k(\xi)}{d\xi} \;.
\end{equation}
The functions $e_i(\xi)$ are polynomials of degree $(N-1)$. These polynomials satisfy, \cite{gerritsma::edge_basis,Kreeft2011,Palha2014}
\begin{equation}
\int_{\xi_{j-1}}^{\xi_j} e_i(\xi) = \left \{ \begin{array}{ll}
1 \quad \quad & \mbox{if } i=j \\
 & \\
0 \quad \quad & \mbox{if } i \neq j
\end{array} \right . \quad \quad i,j = 1,\ldots N \;.
\label{eq:edge_function_property}
\end{equation}
Let a function $\VARUN{g}(\xi)$ be expanded in these functions
\begin{equation}
\VARUN{g}(\xi) = \sum_{i=1}^N \VARUN{b}_i e_i(\xi) \;,
\label{eq:edge_expansion}
\end{equation}
then using (\ref{eq:edge_function_property})
\[ \int_{\xi_{j-1}}^{\xi_j} \VARUN{g}(\xi) = \VARUN{b}_j \;.\]
So the expansion coefficients $b_i$ in (\ref{eq:edge_expansion}) coincide with the integral of $g$ over the edge $[\xi_{i-1},\xi_i]$. We will call these basis functions {\em edge functions} and refer to the expansion (\ref{eq:edge_expansion}) as an {\em edge expansion}, see for instance \cite{bochev_mimetic_ls_2014,Kreeft2011,Palha2014} for examples of nodal and edge expansions.

Let $f(\xi)$ be expanded in terms Lagrange polynomials as in (\ref{eq:nodal_expansion}), then the derivative\footnote{Note that the set of polynomials $\{ h_i'\}$, $i=0,\ldots,N$ is linearly dependent and therefore does not form a basis, while the set $\{ e_i \}$, $i=1,\ldots,N$ is linearly independent and therefore forms a basis for the derivatives of the nodal expansion (\ref{eq:nodal_expansion}).} of $f$ is given by, \cite{gerritsma::edge_basis,Kreeft2011,Palha2014}
\begin{equation}
f'(\xi) = \sum_{i=0}^N a_i h_i'(\xi) = \sum_{i=1}^N ( a_i - a_{i-1}) e_i(\xi) \;.
\label{eq:derivative_nodla2edge}
\end{equation}


If we collect all the expansion coefficients in a column vector and all the basis functions in a row vector we have
\begin{equation}
f(\xi) = [ h_0 \; h_1\;\ldots \,\; h_N ] \, \left [ \begin{array}{c}
a_0 \\
 \\
\vdots \\
 \\
a_N
\end{array} \right ]\;,
\end{equation}
then the derivative is given by\footnote{The matrix $\mathbb{E}^{1,0}$ is the {\em incidence matrix} as was discussed in Sections~\ref{sec:gradient} and \ref{sec:divergence}. It takes the nodal expansion coefficients and maps them to the edge expansion coefficients. The incidence matrix is the topological part of the derivative. It is independent of the order of the method (the polynomial degree N) and the size or the shape of the mesh. The incidence matrix only depends on the topology and orientation of the grid, see \cite{bochev2006principles,kreeft::stokes,Kreeft2011}.} (\ref{eq:derivative_nodla2edge})
\begin{equation}
f'(\xi) = [e_1\;\ldots \,\; e_N ] \,\left (
 \begin{array}{cccccc}
 -1 & 1 & 0 & \dots & & 0 \\
   & \ddots & \ddots & & & 0\\
   &  & -1 & 1 & & 0 \\
   &  &  &  \ddots & \ddots & \\
 0 &  & \dots & 0 & -1 & 1
 \end{array}
 \right ) \left [ \begin{array}{c}
a_0 \\
 \\
\vdots \\
 \\
a_N
\end{array} \right ] = [e_1\;\ldots \,\; e_N ]  \mathbb{E}^{1,0} \left [ \begin{array}{c}
a_0 \\
 \\
\vdots \\
 \\
a_N
\end{array}\right ]\;.
\label{eq:expansion-derivative_f}
\end{equation}
So taking the derivative essentially consists of two step: Apply the matrix $\mathbb{E}^{1,0}$ to the expansion coefficients and expand in a new basis.

\subsection{Two dimensional expansions}

\subsubsection{Expanding $p$ (Direct formulation)}\label{subsec:direct_formulation}
In finite element methods the direct finite element formulation for the anisotropic diffusion problem is given by: For $(\mathbb{K}\nabla p, \bm{n})=0$ along $\Gamma _u$ and $f \in H^{-1} (\Omega)$, find $p \in \VARUN{H_{0, \Gamma_p}^{1}(\Omega)}$ such that
\begin{equation}
\left ( \nabla \tilde{p}, \mathbb{K} \nabla p \right ) = (\tilde{p},f )\;,\;\;\; \forall \tilde{p} \in \VARUN{H_{0, \Gamma_p}^{1}(\Omega)} \;.
\label{eq:FE_direct_formulation}
\end{equation}
\VARUN{where $H_{0, \Gamma_p}^{1} = \lbrace p \in \hone \vert p = 0$ on $\Gamma_p \rbrace$.}

Consider $[-1,1]^2 \subset \mathbb{R}^2$ and let $p(\xi,\eta)$ be expanded as
\begin{equation}
p(\xi,\eta) = \sum_{i=0}^N \sum_{j=0}^N p_{i,j}h_i(\xi) h_j(\eta) \;.
\label{eq:nodal_pressure_expansion}
\end{equation}
From (\ref{eq:Kronecker_delta}) it follows that $p_{i,j} = p(\xi_i,\eta_j)$.
If we take the gradient of $p$ using (\ref{eq:derivative_nodla2edge}) we have
\begin{align}
\nabla p & = \left ( \begin{array}{c}
\sum_{i=1}^N \sum_{j=0}^N ( p_{i,j} - p_{i-1,j}) e_i(\xi) h_j(\eta) \\[1ex]
\sum_{i=0}^N \sum_{j=1}^N ( p_{i,j} - p_{i,j-1}) h_i(\xi) e_j(\eta)
\end{array} \right ) \\
& = \left ( \begin{array}{cccccc}
e_1(\xi) h_0(\eta) & \ldots & e_N(\xi) h_N(\eta) & 0 & \ldots & 0 \\[1ex]
0 & \ldots & 0 & h_0(\xi)e_1(\eta) & \ldots & h_N(\xi)e_N(\eta)
\end{array} \right )
\mathbb{E}^{1,0} \left [
\begin{array}{c}
p_{0,0} \\
\vdots \\
p_{N,N}
\end{array} \right ]
\nonumber \\
& = \left ( \begin{array}{cccccc}
e_1(\xi) h_0(\eta) & \ldots & e_N(\xi) h_N(\eta) & 0 & \ldots & 0 \\[1ex]
0 & \ldots & 0 & h_0(\xi)e_1(\eta) & \ldots & h_N(\xi)e_N(\eta)
\end{array} \right )
\mathbb{E}^{1,0} \mathsf{p} \;.
\end{align}
If we insert this in (\ref{eq:FE_direct_formulation}), we have
\begin{equation}
\VARUN{\left( \mathbb{E}^{{1,0}}\right)} ^T \mathbb{M}_{\mathbb{K}}^{(1)} \mathbb{E}^{1,0} \mathsf{p} = \mathsf{f} \;,
\label{eq:discrete_direct_formulation_spectral}
\end{equation}
where
\begin{equation}
\mathbb{M}_{\mathbb{K}}^{(1)} = \iint_{\Omega} \left (\begin{array}{cc}
e_1(\xi) h_0(\eta) & 0 \\
\vdots & \vdots \\
e_N(\xi) h_N(\eta) & 0 \\
0 & h_0(\xi)e_1(\eta) \\
\vdots & \vdots \\
0 & h_N(\xi)e_N(\eta)
\end{array} \right ) \, \mathbb{K}
\left ( \begin{array}{cccccc}
e_1(\xi) h_0(\eta) & \ldots & e_N(\xi) h_N(\eta) & 0 & \ldots & 0 \\[1ex]
0 & \ldots & 0 & h_0(\xi)e_1(\eta) & \ldots & h_N(\xi)e_N(\eta)
\end{array} \right ) \, \mathrm{d} \Omega \;,
\label{eq:weighted_M1_mass_matrix}
\end{equation}
and $\mathsf{p}$ is the vector which contains the expansion coefficients of $p(\xi,\eta)$ in (\ref{eq:nodal_pressure_expansion}).

The vector $\mathsf{f}$ in (\ref{eq:discrete_direct_formulation_spectral}) is given by
\[
\mathsf{f} = \iint_{\Omega} \left ( \begin{array}{c}
h_0(\xi) h_0(\eta) \\
\vdots \\
h_N(\xi)h_N(\eta)
\end{array} \right)
\,
f(\xi,\eta) \,\mathrm{d} \Omega \;.
\]

If we compare (\ref{eq:discrete_direct_formulation_spectral}) with (\ref{eq:generic_direct_formulation}), we see that the $\mathbb{H}_{\mathbb{K}}^{d-1,1}$-matrix from (\ref{eq:generic_direct_formulation}) is represented in the finite element formulation by the weighted mass matrix $\mathbb{M}_{\mathbb{K}}^{(1)}$ given by (\ref{eq:weighted_M1_mass_matrix}), see also \cite{bochev2006principles,tarhasaari1999some}.

\subsubsection{Expanding $\u$ and $p$ (Mixed formulation)}\label{subsec:mixed_formulation}
The mixed formulation for the anisotropic steady diffusion problem is given by: For $p=0$ along $\Gamma _p$ and for $f \in L^2(\Omega)$, find $u \in H_{0,\Gamma_n}(div;\Omega)$ such that
\begin{equation}
\left \{
\begin{array}{lcrl}
-( \tilde{\u}, \mathbb{K}^{-1} \u ) + ( \nabla \cdot \tilde{\u}, p ) &=& 0 \quad \quad & \forall \tilde{\u} \in H_{0,\Gamma_u}(div;\Omega) \\[1ex]
(\tilde{p}, \nabla \cdot \u ) &=&f \quad \quad & \forall \tilde{p} \in L^2(\Omega)
\end{array}
\right . \;.
\label{eq:continuous_mixed_fomrulation}
\end{equation}
\VARUN{where, $ H_{0,\Gamma_u}(div;\Omega) = \lbrace u \in H(div;\Omega) \vert u \cdot n =0$ along $\Gamma_u \rbrace $.}

In contrast to the pressure expansion in Section~\ref{subsec:direct_formulation} in the direct formulation, (\ref{eq:nodal_pressure_expansion}), in the mixed formulation the pressure is expanded in terms of edge functions
\begin{equation}
p(\xi,\eta) = \sum_{i=1}^N\sum_{j=1}^N p_{i,j} e_i(\xi) e_j(\eta) \;.
\label{eq:volume_pressure_expansion}
\end{equation}
The velocity $\u$ is expanded as
\begin{align}
\u = & \left ( \begin{array}{c} u \\ v \end{array} \right ) =
\left ( \begin{array}{c}
\sum_{i=0}^N \sum_{j=1}^N u_{i,j} h_i(\xi) e_j(\eta) \\[2ex]
\sum_{i=1}^N \sum_{j=0}^N v_{i,j} e_i(\xi) h_j(\eta)
\end{array}
\right ) \label{eq:u_flux_expansion} \\
 & = \left ( \begin{array}{cccccc}
 h_0(\xi)e_1(\eta) & \dots & h_N(\xi) e_N(\eta) & 0 & \dots & 0 \\
 0 & \dots & 0 & e_1(\xi) h_0(\eta) & \dots & e_N(\xi) h_N(\eta)
 \end{array} \right ) \left ( \begin{array}{c}
 u_{0,1} \\
 \vdots \\
 u_{N,N} \\
 v_{1,0} \\
 \vdots \\
 v_{N,N}
 \end{array} \right )\;.\nonumber
\end{align}
Application of the divergence operator to (\ref{eq:u_flux_expansion}) and using (\ref{eq:derivative_nodla2edge}) we obtain
\begin{align}
\nabla \cdot \u  & = \sum_{i=1}^N \sum_{j=1}^N (u_{i,j}-u_{i-1,j} + v_{i,j} - v_{i,j-1} ) e_i(\xi) e_j(\eta) \label{eq:dicsrete_div_velocity}\\
 & = \left ( \begin{array}{ccc}
 e_1(\xi)e_1(\eta) & \ldots & e_N(\xi)e_N(\eta)
 \end{array} \right ) \mathbb{E}^{d,d-1} \left ( \begin{array}{c}
 u_{0,1} \\
 \vdots \\
 u_{N,N} \\
 v_{1,0} \\
 \vdots \\
 v_{N,N}
 \end{array} \right ) \nonumber \\
& =  \left ( \begin{array}{ccc}
 e_1(\xi)e_1(\eta) & \ldots & e_N(\xi)e_N(\eta)
 \end{array} \right ) \mathbb{E}^{d,d-1} \mathsf{u} \;. \nonumber
\end{align}
Note that $\mathbb{E}^{d,d-1}$ is the incidence matrix which also appeared in (\ref{eq:topological_div}) and footnote~\ref{footnote:incidence_matrices}.

If we insert the expansion (\ref{eq:u_flux_expansion}) in $(\tilde{\u},\mathbb{K}^{-1} \u )$ we obtain
\begin{equation}
(\tilde{\u},\mathbb{K}^{-1} \u ) = \tilde{\mathsf{u}}^T \mathbb{M}_{\mathbb{K}^{-1}}^{(d-1)} \mathsf{u} \;,
\label{eq:anisotropic_block}
\end{equation}
with
\begin{align}
&\mathbb{M}_{\mathbb{K}^{-1}}^{(d-1)} = \\
 & \iint_{\Omega} \left ( \begin{array}{cc}
 h_0(\xi)e_1(\eta) & 0 \\
 \vdots & \vdots \\
 h_N(\xi) e_N(\eta) & 0 \\
 0 & e_1(\xi) h_0(\eta) \\
 \vdots & \vdots \\
 0 & e_N(\xi) h_N(\eta)
 \end{array} \right )
 \mathbb{K}^{-1} \left ( \begin{array}{cccccc}
 h_0(\xi)e_1(\eta) & \dots & h_N(\xi) e_N(\eta) & 0 & \dots & 0 \\
 0 & \dots & 0 & e_1(\xi) h_0(\eta) & \dots & e_N(\xi) h_N(\eta)
 \end{array} \right ) \, \mathrm{d}\Omega \;.
\end{align} \label{eq:k_inv}

Note that pressure is expanded in the same basis as the divergence of the velocity field, (\ref{eq:volume_pressure_expansion}) and (\ref{eq:dicsrete_div_velocity}), therefore we can write
\begin{equation}
 ( \tilde{p}, \nabla \cdot \u ) = \tilde{\mathsf{p}}^T \mathbb{M}^{(d)} \mathbb{E}^{d,d-1} \mathsf{u} \;,
 \label{eq:div_block}
 \end{equation}
with
\[ \mathbb{M}^{(\VARUN{d})} = \iint_{\Omega} \left ( \begin{array}{c}
e_1(\xi) e_1(\eta) \\
\vdots \\
e_N(\xi)e_N(\eta)
\end{array} \right )
\left ( \begin{array}{ccc}
e_1(\xi) e_1(\eta) & \dots & e_N(\xi)e_N(\eta)
\end{array} \right )\,\mathrm{d} \Omega \;.
\]
With (\ref{eq:anisotropic_block}) and (\ref{eq:div_block}) we can write (\ref{eq:continuous_mixed_fomrulation}) as
\begin{equation}
\left \{ \begin{array}{lcr}
\mathbb{M}_{\mathbb{K}^{-1}}^{(d-1)}\mathsf{u} +   \mathbb{E}^{{d,d-1}^T} \mathbb{M}^{(d)} \mathsf{p} & = & 0  \\
\mathbb{M}^{(d)} \mathbb{E}^{d,d-1} \mathsf{u} & = & \mathsf{f}
\end{array} \right .\;,
\label{eq:finite_element_mixed_formulation}
\end{equation}
with
\[ \mathsf{f} = \iint_{\Omega} \left ( \begin{array}{c}
e_1(\xi) e_1(\eta) \\
\vdots \\
e_N(\xi)e_N(\eta)
\end{array} \right ) f(\xi,\eta) \, \mathrm{d} \Omega \;.
\]
Comparison of (\ref{eq:finite_element_mixed_formulation}) with (\ref{eq:generic_mixed_formulation}) shows that the topological incidence matrices also appear in the finite element formulation and that the (weighted) mass matrices $\mathbb{M}_{\mathbb{K}^{-1}}^{(d-1)}$ and $\mathbb{M}^{(d)}$ once again play the role of the $\mathbb{H}$-matrices which connect solutions on dual grids.

\VARUN{In this section only the discretization on a single spectral element is discussed. Transformation of the domain $[-1,1]^2$ to more general domains will be discussed in Section~\ref{sec:transformation}. The use of multiple elements follows the general assembly procedure from finite element methods. Results of this approach are presented in Section~\ref{sec:numerical_results}.}

\section{Transformation rules}\label{sec:transformation}
	The basis functions used in the discretization of the different physical field quantities have only been introduced for the reference domain $\tilde{\Omega} = [-1,1]^{2}$. For these basis functions to be applicable in a different domain $\Omega$, it is fundamental to discuss how they transform under a mapping $\Phi: (\xi,\eta)\in\tilde{\Omega}\mapsto (x,y)\in\Omega\subset\mathbb{R}^{2}$. Within a finite element formulation this is particularly useful because the basis functions in the reference domain $\tilde{\Omega}$ can then be transformed to each of the elements $\Omega_{e}$, given a mapping $\Phi_{e}:\tilde{\Omega}\mapsto\Omega_{e}$.
	
	Consider a smooth bijective map $\Phi:(\xi,\eta)\in\tilde{\Omega}\mapsto (x,y)\in\Omega$ such that
	\[
		x = \Phi^{x}(\xi,\eta)\quad \text{and} \quad y = \Phi^{y}(\xi,\eta)\,,
	\]
	and the associated rank two Jacobian tensor $\moperator{J}$
	\[
		\arraycolsep=1.4pt\def\arraystretch{2.4}
		\moperator{J} :=
		\left[
			\begin{array}{cc}
				\dfrac{\partial\Phi^{x}}{\partial\xi} & \dfrac{\partial\Phi^{x}}{\partial\eta} \\
				\dfrac{\partial\Phi^{y}}{\partial\xi} & \dfrac{\partial\Phi^{y}}{\partial\eta}
			\end{array}
		\right]\,.
	\]
	
	The transformation of a scalar function $\varphi$ discretized by nodal values is given by
	\begin{equation}
		\tilde{\varphi}(\xi,\eta) = (\varphi\circ\Phi)(\xi,\eta) \quad \text{and}\quad \varphi(x,y) = (\tilde{\varphi}\circ\Phi^{-1})(x,y), \label{eq:nodal_functions_transformation}
	\end{equation}
	and of a scalar function $\rho$ discretized by surface integrals is given by
	\begin{equation}
		\tilde{\rho}(\xi,\eta) = \det \moperator{J}\,(\rho\circ\Phi)(\xi,\eta) \quad \text{and}\quad \rho(x,y) = \frac{1}{\det\moperator{J}}(\tilde{\rho}\circ\Phi^{-1})(x,y). \label{eq:volume_functions_transformation}
	\end{equation}
	The transformation of vector fields $\bv$ discretized by line integrals is
	\begin{equation}
		\tilde{\bv}(\xi,\eta) = \moperator{J}^{\transpose}(\bv\circ\Phi)(\xi,\eta) \quad \text{and}\quad \bv(x,y) = (\moperator{J}^{\transpose})^{-1}(\tilde{\bv}\circ\Phi^{-1})(x,y), \label{eq:line_functions_transformation}
	\end{equation}
	 and of vector fields $\u$ discretized by flux integrals is
	 \begin{equation}
		\tilde{\u}(\xi,\eta) = \det\moperator{J}\,\,\moperator{J}^{-1}(\u\circ\Phi)(\xi,\eta) \quad \text{and}\quad \u(x,y) = \frac{1}{\det\moperator{J}}\moperator{J}(\tilde{\u}\circ\Phi^{-1})(x,y). \label{eq:flux_functions_transformation}
	\end{equation}
	
		These transformations affect only the mass matrices and not the incidence matrices. This is fundamental to ensure the topological nature of the incidence matrices.
\section{Numerical results}\label{sec:numerical_results}
In this section three test cases are presented to illustrate the accuracy of the discretization scheme developed in this work.
The first test case, \ref{sec:test1}, is an analytical solution taken from \cite{Herbin2008} to assess the convergence rates of the method.
The second test case, \ref{sec:SS}, is the flow through a system of sand and shale blocks with highly heterogeneous permeability in the domain, see for more details \cite{Durlofsky1994}.
The third test case, \ref{sec:test3}, is a highly anisotropic and heterogeneous permeability tensor in the domain, see for more details, \cite{Durlofsky1993}.

\subsection{Manufactured solution} \label{sec:test1}
We first test the  method using the exact solution
\begin{equation}
p _{\mathrm{exact}}(x,y)=\sin(\pi x)\sin(\pi y),
\end{equation}
with the permeability tensor given by
\begin{equation}
\mathbb{K}=\dfrac{1}{\left( x^{2}+y^{2}+\alpha\right) }
\begin{pmatrix}
10^{-3}x^{2}+y^{2}+\alpha   &   \left( 10^{-3}-1\right)xy\\
\left( 10^{-3}-1\right)xy & x^{2}+10^{-3}y^{2} +\alpha
\end{pmatrix}.
\end{equation}
The mixed formulation (\ref{eq:Darcy_problem}) in the form of (\ref{eq:finite_element_mixed_formulation}) is then solved in the domain $ \left( x,y \right)\in \Omega=[0,1]^{2} $ with the source term $ f=-\nabla\cdot ( \mathbb{K}\nabla p _{\mathrm{exact}})  $ and the Dirichlet boundary condition $ \left. p\right| _{\partial \Omega}=0 $. A benchmark of this test case for $ \alpha=0 $ using multiple numerical schemes can be found in \cite{Herbin2008}.

When $ \alpha=0 $, $ \mathbb{K} $ is multi-valued at the origin which makes this test case a challenging one. To see this, we can first convert the Cartesian coordinates $ (x,y) $ to polar coordinates $ (r,\theta) $ by $ x=r\cos\theta $, $ y=r\sin\theta $. Then we have
\begin{equation}
\left. \mathbb{K}\right| _{\alpha=0}=
\begin{pmatrix}
10^{-3}\cos^{2}\theta+\sin^{2}\theta   &   \left( 10^{-3}-1\right)\cos\theta\sin\theta\\
 \left( 10^{-3}-1\right)\cos\theta\sin\theta & \cos^{2}\theta+10^{-3}\sin^{2}\theta
\end{pmatrix}.
\end{equation}
It can be seen that we get different $ \left. \mathbb{K}\right| _{\alpha=0} $ when we approach the origin along different angles, $\theta$.
It must be noted that inverse of $\mathbb{K}$ does not exist at the origin. The inverse of the tensor term appears in \ref{eq:k_inv}. We use Gauss integration and thus the inverse term is not evaluated at the origin.

\begin{figure}[hbt!]
	\centering{
		\subfloat{
			\begin{minipage}[t]{0.48\textwidth}
				\centering
				 \includegraphics[width=0.95\linewidth]{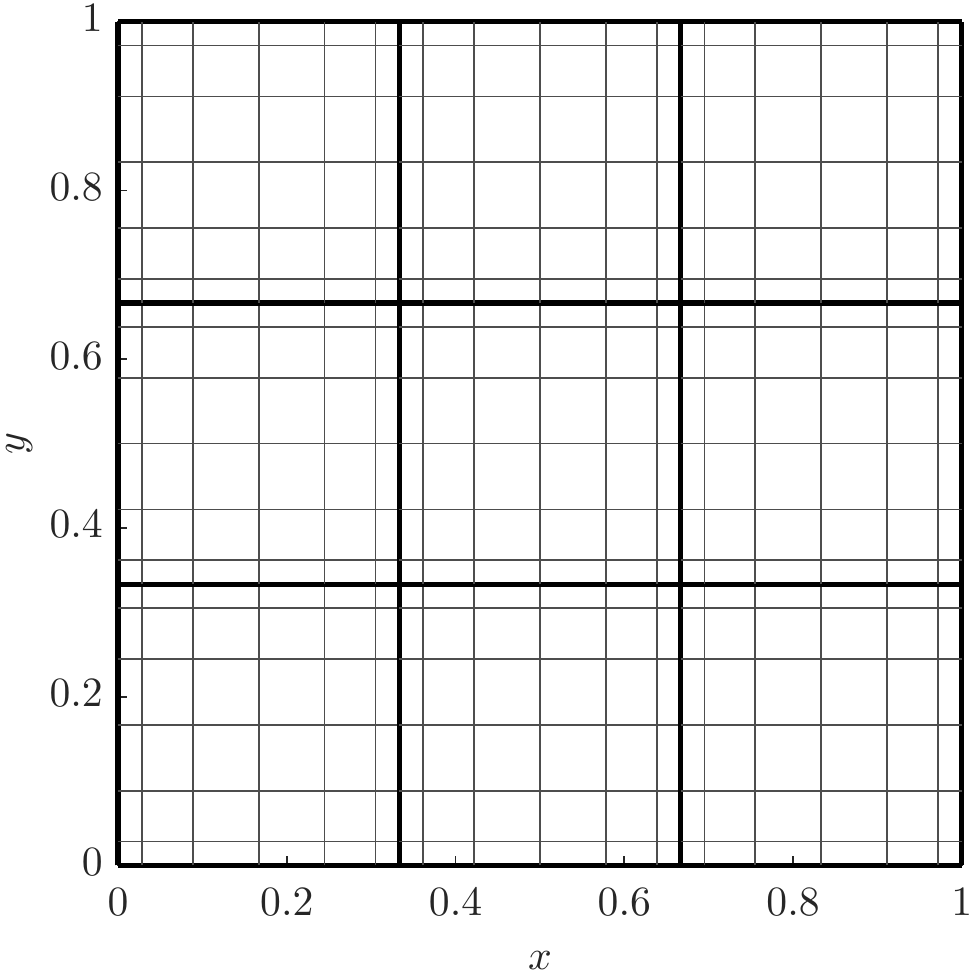}
			\end{minipage}
		}
		\subfloat{
			\begin{minipage}[t]{0.48\textwidth}
				\centering
				 \includegraphics[width=0.95\linewidth]{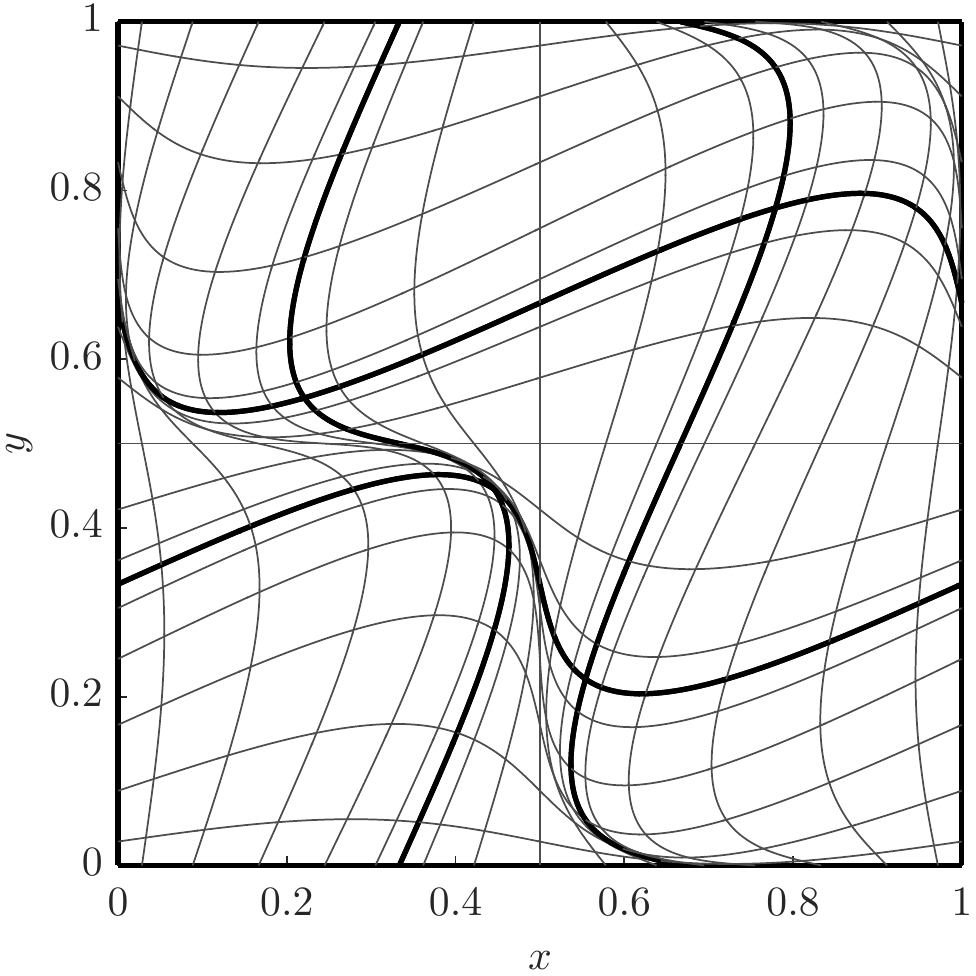}
			\end{minipage}
		}
		\caption{Example meshes with $3\times 3$ elements of polynomial degree $N=6$. \emph{Left:} $c=0$ (orthogonal mesh). \emph{Right:} $c=0.3$ (highly deformed mesh).}
		\label{fig. meshes_testcase_yi}
	}
\end{figure}

The meshes we use here are obtained by deforming the GLL meshes in the reference domain $ \left(  \xi,\eta \right) \in \Omega_{\mathrm{ref}}= [-1,1]^{2}  $ with the mapping, $\Phi$, given as
\begin{equation}
\left\lbrace
\begin{aligned}
x&=\frac{1}{2}+\frac{1}{2}\left( \xi +c\sin(\pi\xi)\sin(\pi\eta)\right) \\
y&=\frac{1}{2}+\frac{1}{2}\left( \eta+c\sin(\pi\xi)\sin(\pi\eta)\right) \\
\end{aligned}\right. ,
\end{equation}
where $ c $ is the deformation coefficient. The two meshes, for $ c=0.0 $ and $ c=0.3 $, are shown in Figure~\ref{fig. meshes_testcase_yi}.

The method is tested for $\alpha\in \lbrace {0,0.01} \rbrace$ and $ c \in \lbrace {0,0.3} \rbrace$.

\begin{figure}[hbt!]
	\centering{
		\subfloat{
			\begin{minipage}[t]{0.48\textwidth}
				\centering
				 \includegraphics[width=1\linewidth]{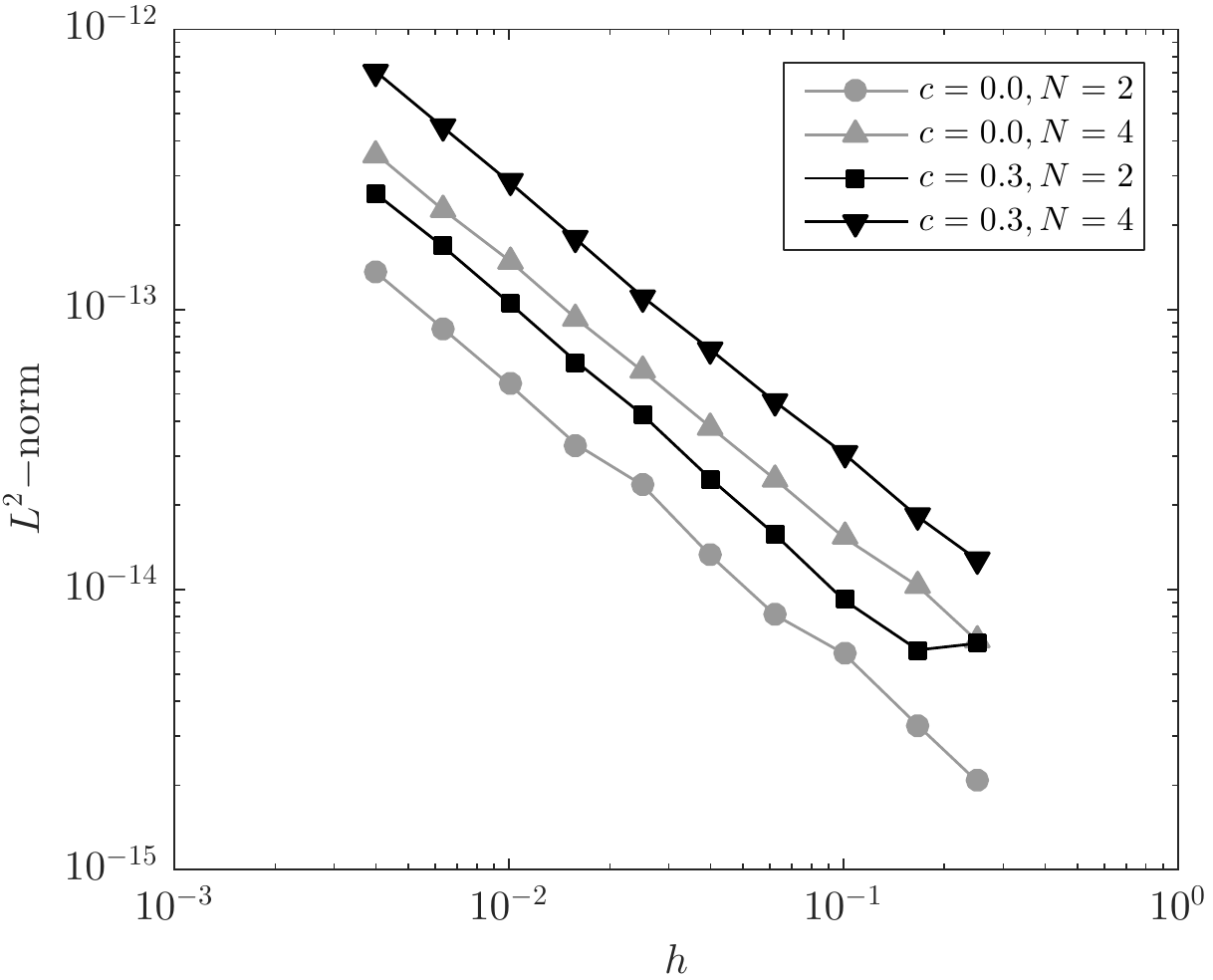}
			\end{minipage}
		}
		\subfloat{
			\begin{minipage}[t]{0.48\textwidth}
				\centering
				 \includegraphics[width=1\linewidth]{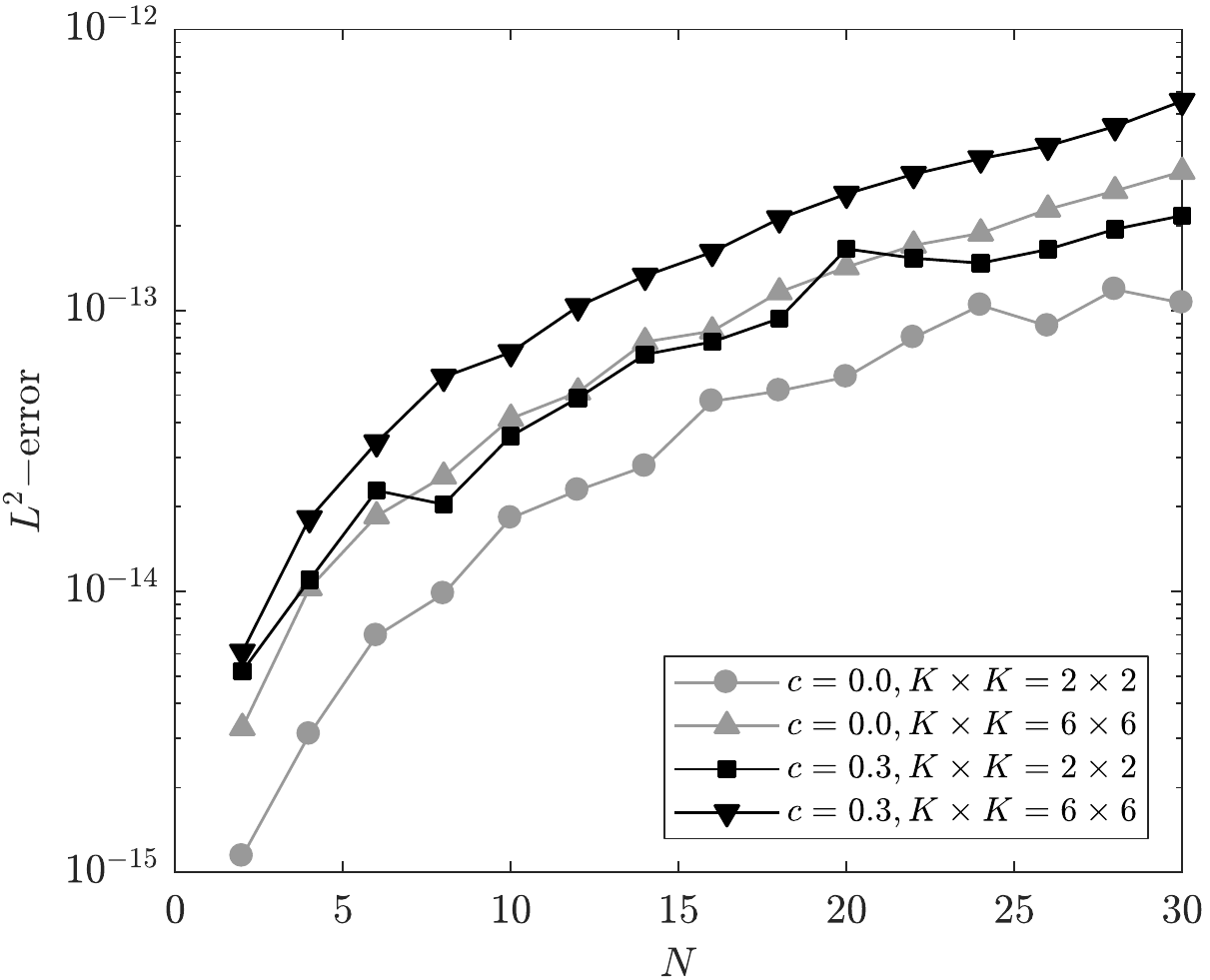}
			\end{minipage}
		}\\
		\subfloat{
		\begin{minipage}[t]{0.48\textwidth}
			\centering
			 \includegraphics[width=1\linewidth]{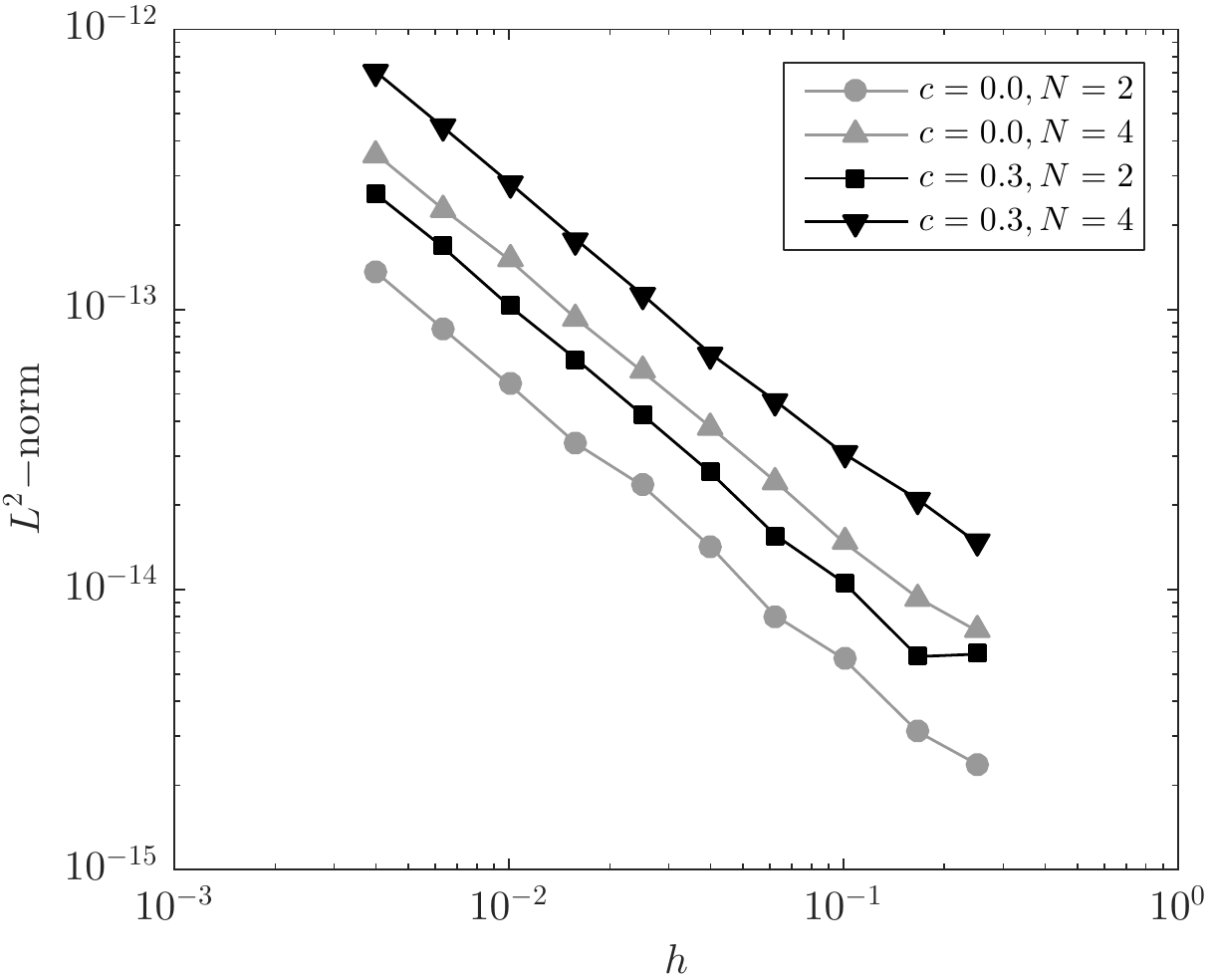}
		\end{minipage}
		}
		\subfloat{
		\begin{minipage}[t]{0.48\textwidth}
			\centering
			 \includegraphics[width=1\linewidth]{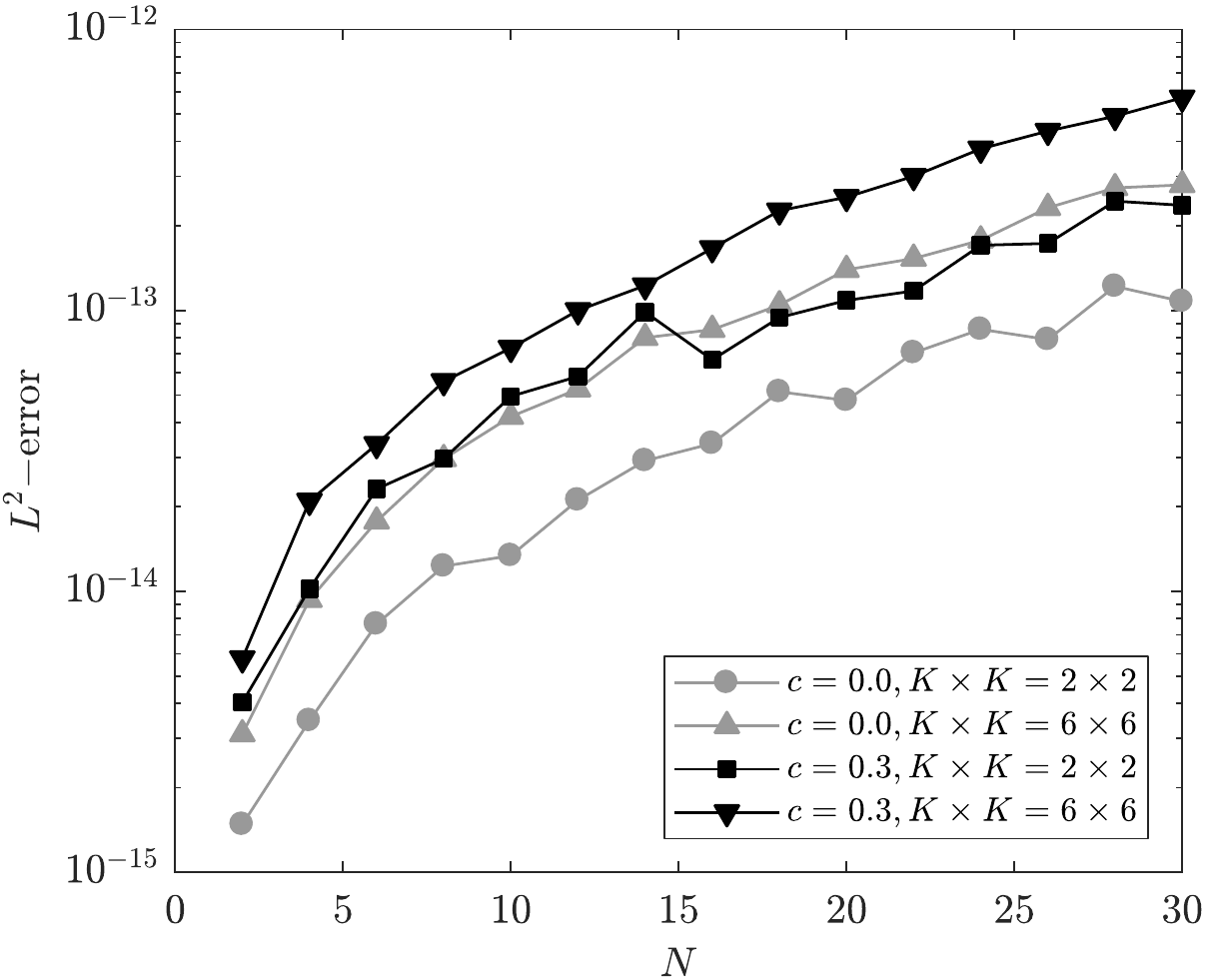}
		\end{minipage}
		}
		\caption{The $ L^{2} $-norm of $ \left( \nabla\cdot\boldsymbol{u}_{h}-f_{h}\right) $.  \emph{Left:} $ K\times K $ elements, $ K=4,...,250 $, and $ N=2,4 $. \emph{Right:} $ 2\times2, 6\times6 $ elements, and $ N = 2,...,30 $. \emph{Top:} $ \alpha=0 $. \emph{Bottom:} $ \alpha=0.01 $.}
		\label{fig. du_f}
	}
\end{figure}

\begin{figure}[hbt!]
	\centering
	\subfloat{
		\begin{minipage}[t]{0.48\textwidth}
			\centering
			 \includegraphics[width=1\linewidth]{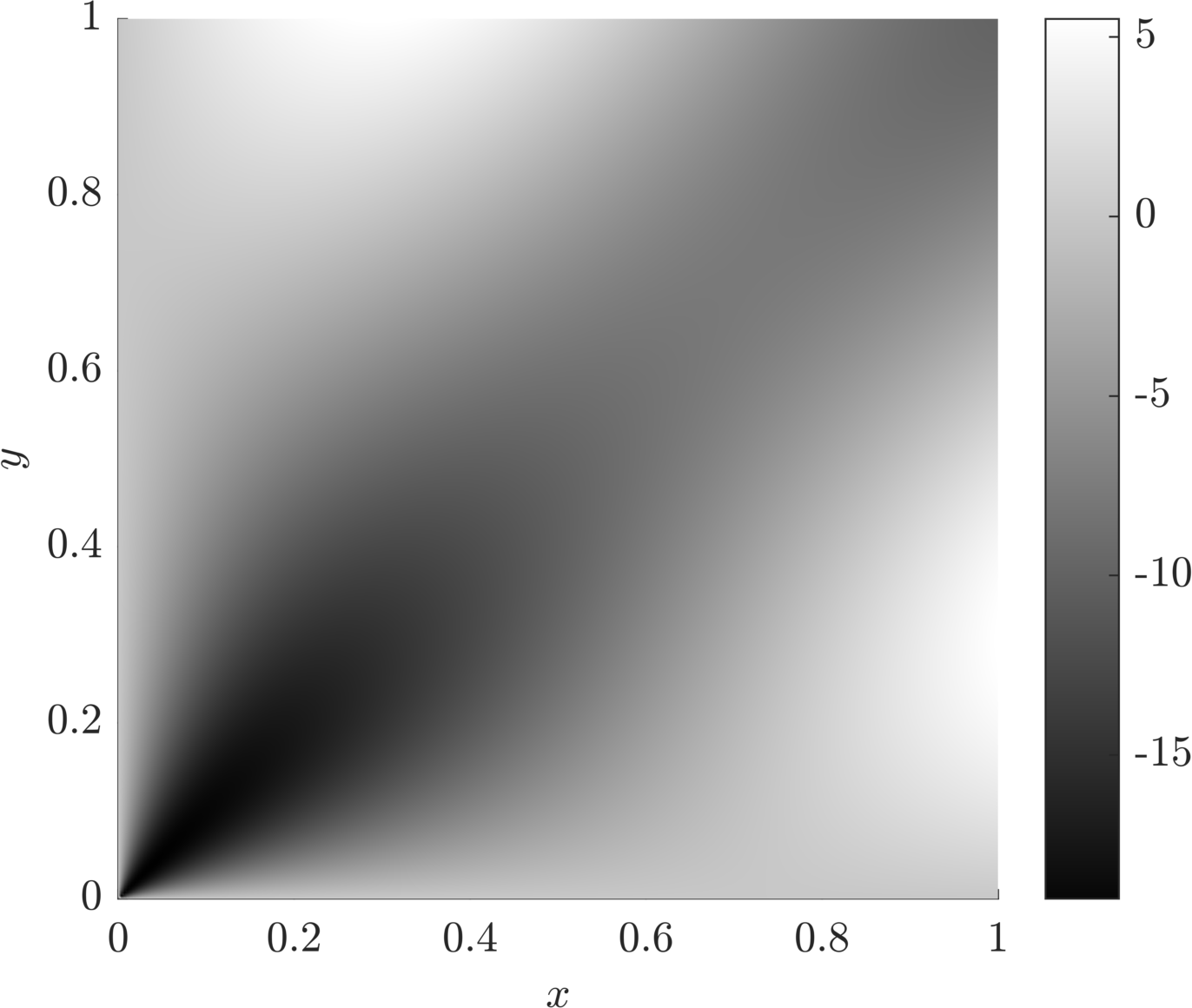}
		\end{minipage}
	}
	\subfloat{
		\begin{minipage}[t]{0.466\textwidth}
			\centering
			 \includegraphics[width=1\linewidth]{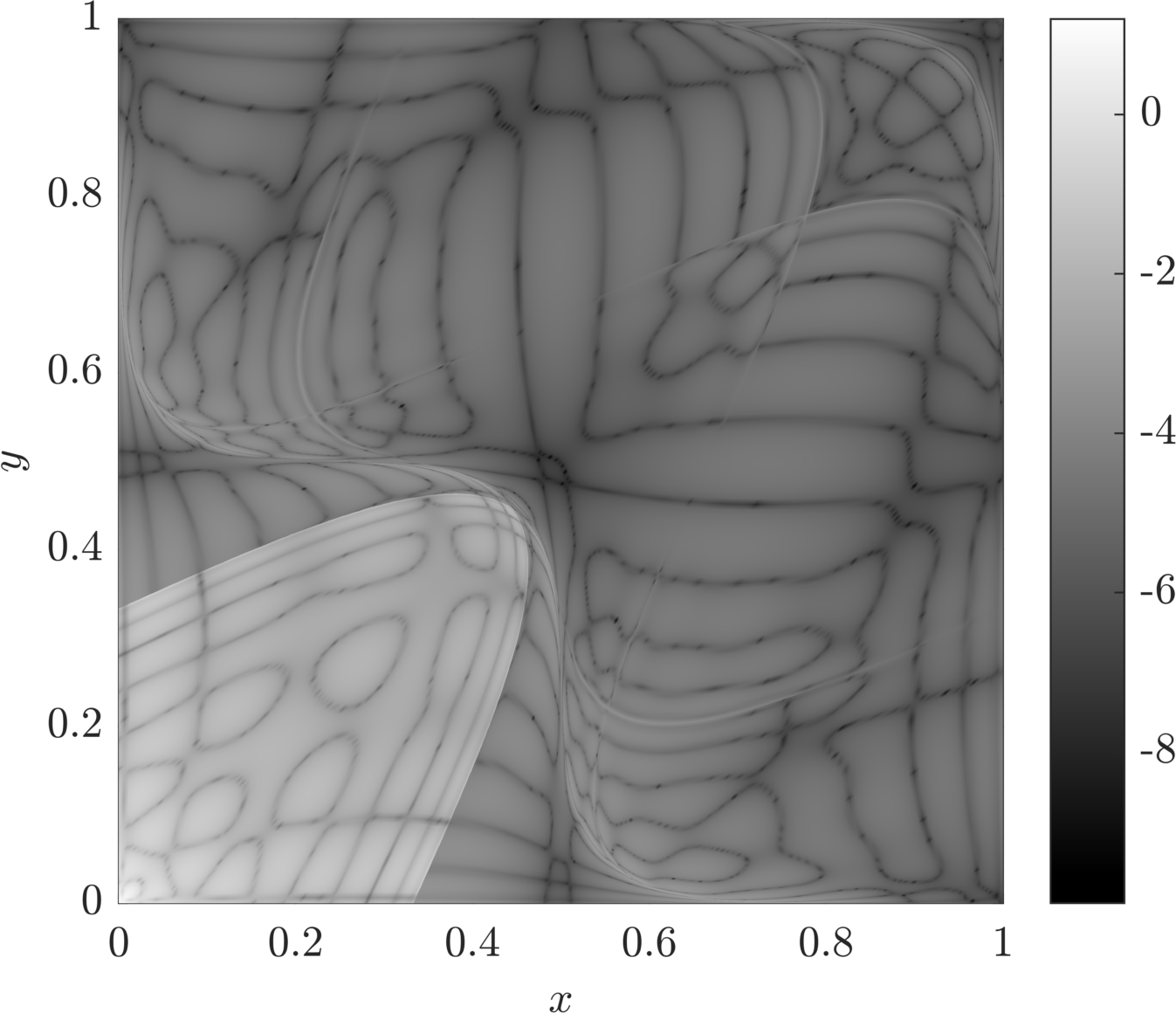}
		\end{minipage}
	}\\
\subfloat{
\begin{minipage}[t]{0.48\textwidth}
	\centering
	 \includegraphics[width=1\linewidth]{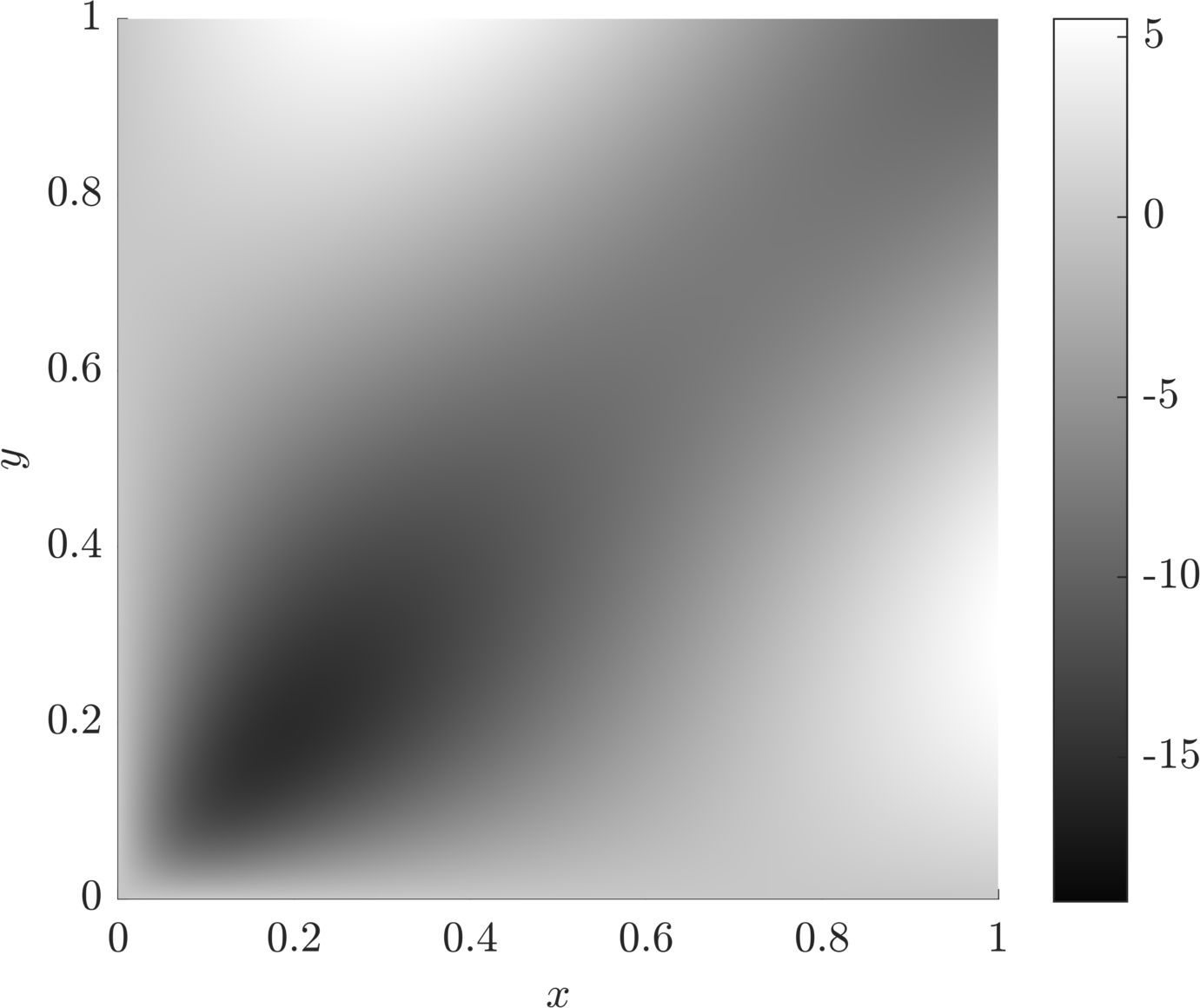}
\end{minipage}
}
\subfloat{
\begin{minipage}[t]{0.466\textwidth}
	\centering
	 \includegraphics[width=1\linewidth]{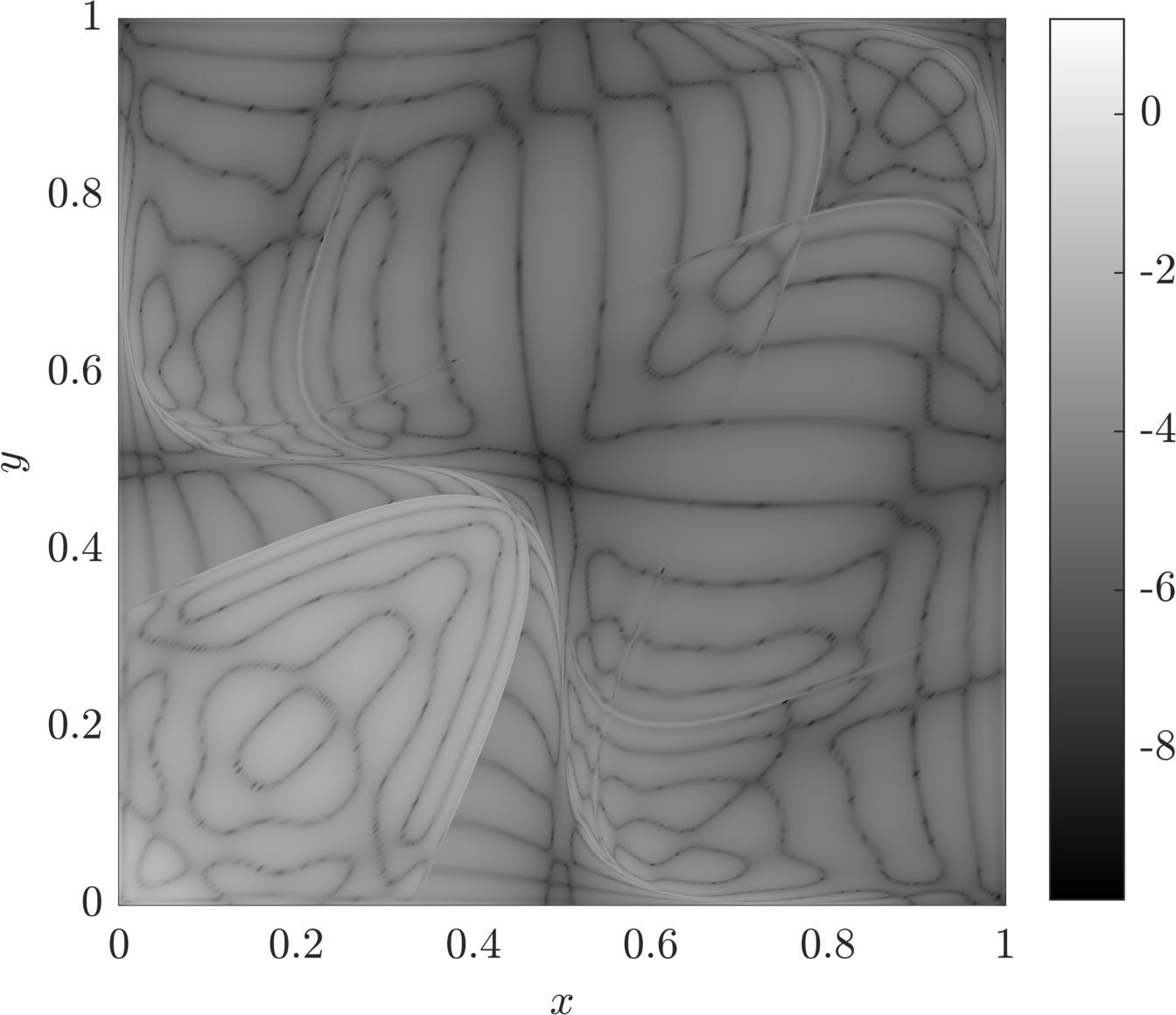}
\end{minipage}
}
	\caption{\emph{Left:} the source term $ f $. \emph{Right:} the $ \log_{10} $ distribution of the projection error of $ f_{h} $ for $ 3\times3 $ elements, $ N=10 $ and $ c=0.3 $. \emph{Top:} $ \alpha=0 $. \emph{Bottom:} $ \alpha=0.01 $.}
	\label{fig:ffield}
\end{figure}

\begin{figure}[hbt!]
	\centering{
		\subfloat{
			\begin{minipage}[t]{0.48\textwidth}
				\centering
				 \includegraphics[width=1\linewidth]{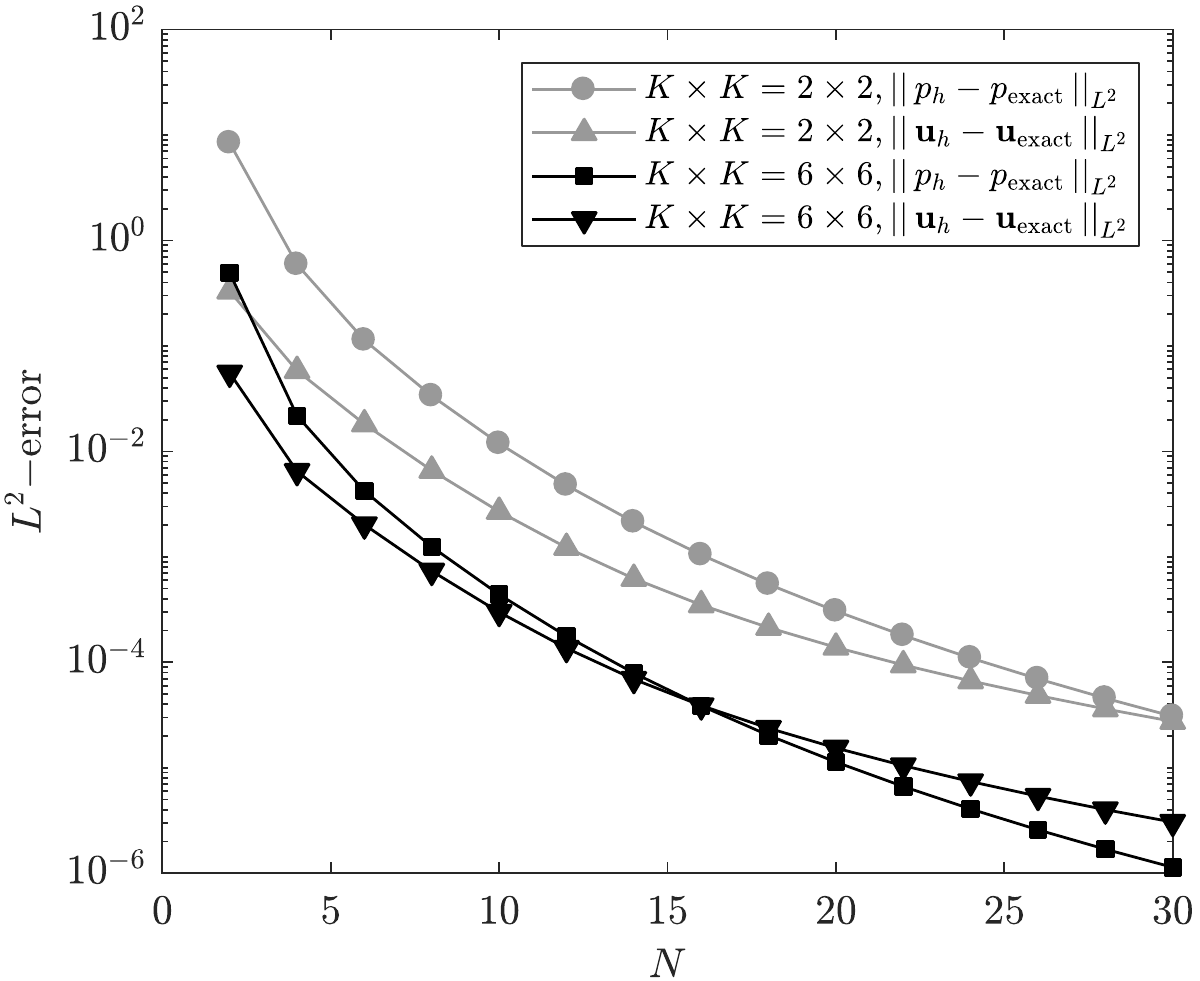}
			\end{minipage}
		}
		\subfloat{
			\begin{minipage}[t]{0.48\textwidth}
				\centering
				 \includegraphics[width=1\linewidth]{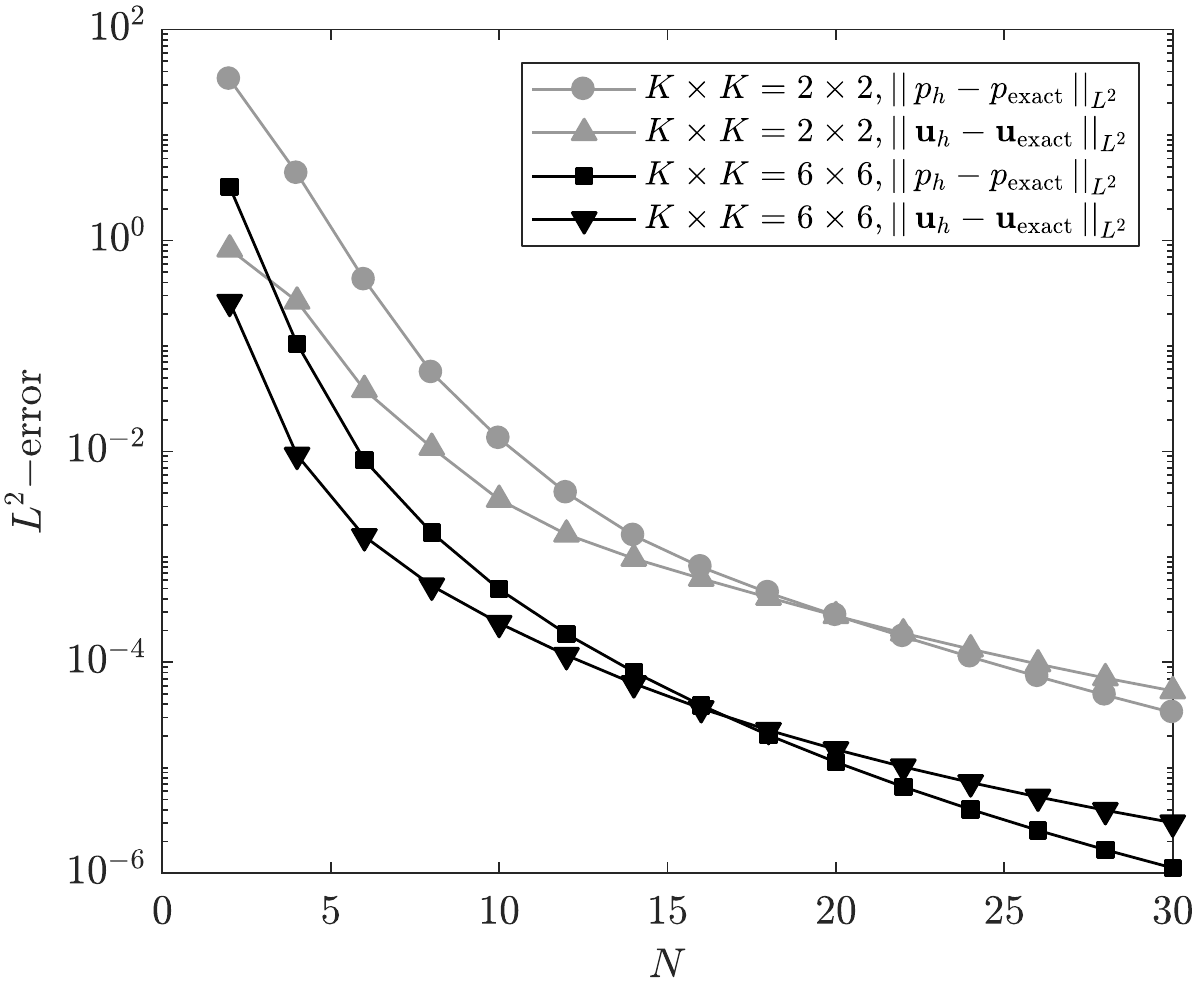}
			\end{minipage}
		}\\
		\subfloat{
		\begin{minipage}[t]{0.48\textwidth}
			\centering
			 \includegraphics[width=1\linewidth]{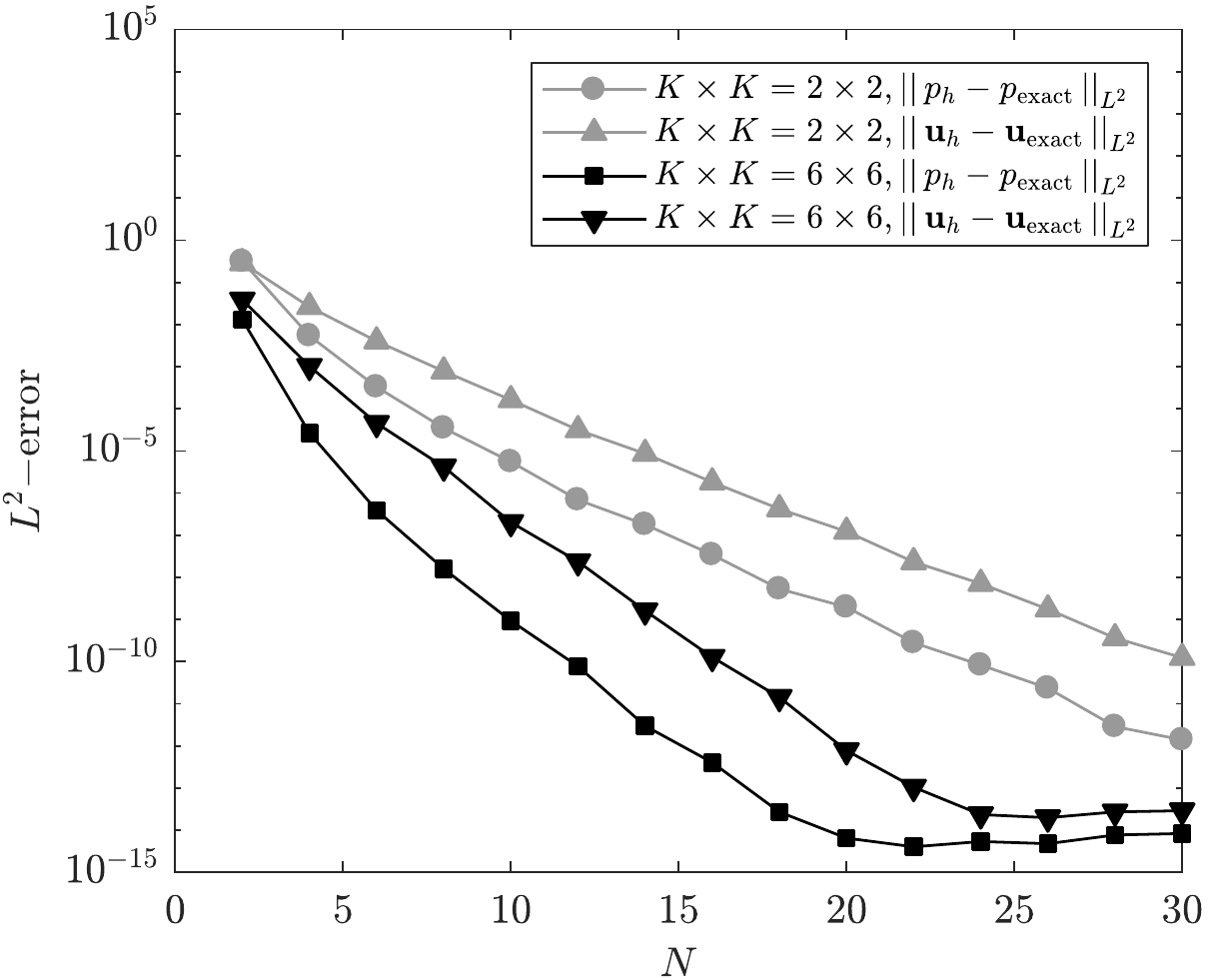}
		\end{minipage}
		}
		\subfloat{
		\begin{minipage}[t]{0.48\textwidth}
			\centering
			 \includegraphics[width=1\linewidth]{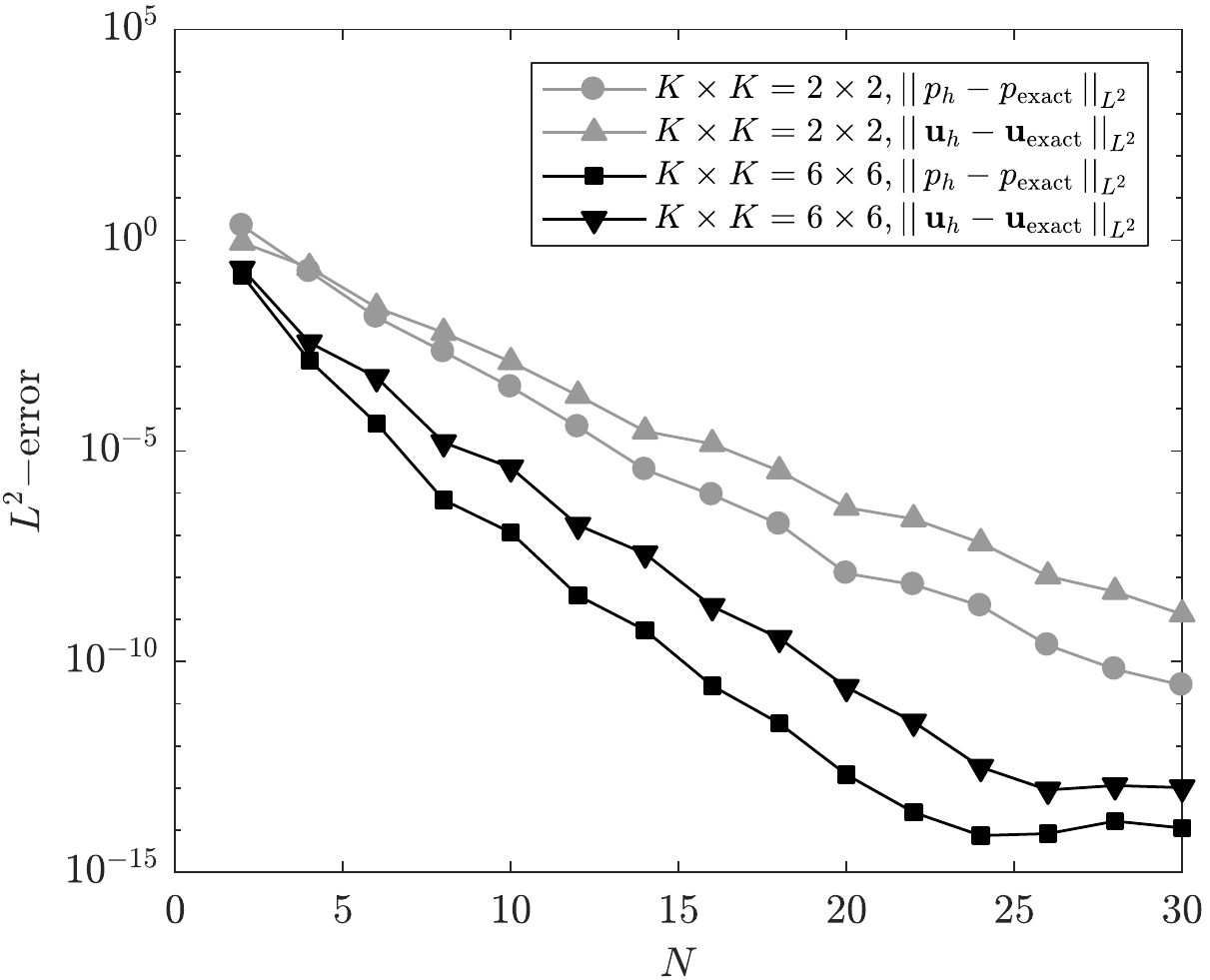}
		\end{minipage}
		}
		\caption{The $ p $-convergence for $ 2\times2, 6\times6 $ elements and $ N=2,...,30 $. \emph{Left:} $ c=0 $. \emph{Right:} $ c=0.3 $. \emph{Top:} $ \alpha=0 $. \emph{Bottom:} $ \alpha=0.01 $.}
		\label{p-convergence}
	}
\end{figure}

\begin{figure}[hbt!]
	\centering{
		\subfloat{
			\begin{minipage}[t]{0.48\textwidth}
				\label{fig. h_convergence_c00}
				\centering
				 \includegraphics[width=1\linewidth]{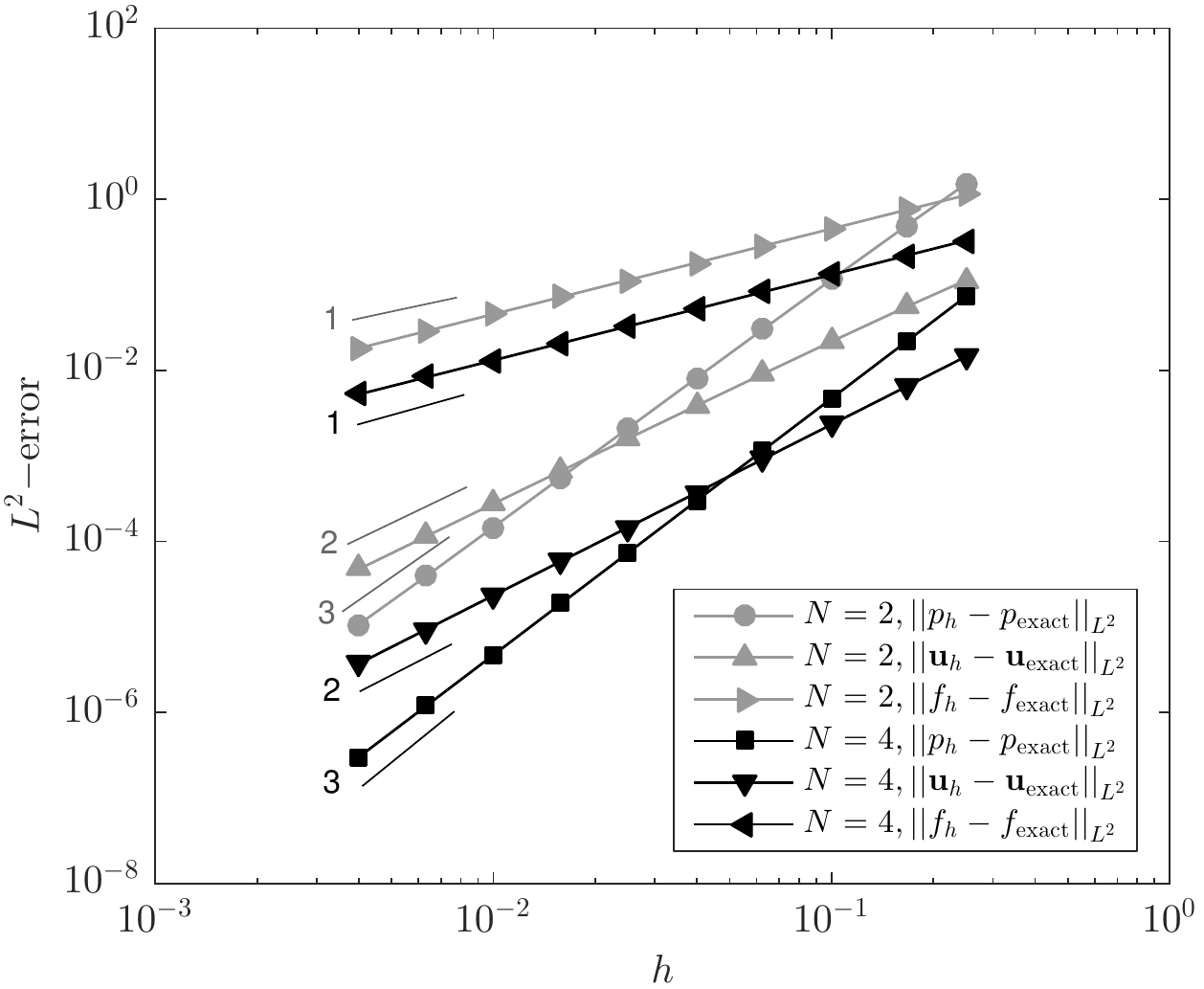}
			\end{minipage}
		}
		\subfloat{
			\begin{minipage}[t]{0.48\textwidth}
				\label{fig. h_convergence_c03}
				\centering
				 \includegraphics[width=1\linewidth]{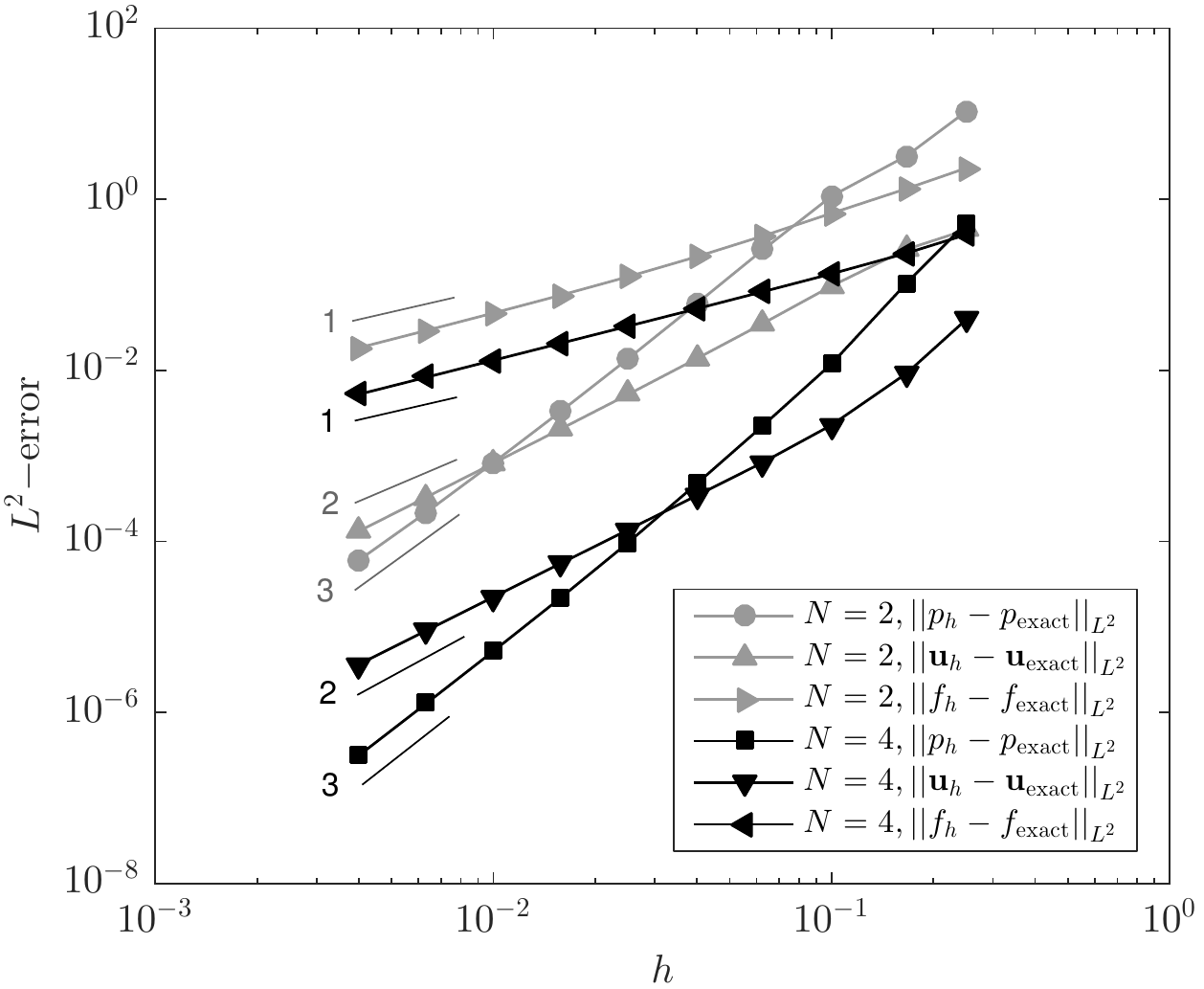}
			\end{minipage}
		}\\
		
		\subfloat{
			\begin{minipage}[t]{0.48\textwidth}
				\label{fig. a001_h_convergence_c00}
				\centering
				 \includegraphics[width=1\linewidth]{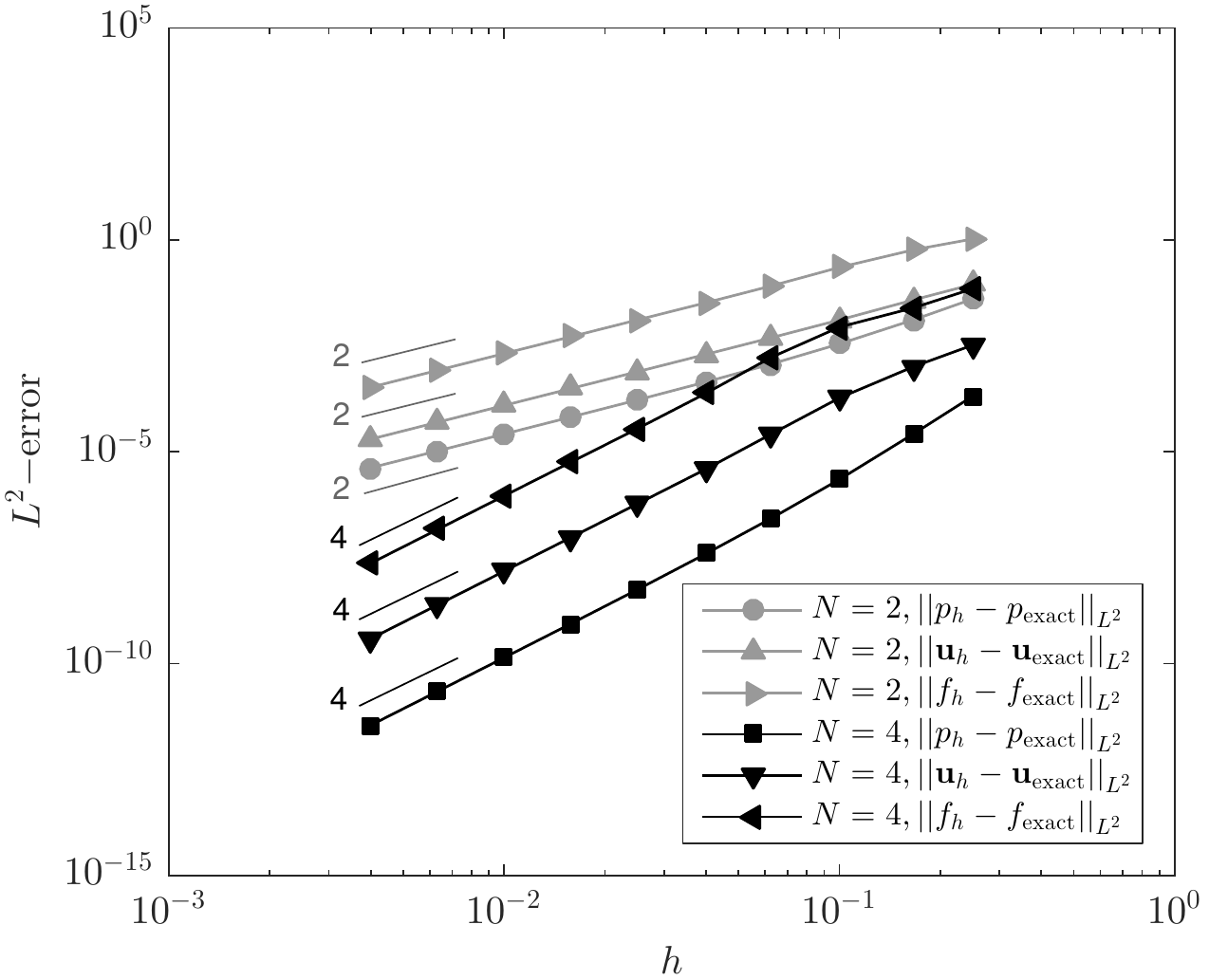}
			\end{minipage}
		}
		\subfloat{
			\begin{minipage}[t]{0.48\textwidth}
				\label{fig. a001_h_convergence_c03}
				\centering
				 \includegraphics[width=1\linewidth]{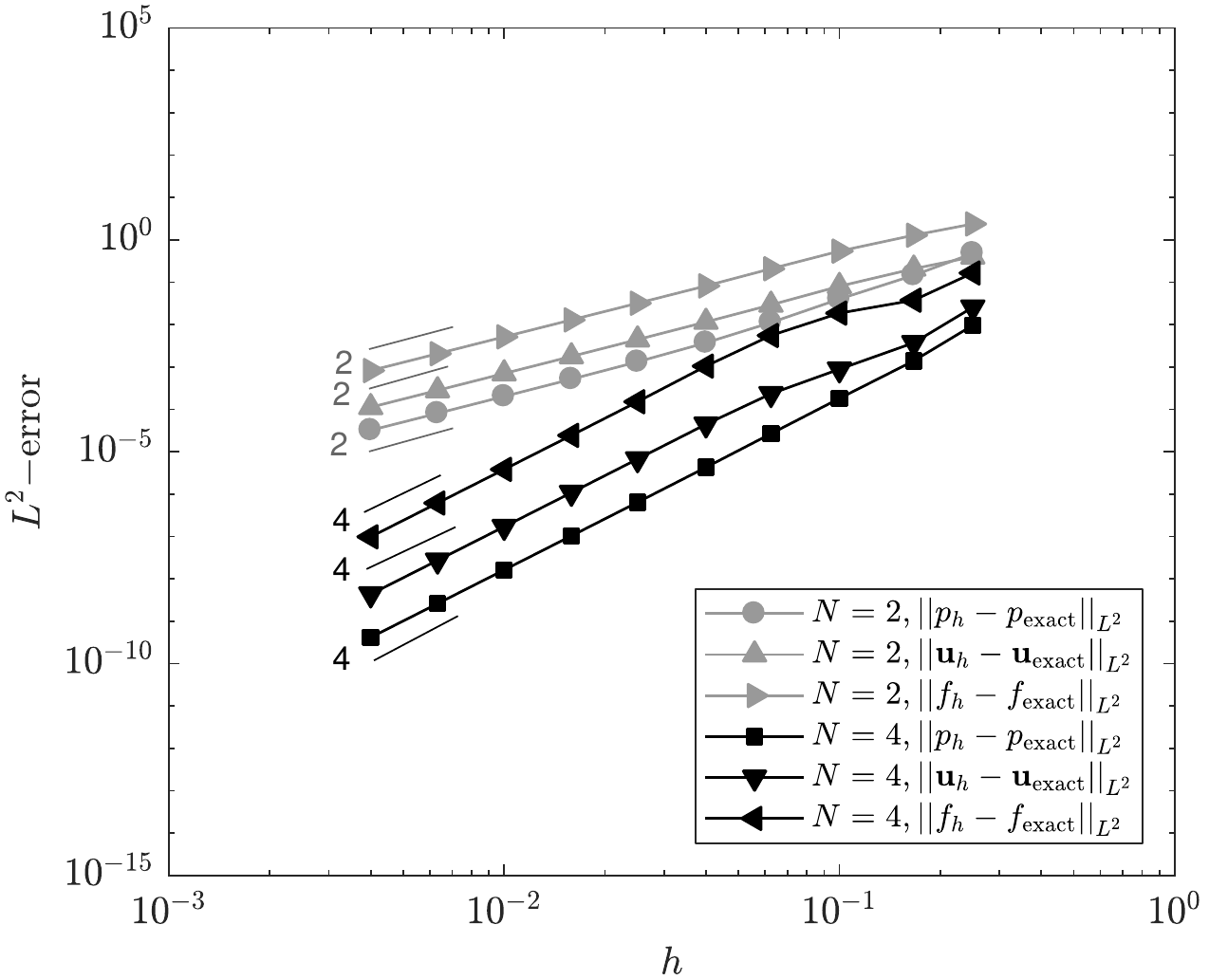}
			\end{minipage}
		}
		\caption{The $ h $-convergence of the $ L^{2} $-error for $ K\times K $ elements, $ K=4,...,250 $ and $ N=2,4 $. \emph{Left:} $ c=0 $. \emph{Right:} $ c=0.3 $. \emph{Top:} $ \alpha=0 $. \emph{Bottom:} $ \alpha=0.01 $.}
		\label{h-convergence}
	}
\end{figure}

In Figure~\ref{fig. du_f}, the results for $ \left| \left|  \nabla\cdot\boldsymbol{u}_{h}-f_{h}\right| \right|_{L^{2}}  $ are presented. They show that the relation $ \nabla\cdot\boldsymbol{u}_{h}=f_{h} $ is conserved to machine precision even on a highly deformed and coarse mesh i.e. of $ 2\times2 $ elements with $ N=2 $ and $ c=0.3 $.

When $ \alpha=0.01 $, $ \mathbb{K} $ is no longer multi-valued at the origin. In this case the source term $ f $ is smooth over the domain, see Figure \ref{fig:ffield} (bottom). For this smooth case, the method displays optimal convergence rates on both the orthogonal mesh and the deformed mesh, i.e. see Figure \ref{p-convergence} (bottom) and Figure \ref{h-convergence} (bottom).

When $ \alpha=0 $, both the $ h $-convergence rate and $ p $-convergence rates are sub-optimal, see Figure \ref{p-convergence} (top) and Figure \ref{h-convergence} (top). This is because $ \mathbb{K} $ is multi-valued and therefore $ f $ becomes singular at the origin when $ \alpha=0 $, see Figure \ref{fig:ffield} (top left).
\clearpage

\subsection{The Sand-Shale system} \label{sec:SS}
\begin{figure}[!htb]
\centering
\includegraphics[scale=0.5]{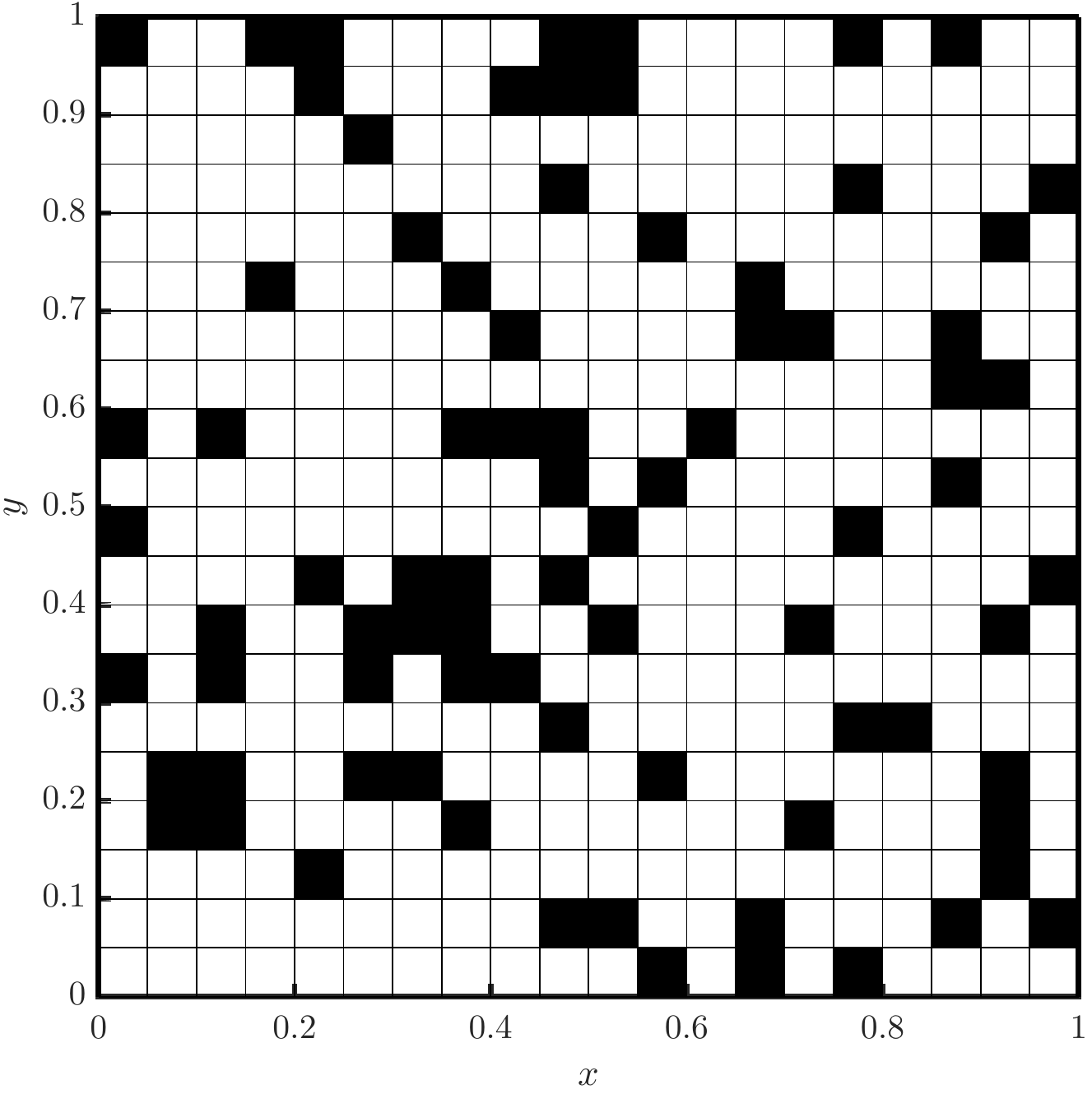}
\caption{The discretized domain for the \VARUN{sand-shale} test case. \emph{Black} blocks are shale blocks with $ k = 10^{-6}$. \emph{White} blocks are sand blocks with $ k = 1 $.}
\label{fig:SSGeometry}
\end{figure}

This example is taken from \cite{Durlofsky1994,HymanShashkovSteinberg97,Kikinzon2017}.
The domain is a \VARUN{2D} unit square, $\Omega = [0,1]^2$, with 80 shale blocks, $\Omega _s$, placed in the domain such that the total area fraction of shale blocks is $A_{shale} = 20 \%$, \VARUN{as shown in Figure~\ref{fig:SSGeometry}}.

We solve the mixed formulation (\ref{eq:finite_element_mixed_formulation}) with $f = 0$ in this domain.
The flux across the top and the bottom boundaries is $\boldsymbol{u} \cdot \boldsymbol{n} = 0$.
The flow is pressure driven with the pressure at the left boundary, $p=1$, and the pressure at the right boundary, $p=0$.
The permeability in the domain is defined as $\mathbb{K}=k\mathbb{I}$, where $k$ is given by:
\begin{equation*}
k = \left\lbrace
\begin{aligned}
&10^{-6} & & \quad \mathrm{in} \quad & & \Omega _s\\
&1 		 & & \quad \mathrm{in} \quad & & \Omega \setminus \Omega _s
\end{aligned}\right. \;.
\end{equation*}

\begin{figure}[!htb]
	\centering
	 \includegraphics[scale=0.5]{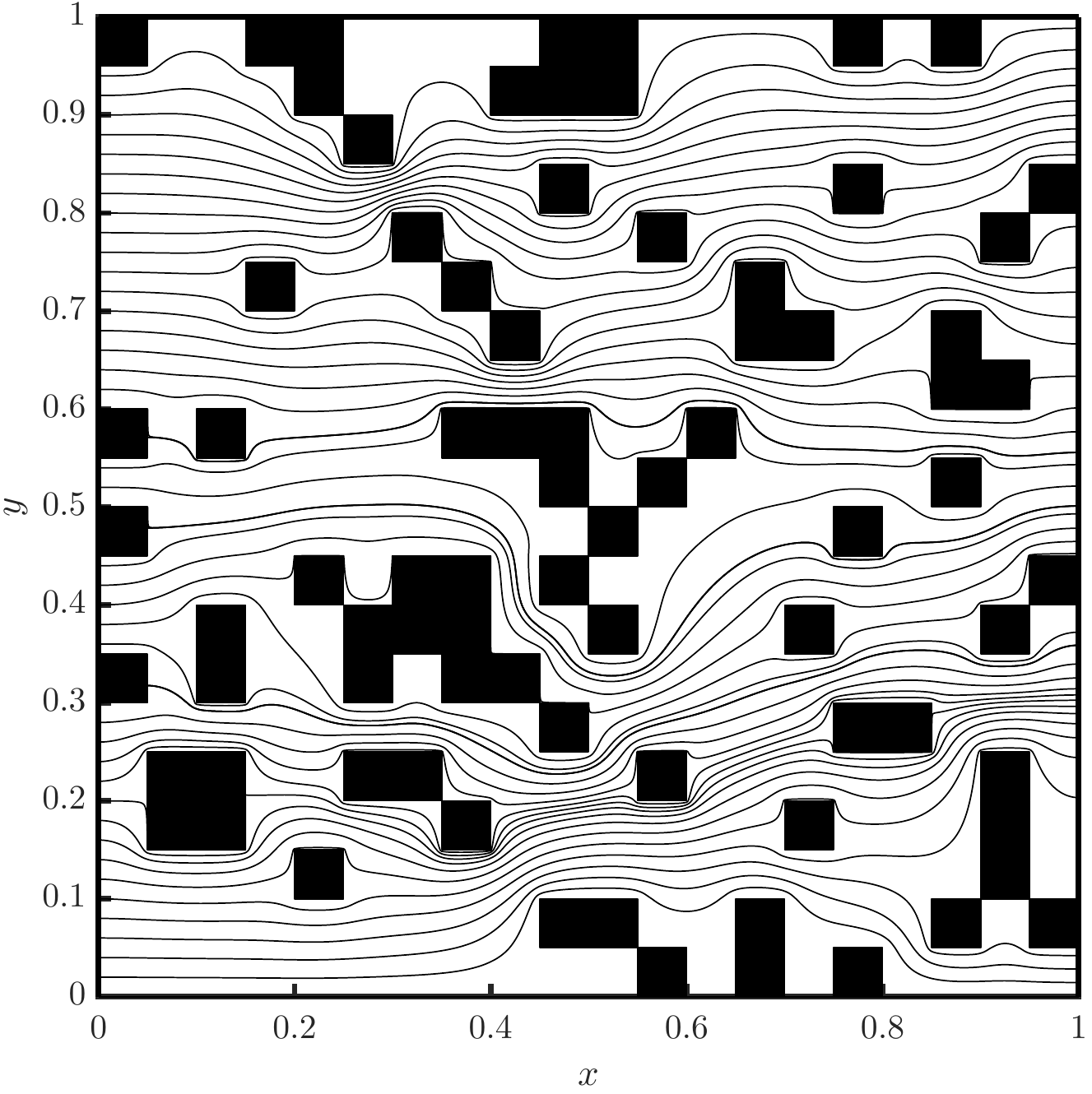}
	\caption{Streamlines through the domain of \VARUN{sand-shale} test case.}
	\label{fig:SSStream}
\end{figure}

For this test case an orthogonal uniform grid of $20\times20$ elements is used.
The polynomial degree is varied to achieve convergence.
Streamlines through the domain for $20 \times 20$ elements and polynomial degree $N=15$ are shown in Figure~\ref{fig:SSStream}.
It can be seen that the streamlines do not pass through, but pass around the shale blocks of low permeability.
\begin{figure}[!htb]
\centering
\includegraphics[scale=0.4]{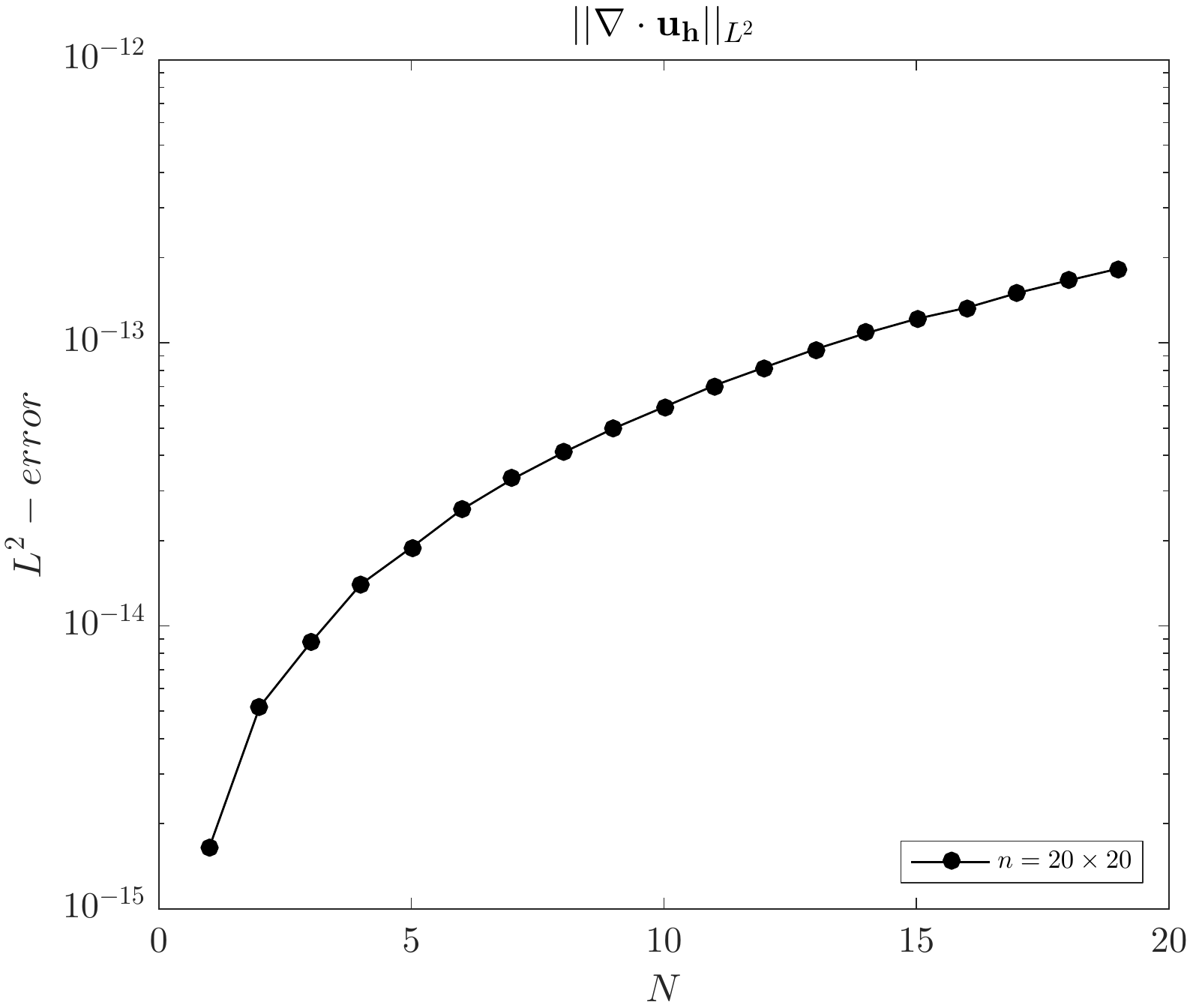}
\caption{ The $L^2$-norm of \VARUN{$\nabla \cdot \boldsymbol{u} _h$} for $20\times20$ elements \VARUN{for a polynomial approximation of} $N=1,...,19$.}
\label{fig:SSDiv}
\end{figure}

The $||\nabla\cdot\boldsymbol{u}_h||_{L^2} $ over the entire domain as a function of polynomial degree is shown in Figure~\ref{fig:SSDiv}.
We observe that $\nabla \cdot \boldsymbol{u}_h = 0$ is satisfied up to machine precision.
\begin{table}[!htb]
\centering
\begin{tabular}{p{1.2cm}  c c}
\hline
N 			& Net flux 		 & No. of unknowns\\
\hline
1 					& 0.49041	& 1240	 \\
2 					& 0.51247	& 4880	 \\
3					& 0.51744	& 10920	 \\
4					& 0.51863	& 19360	 \\
5					& 0.51931	& 30200	 \\
6					& 0.51957	& 43440	 \\
7					& 0.51977	& 59080	 \\
8					& 0.51985	& 77120	 \\
9					& 0.51993	& 97560	 \\
10					& 0.51997	& 120400	 \\
11					& 0.52001	& 145640 \\
12 					& 0.52003	& 173280 \\
13					& 0.52005	& 203320 \\
14					& 0.52007	& 235760	 \\
15					& 0.52008	& 270600	 \\
16					& 0.52009	& 307840	 \\
17					& 0.52009	& 347480	 \\
18					& 0.52010	& 389520	 \\
19					& 0.52010	& 433960	 \\
\hline
\end{tabular}
\caption{Net flux through the left boundary of the \VARUN{sand-shale} domain for $k=10^{-6}$, $20\times 20$ elements, $N=1,...,19$.}
\label{tab:SSError}
\end{table}

The net flux entering the domain (the same as the net flux leaving the domain) is given in Table~\ref{tab:SSError} for varying polynomial degree.
A reference value for this solution is given in \cite{Durlofsky1994} as $0.5205$, \VARUN{and in \cite{Kikinzon2017} as $0.519269$}.
In this work the maximum resolution corresponds to $20 \times 20$ elements and a polynomial degree $N=19$, for which the net flux entering the domain is obtained as $0.52010$.
\begin{figure}[hbt!]
	\centering{
		\subfloat{
			\begin{minipage}[t]{0.48\textwidth}
				\label{fig:k=1}
				\centering
				 \includegraphics[width=1\linewidth]{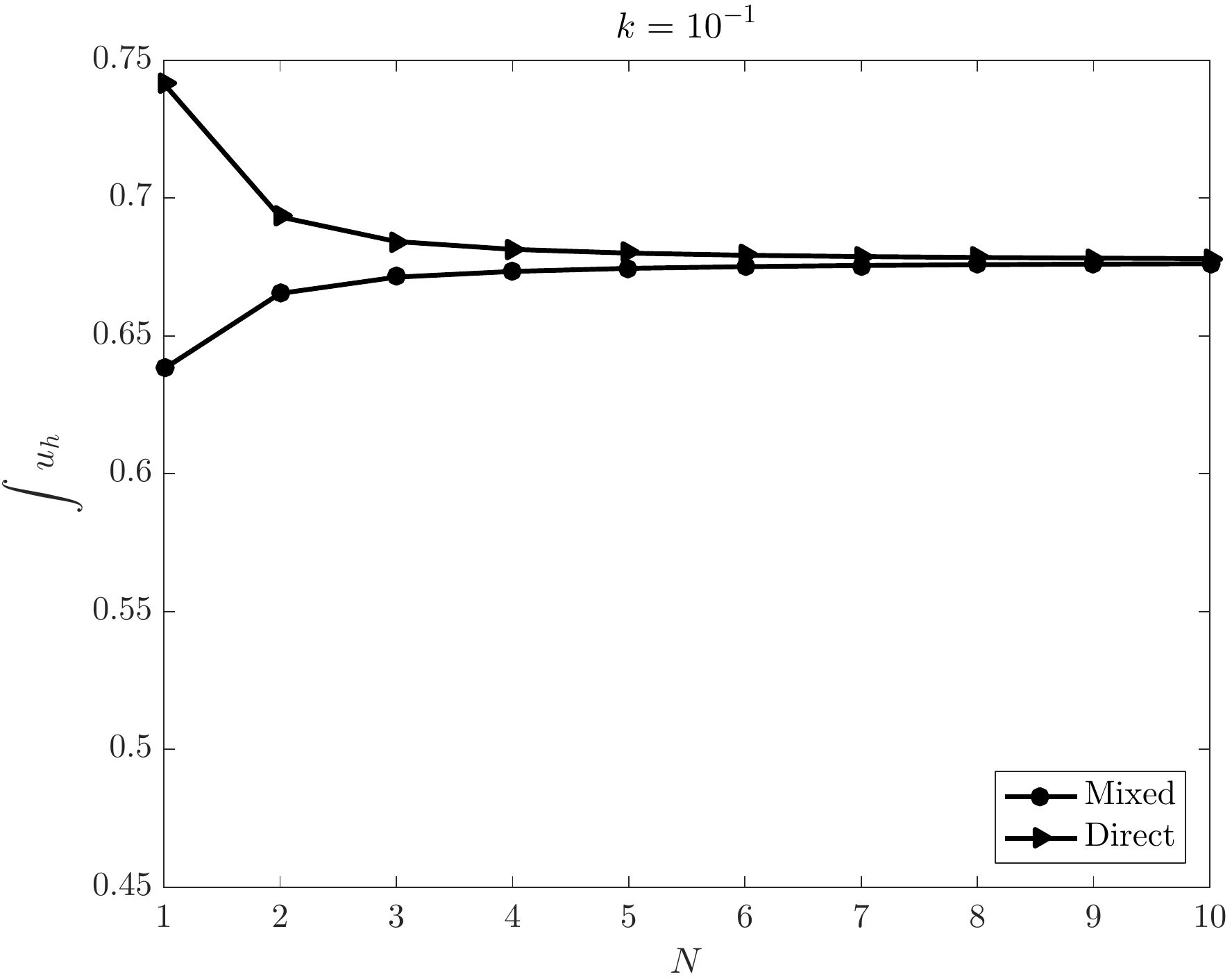}
			\end{minipage}
		}
		\subfloat{
			\begin{minipage}[t]{0.48\textwidth}
				\label{fig:k=2}
				\centering
				 \includegraphics[width=1\linewidth]{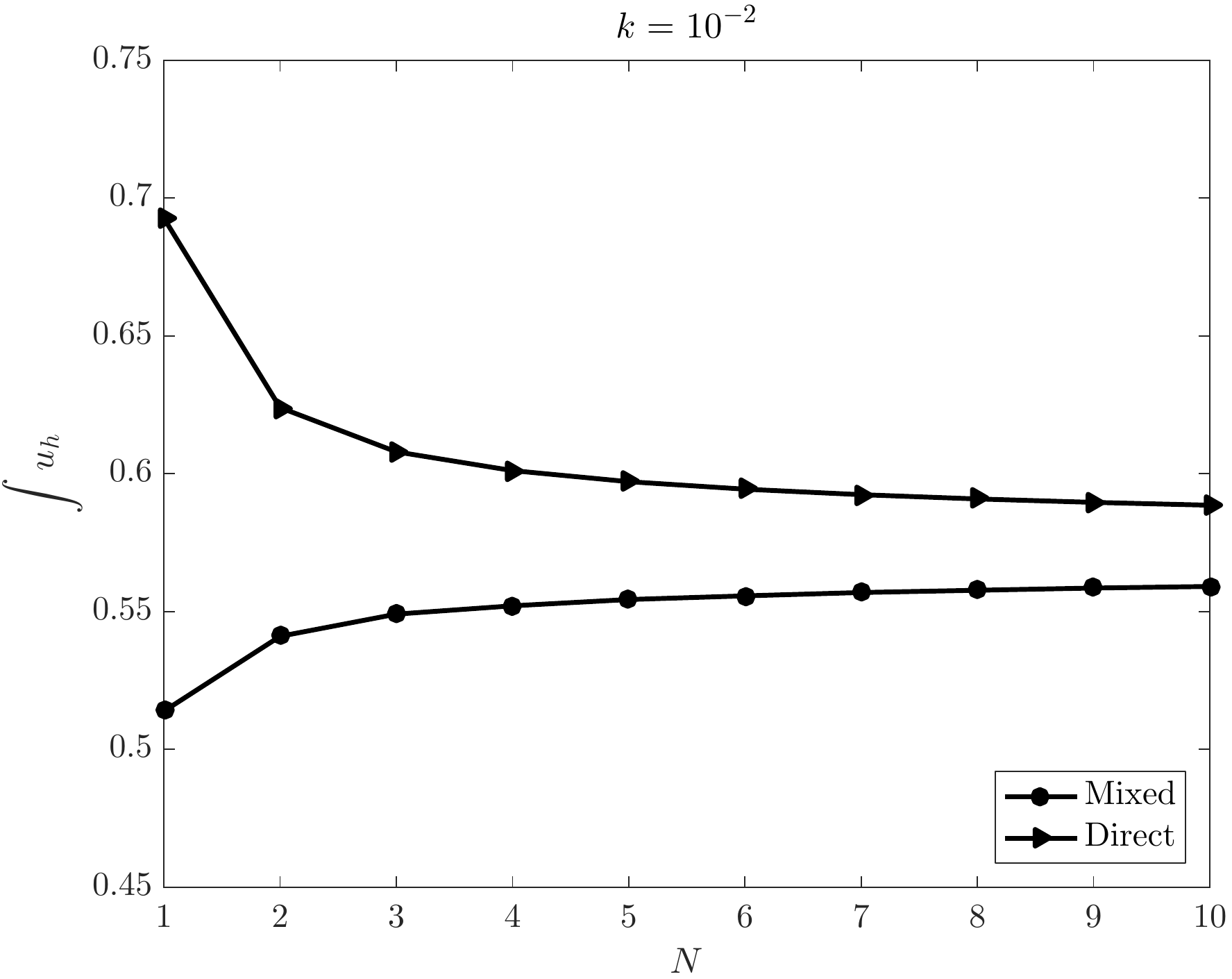}
			\end{minipage}
		}\\
		\subfloat{
			\begin{minipage}[t]{0.48\textwidth}
				\label{fig:k=3}
				\centering
				 \includegraphics[width=1\linewidth]{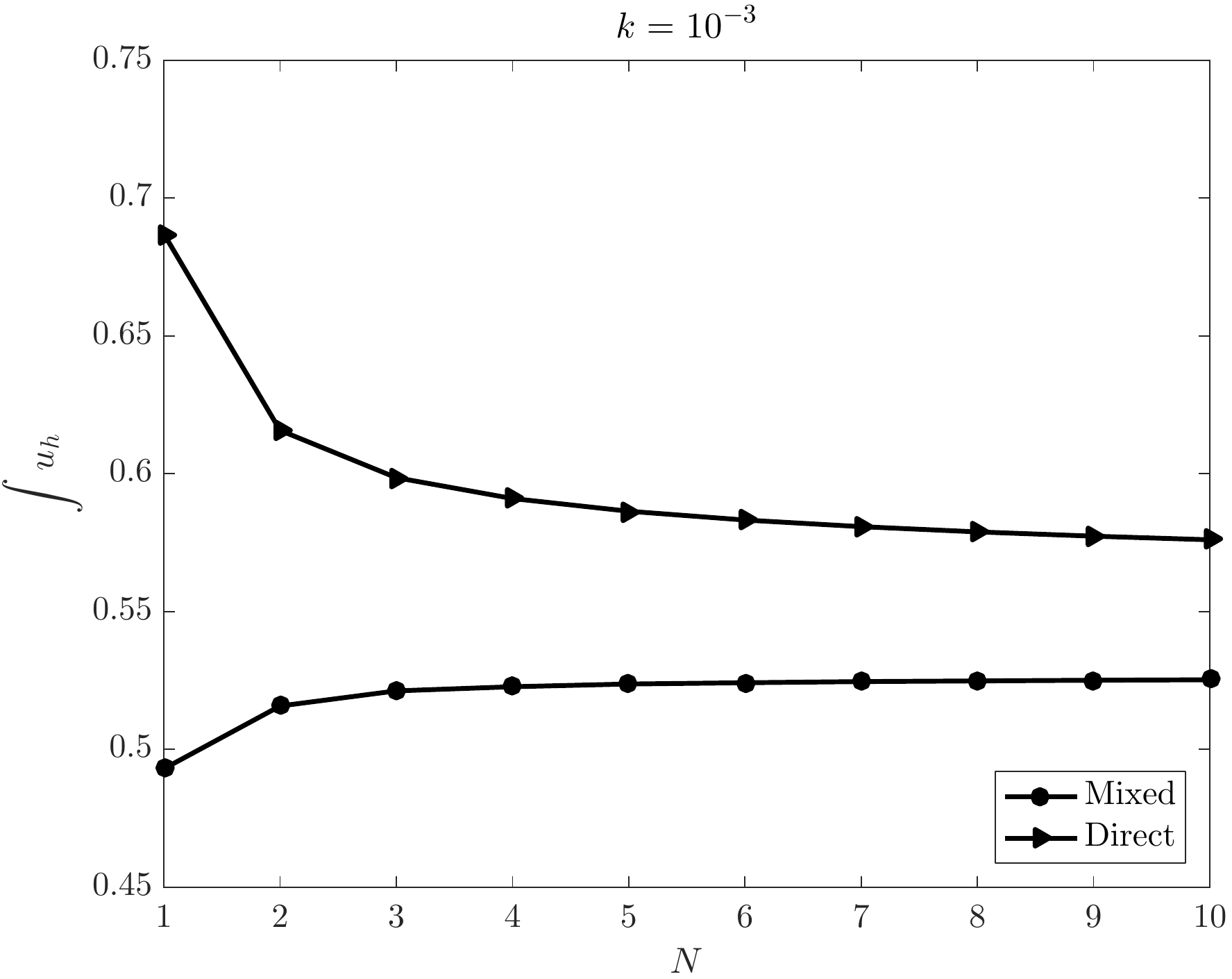}
			\end{minipage}
		}
		\subfloat{
			\begin{minipage}[t]{0.48\textwidth}
				\label{fig:k=4}
				\centering
				 \includegraphics[width=1\linewidth]{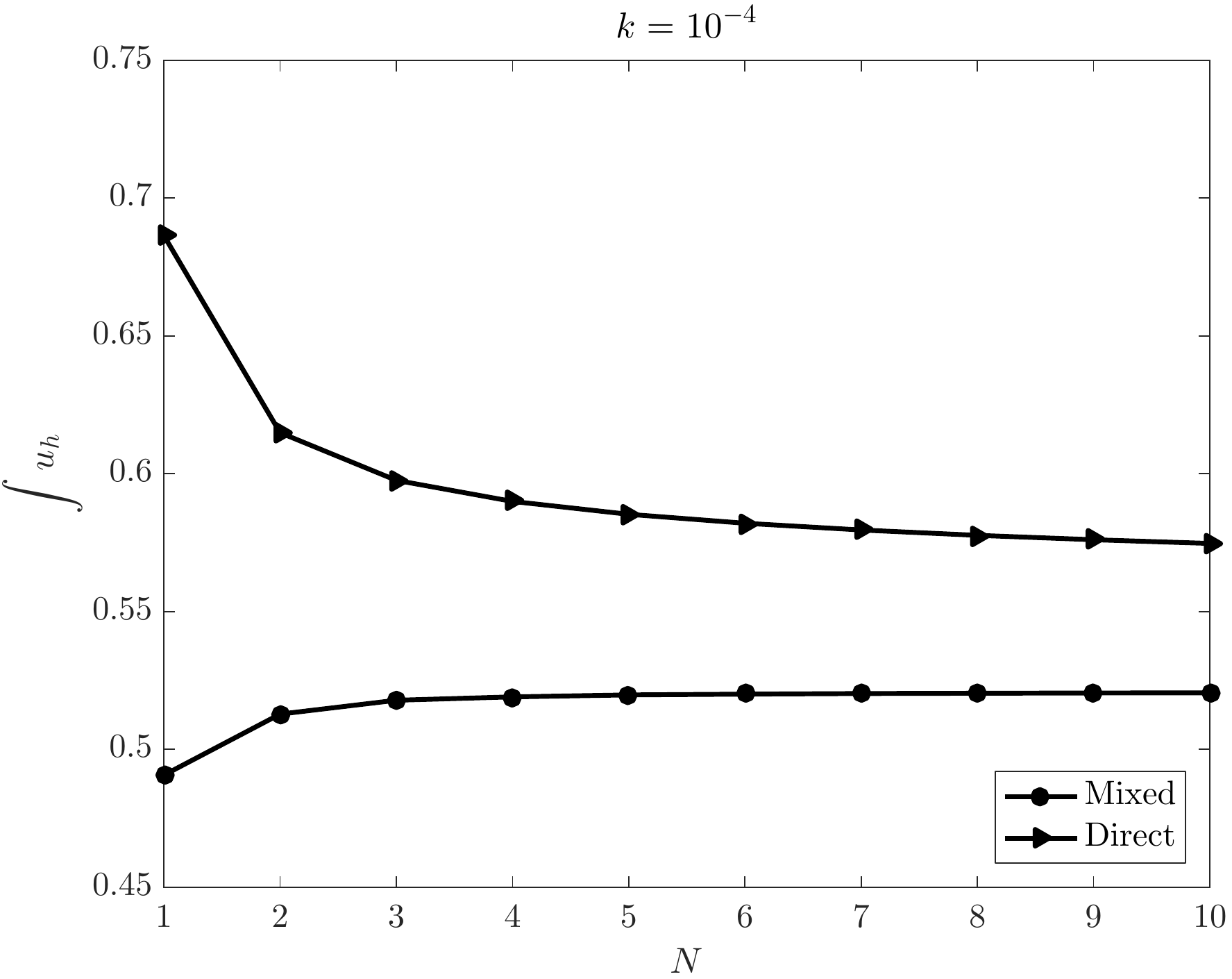}
			\end{minipage}
		}
		\caption{Convergence of the net flux through the left boundary of the \VARUN{sand-shale} domain using the mixed formulation and the direct formulation for $20\times20$ elements, $N=1,...,10$. \emph{Top left:} $k=10^{-1}$. \emph{Top right:} $k=10^{-2}$. \emph{Bottom left:} $k=10^{-3}$. \emph{Bottom right:} $k=10^{-4}$.}
		\label{fig:FluxConvergence}
	}
\end{figure}

In Figure~\ref{fig:FluxConvergence} we compare the net flux entering the \VARUN{sand-shale} domain, calculated using the mixed and the direct formulation of equations, as a function of polynomial degree for different values of $k$ in the shale blocks.
The data for these figures is given in Table~\ref{tab:FluxConvergence}. Note that the direct formulation converges from above towards the correct inflow flux, whereas the mixed formulation converges from below.

\begin{table}[!htb]
\centering
\begin{tabular}{c c c c c c c c c}
\hline
N	&\multicolumn{2}{c}{$k = 10^{-1}$} 	&\multicolumn{2}{c}{$k = 10^{-2}$} 	&\multicolumn{2}{c}{$k = 10^{-3}$}	 &\multicolumn{2}{c}{$k = 10^{-4}$}\\
	& Mixed	  & Direct	    	&Mixed 	& Direct	 & Mixed	 & Direct & Mixed	 & Direct \\
\hline
1	&0.63805   & 0.74149 	&0.51384  & 0.69273   &0.49296  & 0.68699   &0.49066  & 0.68641\\
2	&0.66541   & 0.69316 	&0.54101  & 0.62399   &0.51573  & 0.61572   &0.51279  & 0.61488\\
3	&0.67131   & 0.68423 	&0.54906  & 0.60794   &0.52121  & 0.59856   &0.51782  & 0.59760\\
4	&0.67339   & 0.68139 	&0.55208  & 0.60113   &0.52272  & 0.59099   &0.51904  & 0.58995\\
5	&0.67450   & 0.68003  	&0.55436  & 0.59711   &0.52371  & 0.58639   &0.51975  & 0.58528\\
6	&0.67512   & 0.67926 	&0.55568  & 0.59439   &0.52417  & 0.58320   &0.52003  & 0.58203\\
7	&0.67555   & 0.67877 	&0.55690  & 0.59239   &0.52459  & 0.58079   &0.52026  & 0.57958\\
8	&0.67582   & 0.67844 	&0.55772  & 0.59085   &0.52483  & 0.57890   &0.52036  & 0.57765\\
9	&0.67604   & 0.67821 	&0.55852  & 0.58960   &0.52508  & 0.57734   &0.52046  & 0.57605\\
10	&0.67619   & 0.67803 	&0.55910  & 0.58857   &0.52524  & 0.57603   &0.52051  & 0.57471\\
\hline
\end{tabular}
\caption{Data of net flux through the left boundary of the \VARUN{sand-shale} domain using mixed formulation and direct formulation for $20\times20$ elements, $N=1,...,10$, $k=10^{-1}$ (top-left), $10^{-2}$ (top-right), $10^{-3}$ (bottom-left) and $10^{-4}$ (bottom-right).}
\label{tab:FluxConvergence}
\end{table}
\clearpage
\subsection{The Impermeable-Streak system} \label{sec:test3}
The next example is from \cite{Durlofsky1993,HymanShashkovSteinberg97,Kikinzon2017}. The physical domain is a 2D unit square, $\Omega =[0,1]^2$.
\begin{figure}[hbt!]
	\centering{
		\subfloat{
			\begin{minipage}[t]{0.48\textwidth}
				\centering
				 \includegraphics[width=1\linewidth]{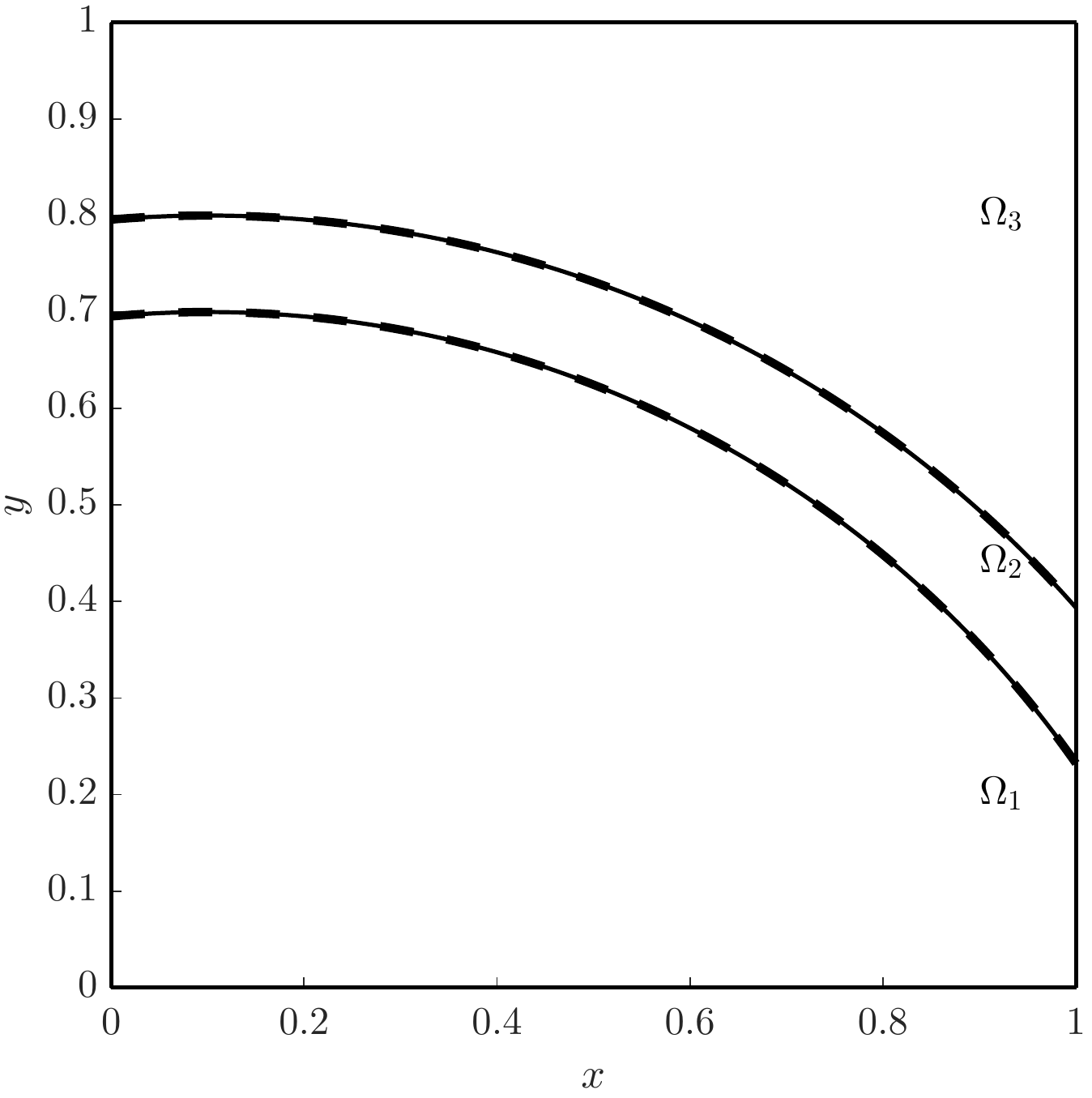}
			\end{minipage}
		}
		\subfloat{
			\begin{minipage}[t]{0.48\textwidth}
				\centering
				 \includegraphics[width=1\linewidth]{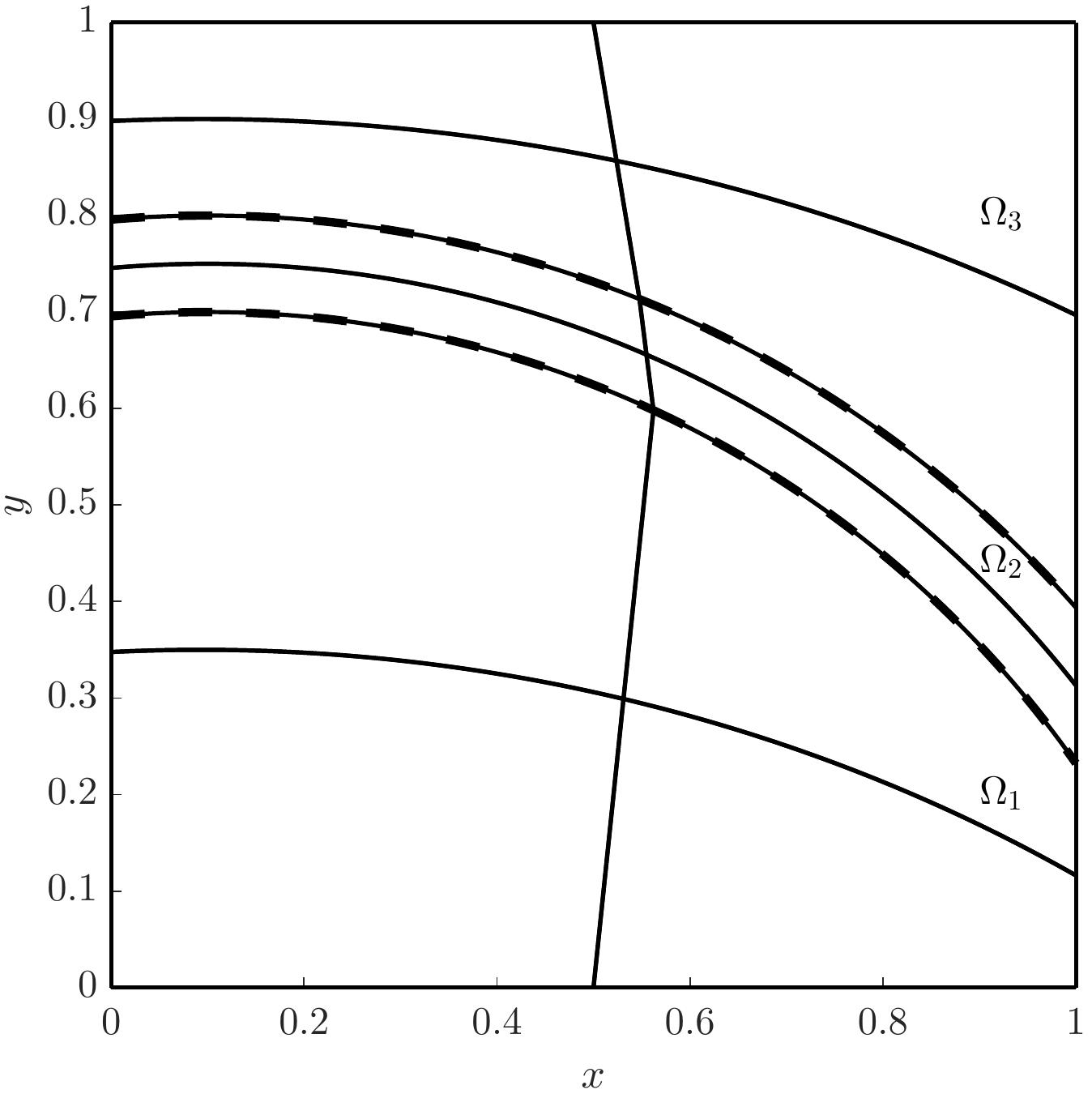}
			\end{minipage}
		}
		\caption{Three regions of the domain for the impermeable streak test case. The regions are separated by the dashed lines. The solid lines indicate the element boundaries. \emph{Left:} $1\times1$ element in each region. \emph{Right:} $2\times2$ elements in each region.}
		\label{fig:StreakGeometry}
	}
\end{figure}
The domain is divided into three different regions, $\Omega _1$, $\Omega _2$, and $\Omega _3 $, as shown in Figure ~\ref{fig:StreakGeometry} (left).
For calculations, each region is further divided into $K \times K $ elements.
Therefore, the total number of elements in the domain is given by $K \times K \times 3$.
In Figure \ref{fig:StreakGeometry} (right) we show the domain with each region divided into $2 \times 2$ elements.

\VARUN{The mixed formulation (\ref{eq:finite_element_mixed_formulation}) is solved, with $f=0$ and mixed boundary conditions, such that at the top and the bottom boundaries the net flux $\boldsymbol{u} \cdot \boldsymbol{n} = 0$, and at the left and the right boundaries, $p=1$ and $p=0$, respectively.} Permeability in $\Omega _1$ and $\Omega _3$ is given by $\mathbb{K}=\mathbb{I}$.
$\Omega_ 2$ has a low permeability and defined such that the component parallel to the local streak orientation is $k_\| = 10^{-1}$, and the component perpendicular to the local streak orientation is $k_\bot = 10^{-3}$.
The analytical expression for the permeability in terms of Cartesian coordinates is given in \cite{HymanShashkovSteinberg97} as,
\begin{equation*}
\begin{aligned}
K_{xx} & = \frac{k_\| (y+0.4)^2 + k_\bot (x-0.1)^2}{(x-0.1)^2+(y+0.4)^2}, \\
K_{xy} & = \frac{-(k_\| - k_\bot)(x-0.1)(y+0.4)}{(x-0.1)^2+(y+0.4)^2} ,\\
K_{yy} & = \frac{k_\| (x-0.1)^2 + k_\bot (y+0.4)^2}{(x-0.1)^2+(y+0.4)^2} .
\end{aligned}
\end{equation*}

\begin{figure}[!htb]
\centering
\includegraphics[scale=0.4]{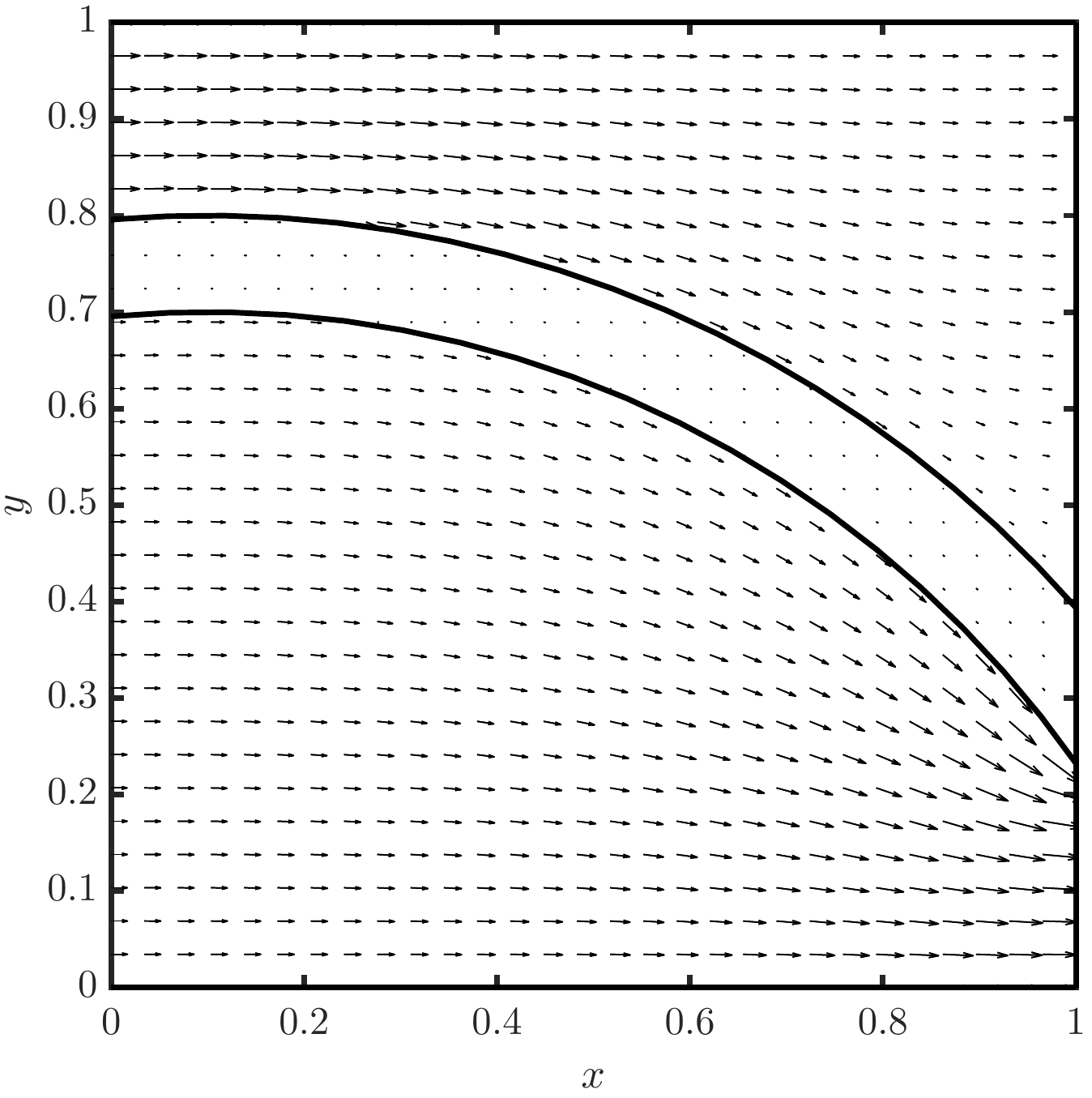}
\caption{Velocity vectors through the domain of permeability streak test case for $ 12\times12 $ elements, $ N=15 $.}
\label{fig:PSQuiver}
\end{figure}

\begin{figure}[!htb]
\centering
\includegraphics[scale=0.4]{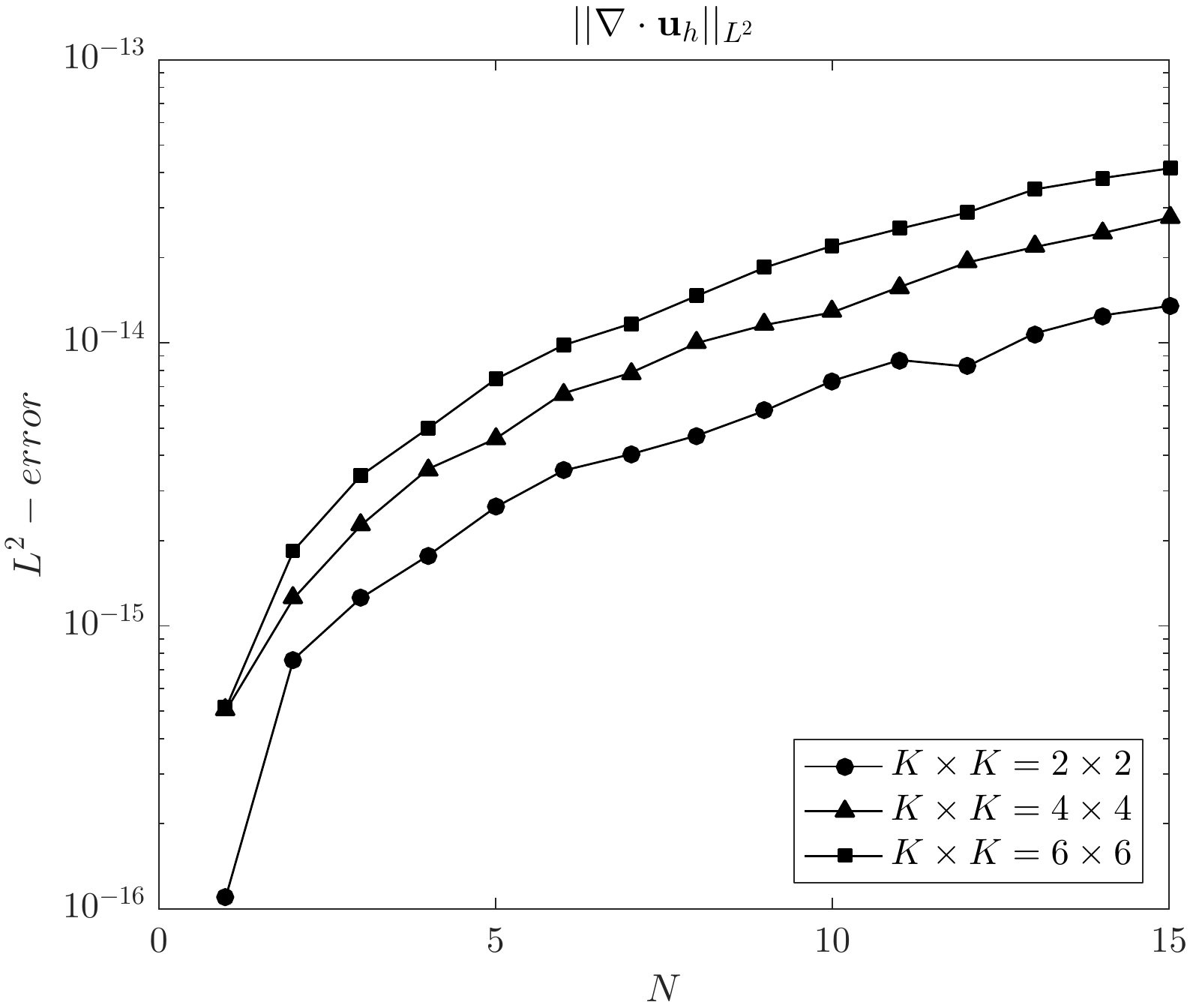}
\caption{The $L^2$-norm of \VARUN{$\nabla \cdot \boldsymbol{u} _h$} for $K\times K$ elements, $ K=2,4,6 $, $N=1,...,15$.}
\label{fig:StreakDiv}
\end{figure}

The flow field in the domain is shown in Figure~\ref{fig:PSQuiver}.
The magnitude of velocity in $\Omega _2$ is small due to low values of the permeability tensor in this region.
The velocity vectors bend in the direction of the permeability streak $\Omega _2$.
The $L^2$-norm of $\nabla \cdot \boldsymbol{u}$ over the entire domain as a function of polynomial degree, $N$, is shown in Figure~\ref{fig:StreakDiv}.
We can see that the flow field is divergence free up to machine precision \VARUN{because $f = 0$}.

The net flux through the system for varying number of elements and polynomial degree is given in Table~\ref{tab:StreakFlux}.
In this work the finest resolution corresponds to \VARUN{$12 \times 12 \times 3$} elements and $N=15$.
For this case the net influx at the left boundary is  $0.75668$.
The net influx and outflux from the region $\Omega _1$, $\Omega _2$ and $\Omega _3$ is given in Tables ~\ref{tab:FluxDomain1}, ~\ref{tab:FluxDomain2}, and ~\ref{tab:FluxDomain3}\VARUN{,} respectively.
The net influx for $\Omega _1$ is larger than the net outflux.
And the net outflux for $\Omega _2$ and $\Omega _3$ is larger than the net influx.
\begin{table}[!htb]
\centering
\begin{tabular}{p{0.6cm} c c c c c}
\hline
N & \multicolumn{5}{c}{Elements Division ($K \times K$)} \\
  			& $4 \times 4$ 	& $6 \times 6$ 	& $8 \times 8$ 	& $10 \times 10$ & $12 \times 12$\\
\hline
1 			&0.74689		&0.74908		&0.75061		 &0.75169		&0.75247\\
2 			&0.75268		&0.75407		&0.75479		 &0.75522		&0.75550\\
3			&0.75479		&0.75548		&0.75582		 &0.75602		&0.75615\\
4			&0.75561		&0.75600		&0.75620		 &0.75631		&0.75639\\
5			&0.75600		&0.75625		&0.75638		 &0.75645		&0.75650\\
6			&0.75621		&0.75639		&0.75648		 &0.75653		&0.75657\\
7			&0.75635		&0.75648		&0.75654		 &0.75658		&0.75660\\
8			&0.75643		&0.75653		&0.75658		 &0.75661		&0.75663\\
9			&0.75649		&0.75657		&0.75661		 &0.75663		&0.75665\\
10			&0.75654		&0.75660		&0.75663		 &0.75665		&0.75666\\
11			&0.75657		&0.75662		&0.75664		 &0.75666		&0.75667\\
12			&0.75659		&0.75663		&0.75665		 &0.75666		&0.75667\\
13			&0.75661		&0.75664		&0.75666		 &0.75667		&0.75668\\
14			&0.75662		&0.75665		&0.75667		 &0.75668		&0.75668\\
15			&0.75663		&0.75666		&0.75667		 &0.75668		&0.75668\\
\hline
\end{tabular}
\caption{Net flux through the left boundary of the permeability streak test case domain for $K\times K$ elements, $ K=4,6,8,10,12 $ and $N=1,...,15$.}
\label{tab:StreakFlux}
\end{table}

\begin{table}[!htb]
\centering
\begin{tabular}{p{0.6cm} c c c c c c c c c c}
	\hline
N & \multicolumn{10}{c}{Elements Division ($K \times K$)} \\
	& \multicolumn{2}{c}{$4 \times 4$} & \multicolumn{2}{c}{$6 \times 6$} & \multicolumn{2}{c}{$8 \times 8$} & \multicolumn{2}{c}{$10 \times 10$} & \multicolumn{2}{c}{$12 \times 12$}\\
			& In flux 	& Out flux & In flux 	& Out flux  & In flux 	& Out flux  & In flux 	& Out flux  & In flux 	& Out flux \\
\hline
1    &0.47155    &0.46497    &0.47342    &0.46657    &0.47483    &0.46786    &0.47585    &0.46881    &0.47659    &0.46952\\
2    &0.47674    &0.46965    &0.47812    &0.47098    &0.47883    &0.47168    &0.47926    &0.47211    &0.47954    &0.47239\\
3    &0.47883    &0.47168    &0.47952    &0.47236    &0.47986    &0.47270    &0.48005    &0.47291    &0.48018    &0.47304\\
4    &0.47964    &0.47249    &0.48004    &0.47289    &0.48023    &0.47309    &0.48035    &0.47321    &0.48042    &0.47329\\
5    &0.48003    &0.47288    &0.48029    &0.47315    &0.48041    &0.47328    &0.48049    &0.47336    &0.48053    &0.47341\\
6    &0.48025    &0.47311    &0.48043    &0.47330    &0.48051    &0.47339    &0.48056    &0.47345    &0.48060    &0.47348\\
7    &0.48038    &0.47325    &0.48051    &0.47339    &0.48057    &0.47346    &0.48061    &0.47350    &0.48063    &0.47353\\
8    &0.48047    &0.47334    &0.48057    &0.47345    &0.48061    &0.47350    &0.48064    &0.47354    &0.48066    &0.47356\\
9    &0.48053    &0.47340    &0.48060    &0.47349    &0.48064    &0.47354    &0.48066    &0.47356    &0.48068    &0.47358\\
10   &0.48057    &0.47345    &0.48063    &0.47352    &0.48066    &0.47356    &0.48068    &0.47358    &0.48069    &0.47360\\
11   &0.48060    &0.47349    &0.48065    &0.47355    &0.48067    &0.47358    &0.48069    &0.47360    &0.48070    &0.47361\\
12   &0.48062    &0.47352    &0.48066    &0.47357    &0.48068    &0.47359    &0.48070    &0.47361    &0.48070    &0.47362\\
13   &0.48064    &0.47354    &0.48067    &0.47358    &0.48069    &0.47360    &0.48070    &0.47362    &0.48071    &0.47363\\
14   &0.48065    &0.47355    &0.48068    &0.47359    &0.48070    &0.47361    &0.48071    &0.47362    &0.48071    &0.47363\\
15   &0.48067    &0.47357    &0.48069    &0.47360    &0.48070    &0.47362    &0.48071    &0.47363    &0.48071    &0.47364\\
\hline
\end{tabular}
\caption{Net flux through the left boundary of the region $ \Omega_{1} $ for $K\times K$ elements, $ K=4,6,8,10,12 $ and $N=1,...,15$.}
\label{tab:FluxDomain1}
\end{table}

\begin{table}[!htb]
\centering
\begin{tabular}{p{0.6cm} c c c c c c c c c c}
	\hline
N & \multicolumn{10}{c}{Elements Division ($K \times K$)} \\
	& \multicolumn{2}{c}{$4 \times 4$} & \multicolumn{2}{c}{$6 \times 6$} & \multicolumn{2}{c}{$8 \times 8$} & \multicolumn{2}{c}{$10 \times 10$} & \multicolumn{2}{c}{$12 \times 12$}\\
			& In flux 	& Out flux & In flux 	& Out flux  & In flux 	& Out flux  & In flux 	& Out flux  & In flux 	& Out flux \\
	\hline
1	&0.00930   &0.01080   &0.00931   &0.01106   &0.00931   &0.01119   &0.00932   &0.01130   &0.00932   &0.01130\\
2   &0.00932   &0.01132   &0.00933   &0.01138   &0.00933   &0.01139   &0.00933   &0.01140   &0.00933   &0.01140\\
3   &0.00933   &0.01139   &0.00933   &0.01140   &0.00933   &0.01140   &0.00933   &0.01140   &0.00933   &0.01139\\
4   &0.00933   &0.01140   &0.00933   &0.01140   &0.00934   &0.01140   &0.00934   &0.01139   &0.00934   &0.01139\\
5   &0.00933   &0.01140   &0.00934   &0.01139   &0.00934   &0.01139   &0.00934   &0.01138   &0.00934   &0.01138\\
6   &0.00934   &0.01139   &0.00934   &0.01139   &0.00934   &0.01138   &0.00934   &0.01137   &0.00934   &0.01137\\
7   &0.00934   &0.01139   &0.00934   &0.01138   &0.00934   &0.01137   &0.00934   &0.01137   &0.00934   &0.01136\\
8   &0.00934   &0.01138   &0.00934   &0.01137   &0.00934   &0.01137   &0.00934   &0.01136   &0.00934   &0.01136\\
9   &0.00934   &0.01138   &0.00934   &0.01137   &0.00934   &0.01136   &0.00934   &0.01136   &0.00934   &0.01135\\
10  &0.00934   &0.01137   &0.00934   &0.01136   &0.00934   &0.01136   &0.00934   &0.01135   &0.00934   &0.01135\\
11  &0.00934   &0.01137   &0.00934   &0.01136   &0.00934   &0.01135   &0.00934   &0.01135   &0.00934   &0.01135\\
12  &0.00934   &0.01136   &0.00934   &0.01136   &0.00934   &0.01135   &0.00934   &0.01135   &0.00934   &0.01134\\
13  &0.00934   &0.01136   &0.00934   &0.01135   &0.00934   &0.01135   &0.00934   &0.01134   &0.00934   &0.01134\\
14  &0.00934   &0.01136   &0.00934   &0.01135   &0.00934   &0.01134   &0.00934   &0.01134   &0.00934   &0.01134\\
15  &0.00934   &0.01135   &0.00934   &0.01135   &0.00934   &0.01134   &0.00934   &0.01134   &0.00934   &0.01134\\
\hline
\end{tabular}
\caption{Net flux through the left boundary of the region $ \Omega_{2} $ for $K\times K$ elements, $ K=4,6,8,10,12 $ and $N=1,...,15$.}
\label{tab:FluxDomain2}
\end{table}

\begin{table}[!htb]
\centering
\begin{tabular}{p{0.6cm} c c c c c c c c c c}
	\hline
N & \multicolumn{10}{c}{Elements Division ($K \times K$)} \\
	& \multicolumn{2}{c}{$4 \times 4$} & \multicolumn{2}{c}{$6 \times 6$} & \multicolumn{2}{c}{$8 \times 8$} & \multicolumn{2}{c}{$10 \times 10$} & \multicolumn{2}{c}{$12 \times 12$}\\
			& In flux 	& Out flux & In flux 	& Out flux  & In flux 	& Out flux  & In flux 	& Out flux  & In flux 	& Out flux \\
	\hline
1	&0.26604    &0.27112    &0.26636    &0.27146    &0.26647    &0.27157    &0.26653    &0.27163    &0.26656    &0.27166\\
2   &0.26662    &0.27172    &0.26663    &0.27172    &0.26663    &0.27172    &0.26663    &0.27172    &0.26663    &0.27172\\
3   &0.26663    &0.27172    &0.26663    &0.27172    &0.26664    &0.27172    &0.26664    &0.27172    &0.26664    &0.27172\\
4-15   &0.26664    &0.27172    &0.26664    &0.27172    &0.26664    &0.27172    &0.26664    &0.27172    &0.26664    &0.27172\\
\hline
\end{tabular}
\caption{Net flux through the left boundary of the region $ \Omega_{3} $ for $K\times K$ elements, $ K=4,6,8,10,12 $ and $N=1,...,15$.}
\label{tab:FluxDomain3}
\end{table}

\section{Future Work}
\VARUN{In the above sections, mixed and direct formulations of mimetic spectral element method are discussed. 
The next step is to explore this framework in the direction of hybrid formulations \cite{brezzi1991mixed,BoffiBrezziFortin2013,Cockburn2016}.
Additionally, the focus will be on developing multiscale methods \cite{Wang2014}, using these formulations, for reservoir modelling applications.}

%
%
%

\bibliographystyle{siam}
\bibliography{./library_darcy.bib}

\begin{thebibliography}{100}

\bibitem{Aarnes2008}
{\sc J.~E. Aarnes, S.~Krogstad, and K.-A. Lie}, {\em {Multiscale mixed/mimetic
  methods on corner-point grids}}, Computational Geosciences, 12 (2008),
  pp.~297--315.

\bibitem{Aavatsmark2002}
{\sc I.~Aavatsmark}, {\em {An Introduction to multipoint flux approximations
  for quadrilateral grids}}, Computational Geosciences, 6 (2002), pp.~405--432.

\bibitem{Aavatsmark2007}
\leavevmode\vrule height 2pt depth -1.6pt width 23pt, {\em {Interpretation of a
  two-point flux stencil for skew parallelogram grids}}, Computational
  Geosciences, 11 (2007), pp.~199--206.

\bibitem{Aavatsmark1998}
{\sc I.~Aavatsmark, T.~Barkve, O.~B{\o}e, and T.~Mannseth}, {\em
  {Discretization on unstructured grids for inhomogeneous, anisotropic media.
  Part I: derivation of the methods}}, SIAM Journal on Scientific Computing, 19
  (1998), pp.~1700--1716.

\bibitem{Aavatsmark1998a}
\leavevmode\vrule height 2pt depth -1.6pt width 23pt, {\em {Discretization on
  unstructured grids for inhomogeneous, anisotropic media. Part II: discussion
  and numerical results}}, SIAM Journal on Scientific Computing, 19 (1998),
  pp.~1717--1736.

\bibitem{Alpak2007}
{\sc F.~O. Alpak}, {\em {A mimetic finite volume discretization operator for
  reservoir simulation}}, in SPE Reservoir Simulation Symposium, Society of
  Petroleum Engineers, 2007.

\bibitem{Alpak2010a}
{\sc F.~O. Alpak}, {\em {A mimetic finite volume discretization method for
  reservoir simulation}}, SPE Journal, 15 (2010), pp.~436--453.

\bibitem{arnold:Quads}
{\sc D.~N. Arnold, D.~Boffi, and R.~S. Falk}, {\em {Quadrilateral H(div) Finite
  Elements}}, SIAM Journal on Numerical Analysis, 42 (2005), pp.~2429--2451.

\bibitem{arnold2006finite}
{\sc D.~N. Arnold, R.~S. Falk, and R.~Winther}, {\em {Finite element exterior
  calculus, homological techniques, and applications}}, Acta Numerica, 15
  (2006), pp.~1--155.

\bibitem{arnold2010finite}
\leavevmode\vrule height 2pt depth -1.6pt width 23pt, {\em {Finite element
  exterior calculus: from Hodge theory to numerical stability}}, Bulletin of
  the American Mathematical Society, 47 (2010), pp.~281--354.

\bibitem{Aziz1993}
{\sc K.~Aziz}, {\em {Reservoir simulation grids: opportunities and problems}},
  Journal of Petroleum Technology, 45 (1993), pp.~658--663.

\bibitem{Babuska1992a}
{\sc I.~Babuska and M.~Suri}, {\em {On locking and robustness in the finite
  element method}}, SIAM Journal on Numerical Analysis, 29 (1992),
  pp.~1261--1293.

\bibitem{Bastian2011a}
{\sc P.~Bastian, O.~Ippisch, and S.~Marnach}, {\em {Benchmark 3D: A Mimetic
  Finite Difference Method}}, in Finite Volumes for Complex Applications VI:
  Problems and Perspectives, Springer, Berlin, Heidelberg, 2011, pp.~961--968.

\bibitem{Bauer2017}
{\sc W.~Bauer and F.~Gay-Balmaz}, {\em {Variational integrators for anelastic
  and pseudo-incompressible flows}},  (2017).

\bibitem{Incropera2011}
{\sc T.~L. Bergman and F.~P. Incropera}, {\em {Fundamentals of heat and mass
  transfer.}}, Wiley, 2011.

\bibitem{bochev_mimetic_ls_2014}
{\sc P.~B. Bochev and M.~Gerritsma}, {\em {A spectral mimetic least-squares
  method}}, Computers and Mathematics with Applications, 68 (2014),
  pp.~1480--1502.

\bibitem{BochevGunzbuger}
{\sc P.~B. Bochev and M.~D. Gunzburger}, {\em Least-Squares Finite Element
  Methods}, Springer series in Applied Mathematical Sciences, Springer-Verlag,
  2009.

\bibitem{bochev2006principles}
{\sc P.~B. Bochev and J.~M. Hyman}, {\em {Principles of mimetic discretizations
  of differential operators}}, IMA Volumes In Mathematics and its Applications,
  142 (2006), p.~89.

\bibitem{bochev:RehabQuad}
{\sc P.~B. Bochev and D.~Ridzal}, {\em {Rehabilitation of the Lowest-Order
  Raviart--Thomas Element on Quadrilateral Grids}}, SIAM Journal on Numerical
  Analysis, 47 (2008), pp.~487--507.

\bibitem{BoffiBrezziFortin2013}
{\sc D.~Boffi, F.~Brezzi, and M.~Fortin}, {\em Mixed Finite Element Methods and
  Applications}, Springer Series in Computational Mathematics, Springer-Verlag,
  2013.

\bibitem{BoffiGastaldi2009}
{\sc D.~Boffi and L.~Gastaldi}, {\em {Some remarks on quadrilateral mixed
  finite elements}}, Computer and Structures, 87 (2009), pp.~751--757.

\bibitem{Bonelle2016}
{\sc J.~Bonelle, E.~Burman, P.~Cantin, and A.~Ern}, {\em {A vertex-based scheme
  on polyhedral meshes for advection-reaction equations with sub-mesh
  stabilization}}, hal-01285957,  (2016), pp.~1--20.

\bibitem{Bonelle2015}
{\sc J.~Bonelle, D.~A. {Di Pietro}, and A.~Ern}, {\em {Low-order reconstruction
  operators on polyhedral meshes: application to compatible discrete operator
  schemes}}, Computer Aided Geometric Design, 35-36 (2015), pp.~27--41.

\bibitem{Bonelle2014}
{\sc J.~Bonelle and A.~Ern}, {\em {Analysis of Compatible Discrete Operator
  schemes for elliptic problems on polyhedral meshes}}, ESAIM: Mathematical
  Modelling and Numerical Analysis, 48 (2014), pp.~553--581.

\bibitem{Bonelle2015a}
\leavevmode\vrule height 2pt depth -1.6pt width 23pt, {\em {Analysis of
  compatible discrete operator schemes for the Stokes equations on polyhedral
  meshes}}, IMA Journal of Numerical Analysis, 35 (2015), pp.~1672--1697.

\bibitem{bossavit_japan_computational_1}
{\sc A.~Bossavit}, {\em {Computational electromagnetism and geometry: (1)
  Network equations}}, Journal of the Japan Society of Applied
  Electromagnetics, 7 (1999), pp.~150--159.

\bibitem{bossavit_japan_computational_2}
\leavevmode\vrule height 2pt depth -1.6pt width 23pt, {\em {Computational
  electromagnetism and geometry: (2) Network constitutive laws}}, Journal of
  the Japan Society of Applied Electromagnetics, 7 (1999), pp.~294--301.

\bibitem{bossavit_japan_computational_3}
\leavevmode\vrule height 2pt depth -1.6pt width 23pt, {\em {Computational
  electromagnetism and geometry: (3) Convergence}}, Journal of the Japan
  Society of Applied Electromagnetics, 7 (1999), pp.~401--408.

\bibitem{bossavit_japan_computational_4}
\leavevmode\vrule height 2pt depth -1.6pt width 23pt, {\em {Computational
  electromagnetism and geometry: (4) From degrees of freedom to fields}},
  Journal of the Japan Society of Applied Electromagnetics, 8 (2000),
  pp.~102--109.

\bibitem{bossavit_japan_computational_5}
\leavevmode\vrule height 2pt depth -1.6pt width 23pt, {\em {Computational
  electromagnetism and geometry: (5) The ``Galerkin Hodge''}}, Journal of the
  Japan Society of Applied Electromagnetics, 8 (2000), pp.~203--209.

\bibitem{bouman::icosahom2009}
{\sc M.~Bouman, A.~Palha, J.~Kreeft, and M.~Gerritsma}, {\em {A Conservative
  Spectral Element Method for Curvilinear Domains}}, in Spectral and High Order
  Methods for Partial Differential Equations, vol.~76 of Lecture Notes in
  Computational Science and Engineering, Springer, 2011, pp.~111--119.

\bibitem{brezzi2010}
{\sc F.~Brezzi and A.~Buffa}, {\em {Innovative mimetic discretizations for
  electromagnetic problems}}, Journal of Computational and Applied Mathematics,
  234 (2010), pp.~1980--1987.

\bibitem{BrezziBuffaLipnikov2009}
{\sc F.~Brezzi, A.~Buffa, and K.~Lipnikov}, {\em {Mimetic finite differences
  for elliptic problems}}, Mathematical Modelling and Numerical Analysis, 43
  (2009), pp.~277--296.

\bibitem{Brezzi2014}
{\sc F.~Brezzi, R.~S. Falk, and L.~{Donatella Marini}}, {\em {Basic principles
  of mixed Virtual Element Methods}}, ESAIM: Mathematical Modelling and
  Numerical Analysis, 48 (2014), pp.~1227--1240.

\bibitem{brezzi1991mixed}
{\sc F.~Brezzi and M.~Fortin}, {\em {Mixed and Hybrid Finite Element Methods}},
  vol.~15 of Springer Series in Computational Mathematics, Springer, 1991.

\bibitem{Brezzi2006}
{\sc F.~Brezzi, K.~Lipnikov, and M.~Shashkov}, {\em {Convergence of mimetic
  finite difference method for diffusion problems on polyhedral meshes with
  curved faces}}, Mathematical Models and Methods in Applied Sciences, 16
  (2006), pp.~275--297.

\bibitem{Brezzi2007}
{\sc F.~Brezzi, K.~Lipnikov, M.~Shashkov, and V.~Simoncini}, {\em {A new
  discretization methodology for diffusion problems on generalized polyhedral
  meshes}}, Computer Methods in Applied Mechanics and Engineering, 196 (2007),
  pp.~3682--3692.

\bibitem{Brezzi2005}
{\sc F.~Brezzi, K.~Lipnikov, and V.~Simoncini}, {\em {A family of mimetic
  finite difference methods on polygonal and polyhedral meshes}}, Mathematical
  Models and Methods in Applied Sciences, 15 (2005), pp.~1533--1551.

\bibitem{Budd2003}
{\sc C.~Budd and M.~Piggott}, {\em {Geometric Integration and its
  Applications}}, in Handbook of Numerical Analysis, vol.~11, North-Holland,
  2003, pp.~35--139.

\bibitem{BuffaDeFalcoSangalli2011}
{\sc A.~Buffa, D.~de~Falco, and G.~Sangalli}, {\em {Isogeometric analysis: new
  stable elements for the Stokes equation}}, International Journal for
  Numerical Methods in Fluids,  (2011).

\bibitem{Canuto_et_al}
{\sc C.~Canuto, M.~Y. Hussaini, A.~Quarteroni, and T.~A. Zang}, {\em {Spectral
  Methods in Fluid Dynamics}}, Springer Verlag, 1988.

\bibitem{Christiansen2011}
{\sc S.~H. Christiansen, H.~Z. Munthe-Kaas, and B.~Owren}, {\em {Topics in
  structure-preserving discretization}}, Acta Numerica, 20 (2011), pp.~1--119.

\bibitem{Cockburn2016}
{\sc B.~Cockburn}, {\em Static Condensation, Hybridization, and the Devising of
  the HDG Methods}, Springer International Publishing, 2016, pp.~129--177.

\bibitem{BeiraodaVeiga2014}
{\sc L.~B. da~Veiga, F.~Brezzi, L.~D. Marini, and A.~Russo}, {\em {The
  Hitchhiker's Guide to the Virtual Element Method}}, Mathematical Models and
  Methods in Applied Sciences, 24 (2014), pp.~1541--1573.

\bibitem{DaVeiga2015}
{\sc L.~B. da~Veiga, F.~Brezzi, L.~D. Marini, and A.~Russo}, {\em
  {$H(\mathrm{div})$ and $H(\mathrm{curl})$ conforming virtual element
  methods}}, Numerische Mathematik,  (2015), pp.~1--30.

\bibitem{BeiraodaVeiga2016}
{\sc L.~B. da~Veiga, F.~Brezzi, L.~D. Marini, and A.~Russo}, {\em {Virtual
  Element Method for general second-order elliptic problems on polygonal
  meshes}}, Mathematical Models and Methods in Applied Sciences, 26 (2016),
  pp.~729--750.

\bibitem{Veiga_mimetic_finite_differences}
{\sc L.~B. da~Veiga, K.~Lipnikov, and G.~Manzini}, {\em {The Mimetic Finite
  Difference Method for Elliptic Problems}}, Springer International Publishing,
  2014.

\bibitem{DaVeiga2015a}
{\sc L.~B. da~Veiga, C.~Lovadina, and G.~Vacca}, {\em {Divergence free Virtual
  Elements for the Stokes problem on polygonal meshes}}, arXiv:1510.01655v1,
  (2015).

\bibitem{desbrun2005discrete}
{\sc M.~Desbrun, A.~N. Hirani, M.~Leok, and J.~E. Marsden}, {\em {Discrete
  exterior calculus}}, arXiv:math/0508341v2,  (2005).

\bibitem{DiPietro2015}
{\sc D.~A. {Di Pietro} and A.~Ern}, {\em {Hybrid high-order methods for
  variable-diffusion problems on general meshes}}, Comptes Rendus Mathematique,
  353 (2015), pp.~31--34.

\bibitem{DiPietro2014}
{\sc D.~A. {Di Pietro}, A.~Ern, and S.~Lemaire}, {\em {An arbitrary-order and
  compact-stencil discretization of diffusion on general meshes based on local
  reconstruction operators}}, Computational Methods in Applied Mathematics, 14
  (2014).

\bibitem{Dodziuk76}
{\sc J.~Dodziuk}, {\em {{F}inite difference approach to the {H}odge theory of
  harmonic functions}}, American Journal of Mathematics, 98 (1976),
  pp.~79--104.

\bibitem{Durlofsky1993}
{\sc L.~J. Durlofsky}, {\em {A Triangle Based Mixed Finite Element — Finite
  Volume Technique for Modeling Two Phase Flow through Porous Media}}, Journal
  of Computational Physics, 105 (1993), pp.~252--266.

\bibitem{Durlofsky1994}
\leavevmode\vrule height 2pt depth -1.6pt width 23pt, {\em {Accuracy of mixed
  and control volume finite element approximations to Darcy velocity and
  related quantities}}, Water Resources Research, 30 (1994), pp.~965--973.

\bibitem{Dziubek2016}
{\sc A.~Dziubek, G.~Guidoboni, A.~Harris, A.~N. Hirani, E.~Rusjan, and
  W.~Thistleton}, {\em {Effect of ocular shape and vascular geometry on retinal
  hemodynamics: a computational model}}, Biomechanics and Modeling in
  Mechanobiology, 15 (2016), pp.~893--907.

\bibitem{Edwards2002}
{\sc M.~G. Edwards}, {\em {Unstructured, control-volume distributed,
  full-tensor finit--volume schemes with flow based grids}}, Computational
  Geosciences, 6 (2002), pp.~433--452.

\bibitem{Edwards1998}
{\sc M.~G. Edwards and C.~F. Rogers}, {\em {Finite volume discretization with
  imposed flux continuity for the general tensor pressure equation}},
  Computational Geosciences, 2 (1998), pp.~259--290.

\bibitem{ElcottTongetal2007}
{\sc S.~Elcott, Y.~Tong, E.~Kanso, P.~Schr{\"{o}}der, and M.~Desbrun}, {\em
  {Stable, circulation-preserving, simplicial fluids}}, ACM TRansactions on
  Graphics, 26 (2007).

\bibitem{Evans2013}
{\sc J.~A. Evans and T.~J. Hughes}, {\em {Isogeometric divergence-conforming
  B-splines for the unsteady Navier-Stokes equations}}, Journal of
  Computational Physics, 241 (2013), pp.~141--167.

\bibitem{Forsyth1990}
{\sc P.~A. Forsyth}, {\em {A control-volume, finite-element method for local
  mesh refinement in thermal reservoir simulation}}, SPE Reservoir Engineering,
  5 (1990), pp.~561--566.

\bibitem{gerritsma::edge_basis}
{\sc M.~Gerritsma}, {\em {Edge Functions for Spectral Element Methods}}, in
  Spectral and High Order Methods for Partial Differential Equations, vol.~76
  of Lecture Notes in Computational Science and Engineering, Springer, 2011,
  pp.~199--207.

\bibitem{gerritsma_mimetic_ls_stokes_2016}
{\sc M.~Gerritsma and P.~B. Bochev}, {\em {A spectral mimetic least-squares
  method for the Stokes equations with no-slip boundary condition}}, Computers
  and Mathematics with Applications, 71 (2016), pp.~2285--2300.

\bibitem{Gerritsma}
{\sc M.~Gerritsma, M.~Bouman, and A.~Palha}, {\em {Least-Squares Spectral
  Element Method on a Staggered Grid}}, in Large-Scale Scientific Computing,
  vol.~5910 of Lecture Notes in Computer Science, Springer, 2010, pp.~653--661.

\bibitem{gerritsmaicosahom2012}
{\sc M.~Gerritsma, R.~Hiemstra, J.~Kreeft, A.~Palha, P.~P. Rebelo, and
  D.~Toshniwal}, {\em {The Geometric Basis of Numerical Methods}}, in Spectral
  and High Order Methods for Partial Differential Equations, vol.~95 of Lecture
  Notes in Computational Science and Engineering, Springer, 2013, pp.~17--35.

\bibitem{Gunasekera1997}
{\sc D.~Gunasekera, J.~Cox, and P.~Lindsey}, {\em {The Generation and
  Application of K-Orthogonal Grid Systems}}, in SPE Reservoir Simulation
  Symposium, Society of Petroleum Engineers, 1997, pp.~199--214.

\bibitem{Hairer2006}
{\sc E.~Hairer, C.~Lubich, and G.~Wanner}, {\em {Geometric Numerical
  Integration}}, Springer, 2006.

\bibitem{Heinemann1991}
{\sc Z.~E. Heinemann, C.~W. Brand, M.~Munka, and Y.~M. Chen}, {\em {Modeling
  reservoir geometry with irregular grids}}, SPE Reservoir Engineering, 6
  (1991), pp.~225--232.

\bibitem{Herbin2008}
{\sc R.~Herbin and F.~Hubert}, {\em {Benchmark on discretization schemes for
  anisotropic diffusion problems on general grids}}, in Finite Volumes for
  Complex Applications V: Problems and Perspectives, Wiley, 2008, pp.~659--692.

\bibitem{Hermeline2000}
{\sc F.~Hermeline}, {\em {A finite volume method for the approximation of
  diffusion operators on distorted meshes}}, Journal of Computational Physics,
  160 (2000), pp.~481--499.

\bibitem{Hiemstra2014}
{\sc R.~Hiemstra, D.~Toshniwal, R.~Huijsmans, and M.~Gerritsma}, {\em {High
  order geometric methods with exact conservation properties}}, Journal of
  Computational Physics, 257 (2014), pp.~1444--1471.

\bibitem{hiptmair2001}
{\sc R.~Hiptmair}, {\em {PIER}}, in Geometric Methods for Computational
  Electromagnetics, vol.~42, EMW Publishing, 2001, pp.~271--299.

\bibitem{Hirani_phd_2003}
{\sc A.~Hirani}, {\em {{D}iscrete {E}xterior {C}alculus}}, PhD thesis,
  California Institute of Technology, 2003.

\bibitem{Hirani2015}
{\sc A.~N. Hirani, K.~B. Nakshatrala, and J.~H. Chaudhry}, {\em {Numerical
  method for Darcy flow derived using discrete exterior calculus}},
  International Journal for Computational Methods in Engineering Science and
  Mechanics, 16 (2015), pp.~151--169.

\bibitem{HymanShashkovSteinberg2002}
{\sc J.~M. Hyman, J.~Morel, M.~Shashkov, and S.~Steinberg}, {\em {Mimetic
  finite difference methods for diffusion equations}}, Computational
  Geosciences, 6 (2002), pp.~333--352.

\bibitem{HymanScovel90}
{\sc J.~M. Hyman and J.~C. Scovel}, {\em {Deriving mimetic difference
  approximations to differential operators using algebraic topology}}, tech.
  rep., Los Alamos National Laboratory, 1990.

\bibitem{HymanShashkovSteinberg97}
{\sc J.~M. Hyman, M.~Shashkov, and S.~Steinberg}, {\em {The numerical solution
  of diffusion problems in strongly heterogeous non-isotropic materials}},
  Journal of Computational Physics, 132 (1997), pp.~130--148.

\bibitem{HYmanSteinberg2004}
{\sc J.~M. Hyman and S.~Steinberg}, {\em {The convergence of mimetic methods
  for rough grids}}, Computers and Mathematics with applications, 47 (2004),
  pp.~1565--1610.

\bibitem{Kikinzon2017}
{\sc E.~Kikinzon, Y.~Kuznetsov, K.~Lipnikov, and M.~Shashkov}, {\em
  {Approximate static condensation algorithm for solving multi-material
  diffusion problems on meshes non-aligned with material interfaces}}, Journal
  of Computational Physics,  (2017).

\bibitem{Kouranbaeva2000}
{\sc S.~Kouranbaeva and S.~Shkoller}, {\em {A variational approach to
  second-order multisymplectic field theory}}, Journal of Geometry and Physics,
  35 (2000), pp.~333--366.

\bibitem{Kraus2015}
{\sc M.~Kraus and O.~Maj}, {\em {Variational integrators for nonvariational
  partial differential equations}}, Physica D: Nonlinear Phenomena, 310 (2015),
  pp.~37--71.

\bibitem{kreeft::stokes}
{\sc J.~Kreeft and M.~Gerritsma}, {\em {Mixed mimetic spectral element method
  for Stokes flow: A pointwise divergence-free solution}}, Journal of
  Computational Physics, 240 (2013), pp.~284--309.

\bibitem{Kreeft2011}
{\sc J.~Kreeft, A.~Palha, and M.~Gerritsma}, {\em {Mimetic framework on
  curvilinear quadrilaterals of arbitrary order}}, arXiv:1111.4304,  (2011),
  p.~69.

\bibitem{Lie2012}
{\sc K.~Lie, S.~Krogstad, I.~S. Ligaarden, J.~R. Natvig, H.~M. Nilsen, and
  B.~Skaflestad}, {\em {Open-source MATLAB implementation of consistent
  discretisations on complex grids}}, Computational Geosciences, 16 (2012),
  pp.~297--322.

\bibitem{Manzini2007}
{\sc G.~Manzini and M.~Putti}, {\em {Mesh locking effects in the finite volume
  solution of 2-D anisotropic diffusion equations}}, Journal of Computational
  Physics, 220 (2007), pp.~751--771.

\bibitem{Marsden2003}
{\sc J.~E. Marsden and M.~West}, {\em {Discrete mechanics and variational
  integrators}}, Acta Numerica 2001, 10 (2003), pp.~357--514.

\bibitem{MullenCraneetal2009}
{\sc P.~Mullen, K.~Crane, D.~Pavlov, Y.~Tong, and M.~Desbrun}, {\em
  {Energy-preserving integrators for fluid animation}}, ACM Transactions on
  Graphics, 28 (2009).

\bibitem{Neuman1977}
{\sc S.~P. Neuman}, {\em {Theoretical derivation of Darcy's law}}, Acta
  Mechanica, 25 (1977), pp.~153--170.

\bibitem{Nicolaides}
{\sc R.~Nicolaides}, {\em {Discrete Discretization of planar div-curl
  problems}}, SIAM Journal on Numerical Analysis, 29 (1992), pp.~32--56.

\bibitem{Nilsen2012}
{\sc H.~M. Nilsen, J.~R. Natvig, and K.-A. Lie}, {\em {Accurate Modeling of
  Faults by Multipoint, Mimetic, and Mixed Methods}}, SPE Journal,  (2012),
  pp.~pp. 568--579.

\bibitem{Palagi1994}
{\sc C.~L. Palagi and K.~Aziz}, {\em {Use of Voronoi grid in reservoir
  simulation}}, SPE Advanced Technology Series, 2 (1994), pp.~69--77.

\bibitem{palha:lssc2009}
{\sc A.~Palha and M.~Gerritsma}, {\em {Mimetic Least-Squares Spectral/$hp$
  Finite Element Method for the Poisson Equation}}, in Large-Scale Scientific
  Computing, vol.~5910 of Lecture Notes in Computer Science, Springer, 2010,
  pp.~662--670.

\bibitem{palha::icosahom2009}
{\sc A.~Palha and M.~Gerritsma}, {\em {Spectral Element Approximation of the
  Hodge-$\star$ Operator in Curved Elements}}, in Spectral and High Order
  Methods for Partial Differential Equations, vol.~76 of Lecture Notes in
  Computational Science and Engineering, Springer, 2010, pp.~283--291.

\bibitem{Palha2017}
{\sc A.~Palha and M.~Gerritsma}, {\em {A mass, energy, enstrophy and vorticity
  conserving (MEEVC) mimetic spectral element discretization for the 2D
  incompressible Navier-Stokes equations}}, Journal of Computational Physics,
  328 (2017), pp.~200--220.

\bibitem{Palha2016}
{\sc A.~Palha, B.~Koren, and F.~Felici}, {\em {A mimetic spectral element
  solver for the Grad–Shafranov equation}}, Journal of Computational Physics,
  316 (2016), pp.~63--93.

\bibitem{palhaAdvectionIcosahom2014}
{\sc A.~Palha, P.~P. Rebelo, and M.~Gerritsma}, {\em {Mimetic Spectral Element
  Advection}}, in Spectral and High Order Methods for Partial Differential
  Equations - ICOSAHOM 2012, vol.~95 of Lecture Notes in Computational Science
  and Engineering, Springer International Publishing, 2014, pp.~325--335.

\bibitem{Palha2014}
{\sc A.~Palha, P.~P. Rebelo, R.~Hiemstra, J.~Kreeft, and M.~Gerritsma}, {\em
  {Physics-compatible discretization techniques on single and dual grids, with
  application to the Poisson equation of volume forms}}, Journal of
  Computational Physics, 257 (2014), pp.~1394--1422.

\bibitem{PavlovMullenetal2010}
{\sc D.~Pavlov, P.~Mullen, Y.~Tong, E.~Kanso, J.~E. Marsden, and M.~Desbrun},
  {\em {Structure preserving discretization of incompressible fluids}}, Physica
  D: Nonlinear phenomena, 240 (2011), pp.~443--458.

\bibitem{Perona1990}
{\sc P.~Perona and J.~Malik}, {\em {Scale-space and edge detection using
  anisotropic diffusion}}, IEEE Transactions on Pattern Analysis and Machine
  Intelligence, 12 (1990), pp.~629--639.

\bibitem{Perot2000}
{\sc J.~B. Perot}, {\em {Conservation properties of unstructured staggered mesh
  schemes}}, Journal of Computational Physics, 159 (2000), pp.~58--89.

\bibitem{perot43discrete}
\leavevmode\vrule height 2pt depth -1.6pt width 23pt, {\em {Discrete
  conservation properties of unstructured mesh schemes}}, Annual Review of
  Fluid Mechanics, 43 (2011), pp.~299--318.

\bibitem{PerotSubramanian2007}
{\sc J.~B. Perot and V.~Subramanian}, {\em {A discrete calculus analysis of the
  Keller Box scheme and a generalization of the method to arbitrary meshes}},
  Journal of Computational Physics, 226 (2007), pp.~494--508.

\bibitem{PerotSubramanian2007a}
\leavevmode\vrule height 2pt depth -1.6pt width 23pt, {\em {Discrete calculus
  methods for diffusion}}, Journal of Computational Physics, 224 (2007),
  pp.~59--81.

\bibitem{perot2006mimetic}
{\sc J.~B. Perot, D.~Vidovic, and P.~Wesseling}, {\em {Mimetic reconstruction
  of vectors}}, IMA Volumes in Mathematics and its Applications, 142 (2006),
  p.~173.

\bibitem{Rapetti2007}
{\sc F.~Rapetti}, {\em {High order edge elements on simplicial meshes}}, ESAIM:
  Mathematical Modelling and Numerical Analysis, 41 (2007), pp.~1001--1020.

\bibitem{Rapetti2009}
\leavevmode\vrule height 2pt depth -1.6pt width 23pt, {\em {Whitney forms of
  higher order}}, SIAM J. Numer. Anal., 47 (2009), pp.~2369--2386.

\bibitem{Rebelo2014}
{\sc P.~P. Rebelo, A.~Palha, and M.~Gerritsma}, {\em {Mixed mimetic spectral
  element method applied to Darcy's problem}}, in Spectral and High Order
  Methods for Partial Differential Equations - ICOSAHOM 2012, vol.~95 of
  Lecture Notes in Computational Science and Engineering, Springer, 2014,
  pp.~373--382.

\bibitem{RobidouxAdjointGradients1996}
{\sc N.~Robidoux}, {\em {A New Method of Construction of Adjoint Gradients and
  Divergences on Logically Rectangular Smooth Grids}}, in Finite Volumes for
  Complex Applications: Problems and Persepctives, {\'{E}}ditions Herm{\`{e}}s,
  Rouen, France, 1996, pp.~261--272.

\bibitem{RobidouxThesis}
\leavevmode\vrule height 2pt depth -1.6pt width 23pt, {\em {Numerical solution
  of the steady diffusion equation with discontinuous coefficients}}, PhD
  thesis, University of New Mexico, Albuquerque, NM, USA, 2002.

\bibitem{robidoux-polynomial}
{\sc N.~Robidoux}, {\em {Polynomial histopolation, superconvergent degrees of
  freedom, and pseudospectral discrete Hodge operators}}, Unpublished:
  http://people.math.sfu.ca/$\sim$nrobidou/public\_html/prints/histogram/histogram.pdf,
   (2008).

\bibitem{Robidoux2011}
{\sc N.~Robidoux and S.~Steinberg}, {\em {A discrete vector calculus in tensor
  grids}}, Computational Methods in Applied Mathematics, 11 (2011), pp.~23--66.

\bibitem{bookShashkov}
{\sc M.~Shashkov}, {\em {Conservative finite-difference methods on general
  grids}}, CRC Press, Boca Raton, FL, USA, 1996.

\bibitem{Sovinec2004}
{\sc C.~Sovinec, A.~Glasser, T.~Gianakon, D.~Barnes, R.~Nebel, S.~Kruger,
  D.~Schnack, S.~Plimpton, A.~Tarditi, M.~Chu, and N.~Team}, {\em {Nonlinear
  magnetohydrodynamics simulation using high-order finite elements}}, Journal
  of Computational Physics, 195 (2004), pp.~355--386.

\bibitem{Steinberg1996}
{\sc S.~Steinberg}, {\em {A discrete calculus with applications of higher-order
  discretizations to boundary-value problems}}, Computational Methods in
  Applied Mathematics, 42 (2004), pp.~228--261.

\bibitem{SteibergZingano2009}
{\sc S.~Steinberg and J.~P. Zingano}, {\em {Error estimates on arbitrary grids
  for 2nd-order mimetic discretization of Sturm-Liouville problems}},
  Computational Methods in Applied Mathematics, 9 (2009), pp.~192--202.

\bibitem{tarhasaari1999some}
{\sc T.~Tarhasaari, L.~Kettunen, and A.~Bossavit}, {\em {{S}ome realizations of
  a discrete {H}odge operator: a reinterpretation of finite element
  techniques}}, IEEE Transactions on Magnetics, 35 (1999).

\bibitem{Taylor1938}
{\sc G.~I. Taylor}, {\em {Production and dissipation of vorticity in a
  turbulent fluid}}, Proceedings of the Royal Society A: Mathematical, Physical
  and Engineering Sciences, 164 (1938), pp.~15--23.

\bibitem{tonti1975formal}
{\sc E.~Tonti}, {\em {On the formal structure of physical theories}}, tech.
  rep., Italian National Research Council, 1975.

\bibitem{Wang2014}
{\sc Y.~Wang, H.~Hajibeygi, and H.~A. Tchelepi}, {\em {Algebraic multiscale
  solver for flow in heterogeneous porous media}}, Journal of Computational
  Physics, 259 (2014), pp.~284--303.

\bibitem{Whitney57}
{\sc H.~Whitney}, {\em {Geometric integration theory}}, Dover Publications,
  Inc., 1957.

\bibitem{Wu2009}
{\sc X.-H. Wu and R.~Parashkevov}, {\em {Effect of grid deviation on flow
  solutions}}, SPE Journal, 14 (2009), pp.~67--77.

\bibitem{Younes2010}
{\sc A.~Younes, P.~Ackerer, and F.~Delay}, {\em {Mixed finite elements for
  solving 2-D diffusion-type equations}}, Reviews of Geophysics, 48 (2010),
  p.~RG1004.

\bibitem{Young1999}
{\sc L.~C. Young}, {\em {Rigorous treatment of distorted grids in 3D}}, in SPE
  Reservoir Simulation Symposium, Society of Petroleum Engineers, 1999.

\bibitem{ZhangSchmidtPerot2002}
{\sc X.~Zhang, Schmidt, D., and J.~B. Perot}, {\em {Accuracy and conservation
  properties of a theree-dimensional unstructured staggered mesh scheme for
  fluid dynamics}}, Journal of Computational Physics, 175 (2002), pp.~764--791.

\end{thebibliography}


\end{document}